\documentclass[12pt]{amsart}
\usepackage{amscd,amssymb,verbatim,mathrsfs,graphicx}
\usepackage[small,nohug,heads=vee]{diagrams}
\usepackage{tikz}
\usepackage[matrix,arrow]{xy}
\usepackage{url}
\usetikzlibrary{decorations.markings, shapes.geometric, patterns}

\newskip\stdskip                      
\newcommand{\dL}{{\Bbb L}}

\renewcommand{\mod}{\operatorname{mod}}

\newcommand{\und}{\underline}

\newcommand{\OO}{\mathcal{O}}

\newcommand{\coker}{\operatorname{coker}}

\newcommand{\KK}{\mathcal{K}}
\newcommand{\RR}{\mathcal{R}}
\newcommand{\bL}{{\bf L}}

\newcommand{\WW}{\mathcal{W}}

\newcommand{\BB}{\mathscr{B}}

\newcommand{\G}{{\Bbb G}}

\newcommand{\mg}{{\frak m}}
\newcommand{\hra}{\hookrightarrow}
\newcommand{\lan}{\langle}
\newcommand{\ran}{\rangle}
\newcommand{\Coh}{\operatorname{Coh}}

\newcommand{\CC}{\mathcal{C}}
\newcommand{\UU}{\mathcal{U}}

\newcommand{\Spec}{\operatorname{Spec}}
\newcommand{\Spf}{\operatorname{Spf}}

\newcommand{\Proj}{\operatorname{Proj}}
\renewcommand{\P}{{\Bbb P}}

\newcommand{\si}{\sigma}

\newcommand{\ga}{\gamma}
\newcommand{\de}{\delta}
\newcommand{\eps}{\epsilon}

\renewcommand{\ker}{\operatorname{ker}}

\numberwithin{equation}{subsection}
\newcommand{\Perf}{\operatorname{Perf}}

\newcommand{\m}{\mathfrak{m}}
\newcommand{\T}{\mathbb{T}}

\newcommand{\A}{{\Bbb A}}

\newtheorem{thm}{Theorem}[subsection]
\newtheorem{prop}[thm]{Proposition}
\newtheorem{lem}[thm]{Lemma}
\newtheorem{cor}[thm]{Corollary}
{  \theoremstyle{definition}
\newtheorem{defi}[thm]{Definition}

\newtheorem{rem}[thm]{Remark}
\newtheorem{rems}[thm]{Remarks}
}

\newcommand{\Pf}{\noindent {\it Proof}}

\newcommand{\Lie}{\operatorname{Lie}}
\newcommand{\ov}{\overline}

\renewcommand{\AA}{\mathscr{A}}
\newcommand{\FF}{\mathcal{F}}

\newcommand{\MM}{\mathcal{M}}
\newcommand{\TT}{\mathcal{T}}

\newcommand{\LL}{\mathcal{L}}

\newcommand{\Om}{\Omega}

\newcommand{\Hom}{\operatorname{Hom}}

\newcommand{\Ext}{\operatorname{Ext}}

\newcommand{\Res}{\operatorname{Res}}

\newcommand{\Aut}{\operatorname{Aut}}

\renewcommand{\a}{\alpha}
\renewcommand{\b}{\beta}
\newcommand{\om}{\omega}
\newcommand{\De}{\Delta}
\newcommand{\la}{\lambda}
\renewcommand{\th}{\theta}
\newcommand{\C}{{\Bbb C}}

\newcommand{\R}{{\Bbb R}}
\newcommand{\Z}{{\Bbb Z}}
\newcommand{\Q}{{\Bbb Q}}

\newcommand{\Ga}{\Gamma}

\newcommand{\wt}{\widetilde}
\newcommand{\ot}{\otimes}

\newcommand{\sub}{\subset}
\newcommand{\ed}{\qed\vspace{3mm}}

\newcommand{\Qcoh}{\operatorname{Qcoh}}

\newcommand{\cha}{\operatorname{char}}
\newcommand{\pa}{\operatorname{\partial}}

\title{Arithmetic mirror symmetry for genus 1 curves \\ with $n$ marked points}

\author[Yank{\i} Lekili and Alexander Polishchuk]{Yank{\i} Lekili \\ Alexander      Polishchuk}
\thanks{Y.L. is supported in part by the Royal Society and the NSF grant DMS-1509141. A.P. is supported in part by the NSF grant DMS-1400390. }

\address{King's College London \vspace{-0.15in}}
\address{University of Oregon}

\begin{document}

\begin{abstract}
    We establish a $\Z[[t_1,\ldots, t_n]]$-linear derived equivalence between the relative Fukaya
    category of the 2-torus with $n$ distinct marked points and the derived category of perfect complexes on the $n$-Tate curve.
    Specialising to $t_1= \ldots =t_n=0$ gives a $\Z$-linear derived equivalence between the Fukaya category
    of the $n$-punctured torus and the derived category of perfect complexes on the standard (N\'eron) $n$-gon.
    We prove that this equivalence extends to a $\Z$-linear derived equivalence between the
    wrapped Fukaya category of the $n$-punctured torus and the derived category of coherent sheaves
    on the standard $n$-gon. The corresponding results for $n=1$ were established in \cite{LP}.  
\end{abstract}
\maketitle

\centerline{\sc Introduction}

Over the past few decades, Kontsevich's Homological Mirror Symmetry (HMS) conjecture
(\cite{KontsevichICM}), which predicts a remarkable equivalence between an $A$-model category
associated to a symplectic manifold and a $B$-model category associated to a complex manifold, has
been studied extensively in many instances: abelian varieties, toric varieties, hypersurfaces in
projective spaces, etc. From the beginning, the case of the symplectic 2-torus played a special
role as this is the simplest instance of a Calabi-Yau manifold where concrete computations in the Fukaya category can be made, and the non-trivial nature of the HMS conjecture becomes manifest. Calculations provided in \cite{PZ},
\cite{Pol00}, \cite{Pol03} and \cite{Pol04} lead to a proof of the HMS for the 2-torus over a Novikov field (see \cite{AS} for
a streamlined exposition). Recently, various versions of the HMS have been verified when
the $A$-model category is a flavour of the Fukaya category of a, possibly non-compact, symplectic 2-manifold
(see for example \cite{SeidelGenus2}, \cite{efimov}, \cite{AAEKO}, \cite{bocklandt}, \cite{HLee}).

In \cite{LP}, the authors have explored an arithmetic refinement of the homological
mirror symmetry which relates exact symplectic topology to arithmetic algebraic
geometry. Concretely, they proved a derived equivalence of the Fukaya category
of the 2-torus, relative to a base-point $D=\{z \}$, with the category of perfect complexes
of coherent sheaves on the Tate curve over the formal disc $\mathrm{Spec\ } \Z [[t]]$.
Setting $t=0$, one obtains a derived equivalence over $\Z$ of the Fukaya
category of the punctured torus with perfect complexes on the nodal projective cubic $y^2 z
+ xyz = x^3$ in $\mathbb{P}^2_{\Z}$.

The main result of this paper is a generalization of the result of \cite{LP}, where we work with the Fukaya category of the 2-torus relative to a divisor $D= \{z_1,\ldots, z_n\}$ for $n>1$. Thus, on the symplectic side we consider the relative Fukaya category $\FF(\T,D)$ of a symplectic $2$-torus $\T$
with $n$ distinct marked points $z_1,\ldots,z_n$ (see Sec.\ \ref{relFuk}).  On the B-side we define a certain family of curves $T_n$ over $\Z[[t_1,\ldots,t_n]]$, which we call the {\it $n$-Tate curve}
(for $n=1$ we get the usual Tate curve; specializing to $t_1=\ldots=t_n$ we get the $n$-gon Tate curve of \cite[VII]{DR}).
The generic fiber of $T_n$ is a smooth elliptic curve, while changing the base to $\Z$ (by setting
$t_i=0$)
we get the curve $G_n$ over $\Z$, which we call following \cite{DR} 
the {\it standard $n$-gon} (or \emph{N\'eron $n$-gon}) over $\Z$. For $n>1$ this is the wheel of $n$ projective lines over $\Z$, while
for $n=1$ this is the projective line with $0$ and $\infty$ identified.

\medskip

\noindent
{\bf Theorem A.} {\it There is a $\Z[[t_1,\ldots,t_n]]$-linear equivalence of triangulated categories between 
the split-closed derived Fukaya category $D^\pi\FF(\T,D)$ of the 2-torus with $n$ 
marked points and the derived category of perfect
complexes $\Perf(T_n)$ on the n-Tate curve $T_n$. 
}

The key ingredient in the proof is the isomorphism between the moduli space of minimal $A_\infty$-structures on
a certain finite-dimensional graded algebra $E_{1,n}$ with a certain moduli space of curves, established in our
previous paper \cite{LPol}. Namely, the algebra $E_{1,n}$ arises as the self-$\Ext$-algebra of the coherent sheaf
$$\OO_C\oplus\bigoplus_{i=1}^n \OO_{p_i},$$
where $C$ is a projective curve of arithmetic genus one, and $p_1,\ldots,p_n$ are distinct smooth points of $C$, such
that $H^1(C,\OO(p_i))=0$ for each $i$. Taking into account higher products coming from the dg-enhancement
of the derived category of coherent sheaves on $C$, we see that each such $n$-pointed curve $(C,p_1,\ldots,p_n)$ gives rise
to an equivalence class of $A_\infty$-structures on $E_{1,n}$. 
Our result in \cite{LPol} implies that every $A_\infty$-structure on $E_{1,n}$ arises in such a way (in the case $n=2$
there is a problem with the characteristic $2$, so
we give a special argument in Section \ref{n=2-case-sec} below). To prove Theorem A we combine this with some direct computations in the relative Fukaya category and with the product formula for theta functions over the $n$-Tate curve.

We denote by $\T_0$ the $n$-punctured torus $\T\setminus\{z_1,\ldots,z_n\}$, and we denote by 
$\FF(\T_0)$ and $\WW(\T_0)$ its exact Fukaya and wrapped Fukaya categories (see Sec.\ \ref{relFuk}).
For every commutative ring $R$ we consider the $R$-linear $A_\infty$-categories $\FF(\T_0)\ot R$ and
$\WW(\T_0)\ot R$ obtained from $\FF(\T_0)$ and $\WW(\T_0)$ by the extension of scalars from $\Z$ to $R$.
We also set $G_{n,R}=G_n\times\Spec(R)$.

\medskip

\noindent
{\bf Theorem B.} {\it Let $R$ be a commutative Noetherian ring.

(i) There is an $R$-linear triangulated equivalence 
\[ D^\pi(\FF(\T_0) \ot R) \simeq  \Perf(G_{n,R}) \]
between the split-closed derived exact Fukaya category of $\T_0$ and the derived
category of perfect complexes on the standard $n$-gon. 

(ii)(=Theorem \ref{wrapped-equivalence-thm}) Assume that $R$ is regular.
Then the equivalence of (i) extends to an $R$-linear triangulated equivalence of derived categories 
\[  D^\pi(\WW(\T_0)\ot R) \simeq D^b(\Coh G_{n,R}) \]  
    between the wrapped Fukaya category of $\T_0$ and the
derived category of coherent sheaves on the standard $n$-gon. If $R=\Z$ or $R$ is a field then
$D^\pi(\WW(\T_0)\ot R)=D^b(\WW(\T_0)\ot R)$.  }

Theorem B(i) is a direct consequence of the proof of Theorem A. However, we give also another argument that
avoids computations, but uses instead a certain characterization of the $A_\infty$-structure associated with the
standard $n$-gon, established in Theorem \ref{B-wheel-char-thm}. To prove Theorem B(ii) we first check that
the natural Yoneda functor from $D^b (\WW(\T_0)\ot R)$ to the derived category of modules over $\FF(\T_0)\ot R$
is fully faithful (see Theorem \ref{A-fully-faithful-thm}), and then identify the image with 
$D^b(\Coh G_{n,R})$ using the equivalence of Theorem B(i).

In addition (working over a field) we prove that the derived categories of
the wrapped Fukaya category $\WW(\T_0)$ and of the Fukaya category $\FF(\T_0)$
are {\it Koszul dual}, in the following sense. In Sec.\ \ref{generators-Aside-sec} we find certain natural generators
\[  L = \bigoplus_{i=0}^n L_i \text{\ \   ,\ \   } \hat{L} = \bigoplus_{i=0}^n \hat{L}_i \] 
for the derived categories of
$\FF(\T_0)$ and $\WW(\T_0)$ respectively.
Denote $\mathscr{A}_0 = CF^*(L,L)$ and $\mathscr{B} =
CF^*(\hat{L}, \hat{L})$. The $\mathscr{A}_0-\mathscr{B}$-bimodule $CF^*(\hat{L},L)$ yields augmentations of
$\mathscr{A}_0$ and $\mathscr{B}$ over the semisimple ring $K = \bigoplus_{i=0}^n k$ (in the generalized
sense considered in \cite[Sec.\ 10]{Keller-derived-dg}), with the property that
\[ \mathrm{RHom}_{\mathscr{A}_0}(K,K) \simeq \mathscr{B} \text{  ,  } \
\mathrm{RHom}_{\mathscr{B}}(K,K) \simeq \mathscr{A}_0. \]
Note that there is a choice in the definition of $\hat{L}_0$: one has to choose one of the punctures. These $n$ choices
lead to different augmentations on $\AA_0$ and are permuted by a natural $\Z/n$-action (see Proposition \ref{aug-prop}).

\emph{Relation to previous work.} As alluded to before, the main results in this paper are 
generalizations of the main results of \cite{LP} to $n>1$. However, we should emphasize that classification
of the relevant moduli of $A_\infty$-structures given in \cite{LPol} and the computations of
the homogeneous coordinate rings are rather more complicated in the case $n>1$. We would like to
note that yet another proof of Theorem B(i) can be given more directly via the corresponding
result for $n=1$ proven in \cite{LP} by using the fact that there are
cyclic covering maps from an $n$-punctured
torus to a once-punctured torus, and the standard $n$-gon to the standard $1$-gon. These coverings
yield models of $\mathcal{F}(\T_0)$ and $\Perf(G_n)$ as a semi-direct product of $\Z/n\Z$ with a
subcategory of the corresponding category associated to the base of the covering
(such covering arguments are well known and go back to \cite{Seidelquartic}).
A sketch based on
this approach to prove Theorem B(i) appeared in the recent \cite{keating} (for $n=3$) during the
writing of this article. However, to complete the proof one needs to show that the two actions of $\Z/n$ 
arising from the covering picture, the one on 
the Fukaya category of the $1$-punctured torus and the one on $\Perf(G_1)$
(both given by tensoring with order $n$ line bundles),
are identified by the equivalence constructed in \cite{LP}.
We have not pursued this approach as it does not directly give a way of proving Theorem A where the weights of the marked points differ.  

In \cite{STZ}, the authors constructed an equivalence between
$\Perf(G_n)$ and a dg-category, called the ``constructible plumbing model''. The construction of
this dg-category is inspired by a suggestion of Kontsevich (\cite{Kon09}) that the Fukaya category
of $\T_0$ can be calculated via a category of constructible sheaves associated to the Lagrangian
skeleton (cf. \cite{NZ}). The authors of \cite{STZ} conjecture that their model is quasi-equivalent
to the $D^\pi \FF(\T_0)$.  Clause (i) of Theorem B implies this conjecture.  

In \cite{AAEKO}, a Landau-Ginzburg (B-model) mirror to $\WW(\T_0)$ was constructed (over $\C$) for some $n$. In \cite{bocklandt} certain non-commutative mirrors to $\WW(\T_0)$
for $n\geq 3$ were given. It may be interesting to compare these more directly to the derived
category of coherent sheaves on the standard $n$-gon.  

\emph{Outline of Sections.} In Section 1, we recall and extend our results from \cite{LPol} about
classification of $A_\infty$-structures on the associative algebra $E_{1,n}$ via the moduli of certain pointed curves
of arithmetic genus 1. Also, in Section \ref{wheel-char-sec} we give a characterization of the $A_\infty$-structure
associated with the standard $n$-gon.
In Section 2, we give a construction of the $n$-Tate curve and compute its
homogeneous coordinate ring. In Section 3, we first give generators for various flavours of Fukaya categories
associated with the pair $(\T,D)$. Then we relate these to generators of various derived categories
associated to the $n$-Tate curve, establishing the homological mirror symmetry results as stated in Theorems $A$ and $B$ by using classification results from
Section 2 and computing the homogeneous coordinate ring of the $n$-Tate curve inside the Fukaya category
of $(\T,D)$. Finally, in Section 4, we prove the Koszul duality result between
$\mathcal{F}(\T_0)$ and $\mathcal{W}(\T_0)$ by using the equivalences from Section 3. 

\noindent
{\it Conventions.} We use the terminology for generation of triangulated subcategories which may be non-standard
for algebraic geometers (but is standard in symplectic geometry): we say that a set of objects $S$ {\it split-generates}
(resp., {\it generates}) a triangulated category $\TT$ if the smallest thick subcategory (resp., triangulated subcategory)
of $\TT$ containing $S$ is the entire $\TT$. For a set of objects $S$ in a triangulated category we denote by
$\lan S\ran$ the thick subcategory split-generated by $S$.
By the {\it elliptic $n$-fold curve} we mean the projective curve of arithmetic genus $1$, which is
the union of $n$ generic projective lines passing through one point in $\P^{n-1}$, for $n\ge 3$,
the union of two $\P^1$'s intersecting at one tacnode point, for $n=2$, and the cuspidal plane cubic, for $n=1$.

\section{Curves of arithmetic genus $1$ with $n$ marked points and $A_\infty$-structures}

\subsection{Extension of some results from \cite{LPol}}\label{results-sec}

Below we always assume that $n\ge 2$.

Recall that in the previous work \cite{LPol} we 
studied the moduli stacks $\UU_{1,n}^{sns}$ of curves $C$ of arithmetic genus $1$
together with $n$ distinct smooth marked points $p_1,\ldots,p_n$ satisfying
\begin{itemize}
    \item $h^0(\mathcal{O}_C(p_i))=1$ for all $i$, and
    \item $\mathcal{O}_C(p_1+\ldots+p_n)$ is ample. 
\end{itemize}
We denote by $\wt{\UU}_{1,n}^{sns}\to \UU_{1,n}^{sns}$ 
the $\G_m$-torsor corresponding to a choice of a nonzero
element $\om\in H^0(C,\om_C)$. 

In the case $n\ge 3$ we identified $\wt{\UU}_{1,n}^{sns}$ with an explicit affine scheme
of finite type over $\Z$, while in the case $n=2$ we proved that over $\Z[1/2]$ one has 
$\wt{\UU}_{1,2}^{sns}\simeq \A^3$ (see \cite[Thm.\ 1.4.2]{LPol}). 

More precisely, in the case $n\ge 3$ we showed that
the ring $\OO(\wt{\UU}^{sns}_{1,n})$ is generated over $\Z$ by the functions defined as follows.
Let $(C,p_1,\ldots,p_n,\om)$ be the universal family, and let $h_{ij}$, for $i\neq j$, 
be elements of $H^0(C,\OO(p_i+p_j))$ such that
$\Res_{p_i}(h_{ij}\om)=1$. We normalize the elements $h_{1i}$ by the condition $h_{1i}(p_2)=0$ for $i\ge 3$
and $h_{12}(p_3)=0$. Then there is a relation of the form
$$h_{12}h_{13}^2-h_{12}^2h_{13}=ah_{12}h_{13}+bh_{12}+ch_{13}+d,$$
and the functions $a,b,c,d$ together with 
\begin{equation}\label{c-ij-eq}
c_{ij}:=h_{1i}(p_j), 
\end{equation}
where $i,j\ge 2$, $i\neq j$, generate the algebra
of functions on $\wt{\UU}^{sns}_{1,n}$ over $\Z$ (see \cite[Sec.\ 1.1]{LPol}).
Note that these functions have the following weights with respect to the $\G_m$-action:
\begin{equation}\label{U-n-gen-weights-eq}
wt(c_{ij})=wt(a)=1, \ wt(b)=wt(c)=2, \ wt(d)=3.
\end{equation}

In the case $n=2$ we showed (see \cite[Sec.\ 1.2]{LPol}) that $\wt{\UU}^{sns}_{1,2}\ot \Z[1/2]$ is isomorphic to the
affine $3$-space over $\Z[1/2]$ with the coordinates $\a$, $\b$ and $\ga$ of weights $2$, $3$ and $4$, so that the
affine part $C\setminus\{p_1,p_2\}$ of the universal curve is given by the equation
$$y^2-yx^2=\a(y-x^2)+\b x+\ga,$$
which is simply the unfolding of the tacnode.

Furthermore, for each $(C,p_1,\ldots,p_n)$ as above we considered
the generator 
$$G=\OO_C\oplus\OO_{p_1}\oplus\ldots\oplus \OO_{p_n}$$
of the perfect derived category of $C$. The standard dg-enhancement of this category allows to
construct a minimal $A_\infty$-algebra, equivalent to $\mathrm{REnd}^*(G)$.

The choice of a nonzero element in $H^0(C,\om_C)$ gives rise to a canonical identification of
the underlying associative algebra, $\Ext^*(G,G)$, with the algebra $E_{1,n}\ot k$, associated
with the following quiver $Q=Q_n$ with relations. 

\begin{figure}[!h] \centering
        \includegraphics[scale=1]{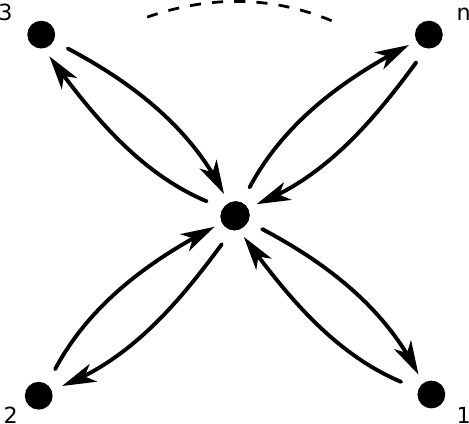}  \caption{The quiver $Q_n$}
            \label{fig0} \end{figure}

The path algebra $\Z[Q]$ is generated by $(A_i,B_i)$, where $A_i$ is the arrow
    from the central vertex to the vertex $i$ and $B_i$ is the arrow  in the opposite direction.
  We define a grading on $\Z[Q]$ by letting $|A_i|=0$ and $|B_i|=1$.
The ideal of relations $J$ is  generated by the relations
    $$B_iA_i=B_jA_j,\ \ A_iB_j=0 \ \text{ for } i\neq j.$$
We set
\begin{equation}\label{E1n-eq}
E_{1,n}:=\Z[Q]/J.
\end{equation}

Let $\MM_\infty=\MM_\infty(E_{1,n})$ be the functor on the category of commutative rings, associating to $R$
the set of gauge equivalence classes of minimal $A_\infty$-structures on $E_{1,n}\ot R$ (strictly unital with respect
to the idempotents $e_i\in E_{1,n}$ corresponding to the vertices) extending the given $m_2$ on $E_{1,n}\ot R$.

The map associating with a point $(C,p_1,\ldots,p_n,\om)\in \wt{\UU}_{1,n}^{sns}(R)$ the corresponding 
$A_\infty$-structure on $E_{1,n}\ot R$ (defined up to a gauge equivalence), extends to a morphism
of functors
\begin{equation}\label{curve-to-ainf-map}
\wt{\UU}_{1,n}^{sns}\to \MM_\infty.
\end{equation}
Namely, we can use the homological perturbation construction associated with a dg-model for
the subcategory generated by $\OO_C\oplus\OO_{p_1}\oplus\ldots\oplus\OO_{p_n}$, described in 
\cite[Sec.\ 3]{P-ainf}. This requires a choice of relative formal parameters $t_i$ at $p_i$ (compatible
with $\om$ so that $\Res_{p_i}(\om/t_i)=1$). However, this affects only the
choice of homotopies in the homological perturbation, and hence, a different choice 
leads to a gauge equivalent $A_\infty$-structure.

Furthermore, the morphism \eqref{curve-to-ainf-map} is compatible with the $\G_m$-action, where
the action of $\la\in\G_m$ on $\MM_\infty$ is given by the rescalings $m_n\mapsto \la^{2-n}m_n$.

We proved in \cite{LPol} that the map \eqref{curve-to-ainf-map} becomes an isomorphism once
we change the base to any field (of characteristic $\neq 2$, if $n=2$).
 Thus, over a field $k$, every minimal $A_\infty$-structure on $E_{1,n}\ot k$
 can be realised, up to gauge equivalence, in the derived category of some curve $(C,p_1,\ldots,p_n)$
 as above. 

The reason we restricted to working over a field in \cite{LPol} was the dependence on the results of
\cite{P-ainf} about moduli of $A_\infty$-structures. Using the new work \cite{P-ainf2} on the relative moduli
of $A_\infty$-structures we can now extend our results, so that everything works over $\Z$, for $n\ge 3$.
For $n=2$ the similar arguments work over $\Z[1/2]$, but we can also use the more ad hoc construction as in
\cite{LP} to establish a partial result we need over $\Z$ (see Section \ref{n=2-case-sec} below).


\begin{thm}\label{M1n-ainf-thm} Assume that $n\ge 3$. 
Then the functor $\MM_\infty$ is represented by an affine scheme of finite type over $\Z$. 
Furthermore, the morphism \eqref{curve-to-ainf-map} is an isomorphism. 
For $n=2$ the same assertions hold over $\Z[1/2]$.
\end{thm}

\Pf . Assume first that $n\ge 3$.
By \cite[Thm.\ 2.2.6]{P-ainf2}, the representability of our functor by an affine scheme follows from the vanishing
\begin{equation}\label{HH-E1n-van-eq}
HH^i(E_{1,n}\ot k)_{<0}=0
\end{equation}
for $i\le 1$ and any field $k$. For $i=1$ this is \cite[Eq.\ (2.2.1)]{LPol}, while for
$i=0$ this is clear. Furthermore, to check that $\MM_\infty$ is a closed subscheme of the scheme $\MM_n$ of
finite type, parametrizing $A_n$-structures, it is enough to have the vanishing
$$HH^2(E_{1,n}\ot k)_{<-d}=0$$
for some $d$ and all fields $k$. This is indeed the case for $d=3$ by \cite[Cor.\ 2.2.6]{LPol}.

Recall that the morphism \eqref{curve-to-ainf-map} is also compatible with the $\G_m$-action, where the
action on the global functions $\la\mapsto (\la^{-1})^*$ has positive weights on generators over $\Z$:
for $\wt{\UU}_{1,n}^{sns}$ the weights are given by \eqref{U-n-gen-weights-eq}, while for $\MM_\infty$
the weight of the coordinates of $m_n$, for $n>2$, is equal to $n-2$.
Thus, the morphism \eqref{curve-to-ainf-map} 
corresponds to a 
homomorphism of non-negatively finitely generated graded algebras over $\Z$
$$f:A\to B,$$ 
such that $A_0=B_0=\Z$, $f\ot\Q$ and $f\ot\Z/p$ are isomorphisms. In addition, we know that
each $B_n$ is a free $\Z$-module. Indeed, for $n\ge 5$ this follows from \cite[Cor.\ 1.1.7]{LPol}, while
for $n=3$ (resp., $n=4$) this follows from the identification of the moduli space with the affine space $\A^4$
(resp., $n=5$) given in \cite[Prop.\ 1.1.5]{LPol}. 
Thus, $C_n=\coker(f_n:A_n\to B_n)$ is a finitely generated
abelian group such that $C_n\ot \Q=0$ and $C_n\ot\Z/p=0$. Hence, $C_n=0$ and so $f_n$ is surjective.
Since $B_n$ is flat over $\Z$, for each $p$ we have an exact sequence
$$0\to (\ker f_n)\ot\Z/p\to A_n\ot\Z/p\to B_n\ot\Z/p\to 0$$
which shows that $\ker f_n\ot \Z/p=0$. Since $\ker f_n\ot\Q=0$, we derive that $\ker f_n=0$.

The argument in the case $n=2$ is similar. The vanishing of \eqref{HH-E1n-van-eq} in the case when $\cha(k)\neq 2$
follows from \cite[Eq.\ (2.1.4), Cor.\ 2.2.2]{LPol}. In the proof of the second assertion we use the identification
of $\wt{\UU}_{1,2}^{sns}\ot\Z[1/2]$ with $\A^3$, with coordinates of weights $2$, $3$ and $4$.
\ed

\begin{rem}
An alternative way to prove the second assertion of Theorem \ref{M1n-ainf-thm} is to mimic the proof of 
\cite[Thm.\ B]{P-ainf2} by first showing that the deformations over $\Z$ of the tacnode point of $\wt{\UU}_{1,2}^{sns}(k)$,
for any field $k$ of characteristic $\neq 2$, match the deformations of the $A_\infty$-structures.
\end{rem}

\subsection{Case $n=2$}\label{n=2-case-sec}

In the cases $n=1$ and $n=2$ the moduli stack $\wt{\UU}_{1,2}^{sns}$ over $\Spec(\Z)$ is not an affine scheme. The case
$n=1$ is considered in detail in \cite{LP}, so here we discuss the case $n=2$.

Given a curve $C$ of arithmetic genus $1$ with smooth marked points $p_1$, $p_2$, such that $H^1(C,\OO(p_i))=0$,
a choice of a nonzero element in $H^0(C,\om_C)$ is equivalent to a choice of a nonzero tangent vector at $p_1$
(see \cite[Lem.\ 1.1.1]{LPol}). Let $t_1$ be a formal parameter at $p_1$ compatible with this tangent vector.
Then there exist elements $x\in H^0(C,\OO(p_1+p_2)$ and $y\in H^0(C,\OO(2p_1))$ such that
$$x\equiv \frac{1}{t_1}+k[[t_1]], \ \ y\equiv \frac{1}{t_1^2}+t_1^{-1}k[[t_1]]$$
at $p_1$, defined uniquely up to adding a constant. It is easy to see that then the elements
$$1,x,x^2,y,xy$$
form a basis of $H^0(C,\OO(3p_1+2p_2))$. Note that $y^2-yx^2\in H^0(C,\OO(3p_1+2p_2))$.
Hence, we should have a relation of the form
$$y^2=yx^2+axy+by+b'x^2+cx+d.$$
Adding a constant to $y$ we can make the term with $x^2$ to disappear, so there is a unique choice of $y$, so that the 
above relation takes form
\begin{equation}\label{tacnode-unfold-eq}
y^2=yx^2+axy+by+cx+d.
\end{equation}
There remains ambiguity in a choice of $x$: changing $x$ to $x-\a$ leads to the transformation
\begin{equation}\label{additive-action-eq}
(a,b,c,d)\mapsto (a+2\a, b+\a a+\a^2, c, d+\a c).
\end{equation}

Let us consider the quotient stack $\A^4/\G_a$, where the action of the additive group on $\A^4$ is given by 
\eqref{additive-action-eq}. Note that the action \eqref{additive-action-eq} is compatible with the $\G_m$-action on
$\A^4$ such that $wt(a)=1$, $wt(b)=2$, $wt(c)=3$, $wt(d)=4$ and the standard $\G_m$-action on $\G_a$ (so $wt(\a)=1$).

\begin{prop} One has a natural isomorphism $\wt{\UU}_{1,2}^{sns}\simeq \A^4/\G_a$, compatible with
the $\G_m$-actions.
\end{prop}

\Pf . Consider the $\G_a$-torsor over $\wt{\UU}_{1,2}^{sns}$ corresponding to a choice
of a function $x\in H^0(C,\OO(p_1+p_2))$ such that the polar part of $x$ at $p_1$ is $\frac{1}{t_1}$. Then
as we have seen above, we can uniquely find $y\in H^0(C,\OO(2p_1))$, such that the defining equation
is of the form \eqref{tacnode-unfold-eq}. Conversely, starting from the affine curve $\Spec(A)$ given by such an equation
we construct the projective curve $C$ by taking $\Proj \RR(A)$, where $\RR(A)$ is the Rees algebra associated
with the filtration by degree, where $\deg(x)=1$, $\deg(y)=2$. As in \cite[Thm.\ 1.2.4]{P-ainf}, one easily
checks that this gives an isomorphism between our $\G_a$-torsor over $\wt{\UU}_{1,2}^{sns}$ and $\A^4$.
\ed

\begin{prop}\label{n=2-ainf-surj-prop} 
The map $\A^4(R)\to \MM_\infty(R)$ associating with $(a,b,c,d)\in R^4$ the $A_\infty$-structure coming from the corresponding
curve in $\wt{\UU}_{1,2}^{sns}$, is surjective in the following cases: (i) $R$ is any field; (ii) $R$ is an integral domain
with the quotient field of characteristic zero.
\end{prop}

\Pf . We mimic the proof of \cite[Thm.\ C]{LP} (see \cite[Sec.\ 5.3]{LP}).

Let $W$ denote the tangent space to $\A^4$ at $0$.
The derivative of the action \eqref{additive-action-eq} at $0$ gives a map
$$d:\Lie(\G_a)\to W,$$
which we can easily compute: $d(\pa_x)=2\pa_a$.

As in \cite[Sec.\ 5.3]{LP}, the map from $\A^4/\G_a$ to the functor of $A_\infty$-structures, induces at the infinitesimal level
a chain map $\kappa$ from the dg Lie algebra $[\Lie(\G_a)\to W]$ (living in degrees $0$ and $1$) to
the shifted Hochschild cochain complex $CH^*(E_{1,2})_{<0}[1]$ (truncated in negative internal degrees).
Similarly to \cite[Thm.\ 5.5]{LP} we claim that this map induces an isomorphism of cohomology in degrees $0$ and $1$,
when tensored with any field $k$, or over $\Z$.

Indeed, first let us check that
$$H^1(\kappa\ot k): \coker(d\ot k)\to HH^2(E_{1,2}\ot k)_{<0}$$
is an isomorphism.
Note that the source can be viewed as the tangent space to the deformations of the tacnode 
curve $(C_{tn},p_1,p_2,\om)$ in $\wt{\UU}_{1,2}^{sns}$, and the map itself as the tangent map to the morphism of deformation functors associating to a deformation of $(C_{tn},p_1,p_2,\om)$ the corresponding family of $A_\infty$-structures on $E_{1,2}$.
Now \cite[Prop.\ 4.3.1]{P-ainf}, applied to families over $k[\eps]/(\eps^2)$ implies that the map $H^1(\kappa\ot k)$ is
injective. On the other hand, we claim that the source and the target have the same dimension, which is equal to $3$ when
$\cha(k)\neq 2$ and is equal to $4$, when $\cha(k)=2$. Indeed, for the source this is easy to see, while for the target this
follows from \cite[Cor.\ 2.2.6]{LPol} in the case $\cha(k)\neq 2$. In the case $\cha(k)=2$ it suffices to show that
$\dim HH^2(E_{1,2}\ot k)_{<0}\le 4$, which follows from Lemma \ref{tacnode-lem} below.
Thus, we conclude that $H^1(\kappa\ot k)$ is an isomorphism.

Now we turn to showing that
$$H^0(\kappa\ot k): \ker(d\ot k)\to HH^1(E_{1,2}\ot k)_{<0}$$
is an isomorphism.
If $\cha(k)\neq 2$ then $\ker(d\ot k)=0$ and $HH^1(E_{1,2}\ot k)_{<0}=0$ (see \cite[Cor.\ 2.2.2]{LPol}).
The interesting case is when $k$ has characteristic $2$. Then $\ker(d\ot k)=\lan\pa_x\ran$ and $\pa_x$ maps under $\kappa$ to  the nonzero element of the one-dimensional space
$$HH^1(E_{1,2})_{<0}\simeq HH^1(C_{tn})_{<0}\simeq H^0(C_{tn},\TT)_{<0}$$
corresponding to the global vector field $\pa_x$ of weight $-1$ on the tacnode $C_{tn}$
(see \cite[Prop.\ 2.1.3, Lem.\ 1.5.2]{LPol} and the proof of \cite[Cor.\ 2.2.2]{LPol}).

The rest of the proof goes as in \cite[Sec.\ 5.3]{LP}.
\ed

We have used the following result.

\begin{lem}\label{tacnode-lem} Let $k$ be a field of characteristic $2$.
Then $\dim HH^2(E_{1,2}\ot k)=4$.
\end{lem}

\Pf . Let $C=C_{tn}$ be the (projective) tacnode curve over $k$, equipped with a pair of smooth points $p_1,p_2$
on each component.
Then we have an isomorphism $HH^*(E_{1,2}\ot k)_*\simeq HH^*(C)$, so that the second grading is induced
by the $\G_m$-action on $C$. Indeed, this follows from the homotopical triviality of the $A_\infty$-structure on 
$E_{1,2}\ot k$ associated with $C$ (see \cite[Lem. 2.1.2]{LPol} and \cite[Prop.\ 4.4.1]{P-ainf}). 

We have an exact sequence
$$0\to H^1(C,\TT)\to HH^2(C)\to HH^2(U)\to 0,$$
where $U=C\setminus\{p_1,p_2\}$ (see \cite[Sec.\ 4.1.3]{LP}).
We claim that $H^1(C,\TT)=0$. Indeed, let $V=C\setminus q$, where $q$ is the singular point.
Then $(U,V)$ is an affine covering of $C$. Let $x_1,x_2$ be the natural coordinates on the components of $U$
(both vanishing at $q$). Derivations of $\OO(U\cap V)$ are just pairs $(P_1(x_1,x_1^{-1})\pa_{x_1},P_2(x_2,x_2^{-1})\pa_{x_2})$.
Derivations of $\OO(U)$ are those pairs, for which there exist constants $a,b$ such that
$P_i\equiv a+bx_i\mod x_i^2k[x_1]$, for $i=1,2$ (see \cite[Lem.\ 1.5.2]{LPol}). On the other hand, when $P_i$ are linear combinations of $x_i^n$ for $n\le 2$
then the corresponding derivation extends to $V$. This immediately implies that every derivation of $\OO(U\cap V)$ is
a sum of those extending either to $U$ or to $V$, hence, $H^1(C,\TT)=0$.
Finally, $U$ is an affine plane curve $k[x,y]/(y^2-yx^2)$, so $HH^2(U)$ is given by the corresponding Tjurina algebra
$k[x,y]/(x^2,y^2)$, which is $4$-dimensional.
\ed

\subsection{$A_\infty$-characterization of the wheel of projective lines}\label{wheel-char-sec}

For a commutative ring $R$ we consider 
$$G_{n,R}:=G_n\times \Spec(R),$$ 
the {\it standard $n$-gon over} $R$. Note that it has natural smooth $R$-points
$p_1,\ldots,p_n$ (corresponding to the point $1\in \P^1$ on each component, where the points $0$ and $\infty$
are used for gluing). Furthermore, there is a natural choice of a section $\om$ of the dualizing sheaf of $G_{n,R}$
over $R$ (see \cite[Ex.\ 1.1.9]{LPol}), so
we can view $(G_{n,R},p_1,\ldots,p_n,\om)$ as a family in $\wt{\UU}_{1,n}^{sns}(R)$. By abuse of notation we
 will sometimes refer to this family simply as $G_{n,R}$.

Now let $k$ be a field.
We are going to give several characterizations of the equivalence class of minimal $A_\infty$-structures on $E_{1,n}\ot k$
associated with $G_{n,k}$ (via the morphism \eqref{curve-to-ainf-map}).

Note that for every subset $S\sub\{1,\ldots,n\}$ we have a natural subquiver in $Q_n$ such that the
corresponding subalgebra is isomorphic to $E_{1,|S|}$. In particular, we have $n$ subquivers $Q_1(i)\sub Q_n$ (where $i=1,\ldots,n$) that give embeddings of $E_{1,1}$ into $E_{1,n}$. 
Now given a minimal $A_\infty$-structure $m_\bullet$ on $E_{1,n}$, for each $i$ we have a well defined restriction $m_\bullet|_{Q_1(i)}$ which is a minimal $A_\infty$-structure on $E_{1,1}$ (recall that
we consider $A_\infty$-structures unital with respect to the idempotents in $E_{1,n}$).
On the other hand, every such $m_\bullet$ gives a structure of right $A_\infty$-module on $P_i=e_iE_{1,n}$, $i=0,\ldots,n$
(where $e_0,\ldots,e_n$ are the idempotents in $E_{1,n}$ corresponding to the vertices in $Q_n$).

\begin{thm}\label{B-wheel-char-thm} 
Let $k$ be a field, and let $m^{wh}_\bullet$ be the minimal $A_\infty$-structure on $E_{1,n}\ot k$ associated with
$G_{n,k}$, where $n\ge 2$. 
Then $m^{wh}_\bullet$ is characterized uniquely (among the $A_\infty$-structures we consider in Theorem \ref{M1n-ainf-thm})
up to gauge equivalence and up to $\G_m$-action,
by the following conditions (i) and either (ii) or (ii'):

\noindent
(i) for every $i=1,\ldots,n$, the restriction $m^{wh}_\bullet|_{Q_1(i)}$ is not homotopically trivial;

\noindent
(ii) $\dim HH^2(E_{1,n},m^{wh}_\bullet)=n$;

\noindent
(ii') the subcategories $\lan P_0,P_i\ran$ split-generated by the right $A_\infty$-modules $P_0$ and $P_i$
(where the $A_\infty$-structure comes from $m^{wh}_\bullet$) are all distinct for $i=1,\ldots,n$.

Furthermore, for all minimal $A_\infty$-structures $m_\bullet$ on $E_{1,n}$ satisfying (i), one has
$$\dim HH^2(E_{1,n},m_\bullet)\le n.$$
\end{thm}

\begin{lem}\label{HH-wheel-lem} 
One has $\dim HH^2(G_{n,k})=n$.
\end{lem}

\Pf . Let us write $G_n$ instead of $G_{n,k}$ for brevity.
We have an isomorphism 
$$HH^2(G_n)\simeq \Ext^1(\bL_{G_n/k},\OO_{G_n})\simeq \Ext^1(\Om_{G_n /k},\OO_{G_n}),$$
where $\bL_{G_n /k}$ is the cotangent complex, which in this case is isomorphic to $\Om_{G_n /k}$
(since $G_n$ is a locally complete intersection). Let us pick $n$ smooth points $p_1,\ldots,p_n\in G_n$, one on each component,
and let $D=p_1+\ldots+p_n$.
Then the exact sequence
$$0\to \OO_{G_n}(-D)\to \OO_{G_n} \to \OO_D\to 0$$
induces a long exact sequence
\begin{align*}
    &0\to\Hom(\Om_{G_n /k},\OO_{G_n}(-D))\to\Hom(\Om_{G_n /k},\OO_{G_n})\to \Hom(\Om_{G_n /k},\OO_D)\to \\
    &\Ext^1(\Om_{G_n /k},\OO_{G_n}(-D))\to
    \Ext^1(\Om_{G_n/k},\OO_{G_n})\to \Ext^1(\Om_{G_n/k},\OO_D)\to\ldots
\end{align*}
Since $\Om_{G_n/k}$ is locally free near $D$, we have $\Ext^1(\Om_{G_n/k},\OO_D)=0$.
On the other hand, we have 
$$\Hom(\Om_{G_n/k},\OO_{G_n}(-D))=H^0(G_n,\TT(-D))=0,$$
while $H^0(G_n,\TT)$ is $n$-dimensional. Hence the restriction map
$$H^0(G_n,\TT)\to H^0(G_n,\TT|_D)$$
is an isomorphism, so from the above long exact sequence we obtain an isomorphism
$$\Ext^1(\Om_{G_n/k},\OO_{G_n}(-D))\simeq\Ext^1(\Om_{G_n/k},\OO_{G_n}).$$
But the space $\Ext^1(\Om_{G_n/k},\OO_{G_n}(-D))$ is the tangent space to $\ov{\MM}_{1,n}$ at the point
$(G_n,p_\bullet)$, so it is $n$-dimensional (see \cite{DM}).
\ed

\noindent
{\it Proof of Theorem \ref{B-wheel-char-thm}.}
By Theorem \ref{M1n-ainf-thm} (resp., Proposition \ref{n=2-ainf-surj-prop} for $n=2$), any minimal $A_\infty$-structure on $E_{1,n}$ comes from a curve
$(C,p_1,\ldots,p_n,\om)$ defining a point of $\wt{\UU}^{sns}_{1,n}(k)$. 
Recall that over an algebraically closed field any such pointed curve
coincides with its minimal subcurve of arithmetic genus $1$ without disconnecting nodes
(see the proof of \cite[Thm.\ 1.5.7]{LPol}). Hence, by \cite[Lem.\ 3.3]{Smyth-I},
$$(\ov{C},\ov{p}_1,\ldots,\ov{p}_n):=(C,p_1,\ldots,p_n)\times_k \ov{k}$$ 
(where $\ov{k}$ is an algebraic closure of $k$) 
is either an elliptic $m$-fold curve with $m\le n$, or the standard $m$-gon $G_{m,\ov{k}}$ with $m\le n$
(the case $m=1$ being the irreducible nodal curve), or a smooth elliptic curve. 

Recall that by definition of
$\wt{\UU}^{sns}_{1,n}$ there is at least one marked point on each irreducible component of $\ov{C}$. 
Assume first that $\ov{C}$ is an elliptic $m$-fold curve. Without loss of generality we can assume that
the points $\ov{p}_1,\ldots,\ov{p}_m$ all lie on different irreducible components of $\ov{C}$. Then
the group $\Aut(\ov{C},\ov{p}_1,\ldots,\ov{p}_m)$ is isomorphic to $\ov{k}^*$, hence, by the Hilbert Theorem 90,
we get that $(C,p_1,\ldots,p_m)$ is isomorphic to the elliptic $m$-fold curve over $k$.
Now we claim that restricting $m_\bullet$
to the subquiver $Q_m\sub Q_n$ corresponding to the points $p_1,\ldots,p_m$,
we obtain a homotopically trivial $A_\infty$-structure. 
Indeed, this follows from the fact that the operation of restricting to $Q_m$ corresponds to the natural morphism
$\wt{\UU}^{sns}_{1,n}\to \wt{\UU}^{sns}_{1,m}$ forgetting the last $n-m$ marked points (see \cite[Rem.\ 2.2.10]{LPol}),
together with the fact that the elliptic $m$-fold curve in $\wt{\UU}^{sns}_{1,m}$ corresponds to the trivial
$A_\infty$-structure. 

Hence, for any $A_\infty$-structure satisfying (i), the curve $\ov{C}$ is either smooth or isomorphic
to $G_{m,\ov{k}}$ with $m\le n$. Note that in the former case $C$ itself is smooth, while in the latter case, using the fact that
the standard $m$-gon, equipped with one marked point on each component, has no automorphisms, we see
that $C$ itself is isomorphic to $G_{m,k}$ over $k$. In both cases applying forgetful morphisms to
$\wt{\UU}^{sns}_{1,1}$ we get a curve which is either smooth or nodal, hence condition (i) is satisfied. 

Since for smooth $C$ we have $\dim HH^2(C)=1$, the
characterization of $m^{wh}_\bullet$ using conditions (i) and (ii) now follows from Lemma \ref{HH-wheel-lem}.

Note that if $C$ is irreducible then for any pair of points $p_1,p_2\in C$ we have
$$\Perf(C)=\lan \OO_C,\OO_{p_1}\ran=\lan \OO_C,\OO_{p_2}\ran$$
(see the proof of \cite[Prop.\ 4.3.1]{P-ainf} or Corollary \ref{Perf-C-gen-cor} below),
so the condition (ii') does not hold in this case.
Thus, for the characterization using conditions (i) and (ii') we need to show that for a standard $m$-gon $C$ (with $m\ge 2$)
and a pair of smooth points $p_1,p_2\in C$, the subcategories $\CC_1=\lan \OO_C,\OO_{p_1}\ran$ and
$\CC_2=\lan \OO_C,\OO_{p_2}\ran$ are the same if and only if $p_1$ and $p_2$ lie on the same irreducible component of $C$.
Assume first that $p_1\in C_1$, $p_2\in C_2$, where $C_1$ and $C_2$ are different components of $C$.
Consider the (derived) restriction functor $i_1^*:\Perf(C)\to D^b(C_1)$. Then $i_1^*\CC_2\sub \lan \OO_{C_1}\ran$
which does not contain $i_1^*\OO_{p_1}\simeq\OO_{p_1}$, so $\CC_1\neq\CC_2$.
Now assume that $p_1$ and $p_2$ lie on the same component $C_1$. There exists a morphism
$f:C\to \ov{C}$, where $\ov{C}=G_1$ 
is the irreducible rational curve with one node, contracting all components different from $C_1$, and such that $f|_{C_1}:C_1\to\ov{C}$ is the normalization map. Let $x_i=f(p_i)\in \ov{C}$, $i=1,2$. Then
$\OO_{p_i}\simeq f^*\OO_{x_i}$ for $i=1,2$. Furthermore, we have
$$\Perf(\ov{C})=\lan \OO_{\ov{C}}, \OO_{x_1}\ran=\lan \OO_{\ov{C}}, \OO_{x_2}\ran.$$
Since $\Perf(\ov{C})$ is idempotent-complete,
it is enough to check that the functor $f^*:\Perf(\ov{C})\to \Perf(C)$ is fully faithful.
By the projection formula, this would follow from the equality $Rf_*\OO_{C}\simeq \OO_{\ov{C}}$.
Let $q\in\ov{C}$ be the node. Then $f$ is an isomorphism over the complement to $q$, so
$R^1f_*\OO_C$ is supported at $q$. It is also easy to see that $R^0f_*\OO_C\simeq \OO_{\ov{C}}$.
Hence, from the Leray spectral sequence we deduce that
$$H^1(C,\OO)\simeq H^1(\ov{C},\OO)\oplus H^0(\ov{C}, R^1f_*\OO_C).$$
Since both $H^1(C,\OO)$ and $H^1(\ov{C},\OO)$ are $1$-dimensional this implies the vanishing of $R^1f_*\OO_C$.
\ed

\section{The $n$-Tate curve}

\subsection{Construction}

The $n$-Tate curve we are going to construct will be a family of curves over $\Z[[t_1,t_2,\ldots,t_n]]$,
generically smooth and with the standard $n$-gon $G_n$ as the
specialization at $t_1=\ldots=t_n=0$.
This is a natural generalization of the construction of the Tate curves in \cite{DR} (obtained as
the specialization $t_1=\ldots=t_n$ from our construction). In particular, for $n=1$ we get the standard Tate curve. 

The construction goes through the same steps as in the case of the usual Tate curve:
we first construct a formal scheme over $\Z[[t_1,t_2,\ldots,t_n]]$ as the quotient of the formal completion of certain toric scheme 
$\TT_{\infty,n}$ by the action of $\Z$, and then apply Grothendieck's existence theorem.

We start by describing the scheme $\TT_{\infty,n}$ in terms of gluing open affine pieces $U_i$ numbered by $i\in\Z$.
Consider the periodic set of independent variables $t_i$, $i\in\Z$, where $t_{i+n}=t_i$, and set
$$\Z[t]_n=\Z[\ldots,t_i,t_{i+1},\ldots].$$ 
Then we define open affine pieces by 
$$U_i=\Spec \Z[t]_n[X_i,Y_{i+1}]/(X_iY_{i+1}-t_i).$$ 
The intersections $V_i=U_{i-1}\cap U_i$ correspond to setting $X_iY_i=1$,
so $V_i$ is the distinguished open $X_i\neq 0$ in $U_i$ (resp., $Y_i\neq 0$ in $U_{i-1}$). Thus,
$$V_i=\Spec \Z[t]_n[X_i,X_i^{-1}].$$
Let $\TT_{\infty,n}$ be the scheme over $\Z[t]_n$ obtained by gluing these open pieces (note that 
$U_i$ and $U_j$ will intersect not only for $|i-j|=1$).
Note that the central fiber $\TT_0$ (obtained by setting $t_1=\ldots=t_n=0$) is just the infinite chain of projective lines.
On the other hand,
since $U_i$ are just affine spaces over $\Z$, we see that $\TT_{\infty,n}$ is an integral scheme.
Note that we have the following relations between rational functions on $\TT_{\infty,n}$:
\begin{equation}\label{XYt-rel}
X_{i-1}=X_it_{i-1}, \ \ Y_{i+1}=Y_it_i, \text{ for } i\in\Z.
\end{equation}
In particular, all $X_i$ with $i\le i_0$ and all $Y_i$ with $i>i_0$ are regular functions on $U_{i_0}$. 

We have a natural action of $\Z$ on $\TT_{\infty,n}$, so that the automorphism $\tau$ corresponding
to $1\in\Z$ maps $U_{i-1}$ to $U_{i}$, and
$$\tau^*t_j=t_{j-1}, \ \ \tau^*X_i=X_{i-1}, \ \ \tau^*Y_i=Y_{i-1}.$$

Next, we take the formal neighborhood of $t_1=\ldots=t_n=0$ in $\TT_{\infty,n}$, 
and take the quotient by the action of the subgroup $n\Z\sub \Z$.
Note that this subgroup acts trivially on $\Z[t]_n$, so we get a formal curve 
$\hat{\TT}_n$ over 
$\Z[[t]]_n=\Z[[t_1,\ldots,t_n]]$.

Note that for each $i\in\Z$ we have a section of the central fiber
\begin{equation}\label{sigma-i-eq}
\si_i:\Spec(\Z)\to \TT_0\cap V_i 
\end{equation}
given by $X_i=-1$.

\subsection{Toric point of view}

As in the case $n=1$, we can view $\TT_{\infty,n}$ as a toric $\Z$-scheme of infinite type
associated with an infinite fan. Namely, we observe that there is an open embedding of a torus $T$ (over $\Z$) into $\TT_{\infty,n}$
given by 
$$T=\Spec(\Z[t_1,t_1^{-1},\ldots,t_n,t_n^{-1},X_0,X_0^{-1}])\sub V_0.$$
We have 
$$X_i|_T=\begin{cases}\frac{X_0}{t_0t_1\ldots t_{i-1}}, & i>0, \\ t_it_{i+1}\ldots t_{-1}X_0, &i<0,\end{cases}$$
(and $Y_i=X_i^{-1}$ on $T$). Let us consider the $(n+1)$-dimensional real vector space $N_\R$ with
the coordinates $(x,y_1,\ldots,y_n)$. We identify indices of coordinates $y_i$ with $\Z/n$.
We define simplicial cones $C_i\sub N_\R$ for $i\in\Z$ by 
$$C_i=\begin{cases}
\{(x,y_1,\ldots,y_n)\in\R\times \R_{\ge 0}^n \ |\ y_0+y_1+\ldots+y_{i-1}\le x\le y_0+y_1+\ldots+y_i\}, & i>0,\\
\{(x,y_1,\ldots,y_n)\in\R\times \R_{\ge 0}^n \ |\ 0\le x\le y_0\}, & i=0,\\
\{(x,y_1,\ldots,y_n)\in\R\times \R_{\ge 0}^n \ |\ -y_{-1}\le x\le 0\}, & i=-1,\\
\{(x,y_1,\ldots,y_n)\in\R\times \R_{\ge 0}^n \ |\ -y_i-\ldots-y_{-1}\le x\le -y_{i+1}-\ldots-y_{-1}\}, & i<-1.
\end{cases}
$$
Let $M_\R$ be the dual vector space to $N_\R$ with the basis $(f,e_1,\ldots,e_n)$ dual to the coordinates
$(x,y_1,\ldots,y_n)$, and let $M\sub M_\R$ be the $\Z$-lattice generated by the basis vectors.
We identify elements of $M$ with characters of the torus $T$ by 
$t_i=z^{e_i}$, $X_0=z^f$. Then 
$$X_i=\begin{cases} z^{f-e_0-\ldots-e_{i-1}}, & i>0, \\ z^f, & i=0, \\ z^{f+e_i+e_{i+1}+\ldots+e_{-1}}, & i<0,
\end{cases}$$
and we get
$$U_i=\Spec(\Z[C_i^\vee\cap M]),$$
Thus, we can identify $\TT_{\infty,n}$ with the toric scheme associated with the fan generated by the cones
$C_i$.

The automorphism $\tau$ of $\TT_{\infty,n}$ corresponds to the action of the linear automorphism
\begin{equation}\label{tau-N-eq}
\tau_N:N_\R\to N_\R:(x,y_0,\ldots,y_{n-1})\mapsto (x+y_{n-1},y_{n-1},y_0,\ldots,y_{n-2}).
\end{equation}
which preserves the lattice $N$ and sends $C_i$ to $C_{i+1}$.

\subsection{Polarization}

Let us define the line bundle $L$ over $\TT_{\infty,n}$ by setting $L|_{U_i}=\OO_{U_i}z_i$ (where $z_i$ is
a formal symbol) and defining the transitions on $V_i=U_{i-1}\cap U_i$ by the rule
\begin{equation}\label{X-z-rel}
z_{i-1}=X_iz_i.
\end{equation}

We can realize $L$ as a subsheaf in the sheaf of rational functions $\KK_{\TT_{\infty,n}}$  by identifying $z_i$ with the
character $z^{w_i}$ of the torus $T$, where $(w_i\in M)_{i\in\Z}$ is the unique collection of lattice points such that
$z^{w_{i-1}}=X_iz^{w_i}$ and $w_0=0$.
More explicitly, we have
$$w_i=\begin{cases} -if+ie_0+(i-1)e_1+\ldots+e_{i-1}, & i>0,\\ 0, & i=0,\\ f, & i=-1, \\
-if+e_{i+1}+2e_{i+2}+\ldots+(-i-1)e_{-1}, & i<-1.\end{cases}
$$

We need to lift the $\Z$-action on $\TT_{\infty,n}$ to $L$. For this let us
consider the affine transformation $\tau$ of $M_\R=\R\times\R^n$ given by 
\begin{equation}\label{tau-eq}
\tau(v)=\tau_M(v)+f, 
\end{equation}
where $\tau_M:M_\R\to M_\R$ is the linear transformation dual to $\tau_N$ (see \eqref{tau-N-eq}), so that
$$\tau_M(e_i)=e_{i-1}, \ \tau_M(f)=f+e_{-1}.$$
Then $\tau$ preserves the lattice $M$ and satisfies $\tau(w_i)=w_{i-1}$, hence, it gives the required lifting of the $\Z$-action
to $L$.

Let $C_i^\vee\sub M_\R$ be the cone dual to $C_i$, and let
$$\De=\cap_{i\in\Z}(w_i+C_i^\vee).$$
Note that $\tau_M(C_i^\vee)=C_{i-1}^\vee$, so $\tau(w_i+C_i^\vee)=w_{i-1}+C_{i-1}^\vee$, hence,
$$\tau(\De)=\De.$$

Let us consider the piecewise linear function $\phi:\R\to\R_{\ge 0}$ given by
\begin{equation}\label{phi-def}
\phi(t)=tk-\frac{k(k+1)}{2} \ \text{ for} \ k\le t\le k+1, \ k\in\Z.
\end{equation}
Note that this function has the property $\phi(\frac{1}{m}\Z)\sub \frac{1}{m}\Z$, which will be useful for the theory
of theta functions (see Section \ref{mult-theta-sec}). It also satisfies the quasi-periodicity 
\begin{equation}\label{phi-t+1-eq}
\phi(t+1)=\phi(t)+t.
\end{equation}

\begin{lem}\label{wi-coef-lem} We have the following formula for $w_i$:
$$w_i=\sum_{j=0}^{n-1} n\phi(\frac{j-i}{n})e_j - if.$$
\end{lem}

\Pf . Let $w'_i$ denote the right-hand side of this formula. Since $w'_0=w_0=0$, it is enough to check that 
$\tau(w'_i)=w'_{i-1}$. We have
\begin{align*}
&\tau(\sum_{j=0}^{n-1} n\phi(\frac{j-i}{n})e_j - if)=
(n\phi(\frac{-i}{n})-i)e_{n-1}+
\sum_{j=1}^{n-1} n\phi(\frac{j-i}{n})e_{j-1} - (i-1)f=\\
&(n\phi(\frac{-i}{n})-i)e_{n-1}+
\sum_{j=0}^{n-2} n\phi(\frac{j-i+1}{n})e_j - (i-1)f,
\end{align*}
so the statement reduces to
$$n\phi(\frac{-i}{n})-i=n\phi(\frac{n-i}{n}),$$
which follows from \eqref{phi-t+1-eq}.
\ed

\begin{lem}\label{phi-lem} 
(i) Let $(x^*,y^*_0,\ldots,y^*_{n-1})$ be the coordinates on $M_\R$. Then the cone $C_i^\vee\sub M_\R$
is described by the inequalities 
\begin{equation}
y^*_j+(\lfloor \frac{i-j}{n}\rfloor+1)x^*\ge 0, \ \ 
y^*_j+(\lfloor \frac{i-j-1}{n}\rfloor+1)x^*\ge 0, \ \ j=0,\ldots,n-1.
\end{equation}

\noindent (ii) The set $\De\sub M_\R$ is described by the inequalities 
$$y^*_j\ge n\phi(\frac{x^*+j}{n}), \ j=0,\ldots, n-1.$$
\end{lem}

\Pf . (i) It suffices to check that the cone $C_i$ is generated by the vectors
$$e_j^*+(\lfloor \frac{i-j}{n}\rfloor+1)f^*, \ \ e_j^*+ (\lfloor \frac{i-j-1}{n}\rfloor+1)f^*, \ \ j=0,\ldots,n-1,$$
where $(f^*,e_0^*,\ldots,e_j^*)$ is the standard basis in $N_\R$.
In the case $i=0$ these are the vectors $e_0^*+f^*,e_0^*,\ldots,e_{n-1}^*$, so the assertion is clear.
The general case follows easily using the automorphism $\tau_N$ that sends $C_i$ to $C_{i+1}$.

\noindent
(ii) Let us unravel the conditions $v-w_i\in C_i^\vee$, $i\in\Z$. Using part (i) and Lemma \ref{wi-coef-lem} we get the
following inequalities on the coordinates $(x^*,y^*_0,\ldots,y^*_{n-1})$ of $v\in M_\R$:
$$y^*_j-n\phi(\frac{j-i}{n})+(\lfloor \frac{i-j}{n}\rfloor+1)(x^*+i)\ge 0, \ \ 
y^*_j-n\phi(\frac{j-i}{n})+(\lfloor \frac{i-j-1}{n}\rfloor+1)(x^*+i)\ge 0, 
$$
for $j=0,\ldots,n-1$.
At this point it is convenient to introduce the function
$$\psi(q,t)=qt-\frac{q(q+1)}{2},$$
so that 
\begin{equation}\label{max-phi-eq}
\phi(t)=\max_{q\in\Z}\psi(q,t)=\psi(\lfloor t\rfloor,t).
\end{equation}
Using the identities
$$n\psi(q,\frac{t+j}{n})=n\psi(q,\frac{j-i}{n})+q(t+i) \ \ \ \text{ and}$$
$$n\phi(q)+(q-1)(t+j-nq)=n\psi(q-1,\frac{t+j}{n}), \ \ q\in\Z,$$
we can rewrite the above collection of inequalities as
$$y^*_j\ge n\psi(q,\frac{x^*+j}{n})$$
for $j=0,\ldots,n-1$, $q\in\Z$. It remains to apply \eqref{max-phi-eq} again.
\ed

\begin{lem} (i) For each $i$ the section $z_i\in H^0(U_i,L)$ extends uniquely to a global section of $L$.

\noindent
(ii) The sections $(z_i)$ form a basis of $H^0(\TT_{\infty,n},L)$ as $\Z[t]_n$-module.
\end{lem}

\Pf . (i) This is equivalent to showing that $w_i\in M\cap\De$. This is obvious for $w_0=0$. The general
case follows by the $\Z$-action.

\noindent
(ii) The group $H^0(\TT_{\infty,n},L)$ has a $\Z$-basis spanned by $z^u$, where $u\in M\cap\De$.
Note that the multiplication by $t_i$ corresponds to adding $e_i$ to $u$. In particular,
the group $H^0(\TT_{\infty,n},L)$ decomposes into a direct sum of subgroups spanned by $z^u$, where
$u$ has a fixed $f$-component. Thus, it suffices to check that for any $u=-if+m_1e_1+\ldots+m_ne_n\in M\cap\De$
one has $u-w_i\in \Z_{\ge 0}e_1+\ldots+\Z_{\ge}e_n$ (note that $w_i\in M\cap\De$ by part (i)). 
But this follows from the inclusion $C_i^\vee\sub \R\times \R_{\ge 0}^n$, which in turn follows from 
the fact that $C_i$ maps surjectively onto $\R_{\ge 0}^n$ under the projection
$(x,y_1,\ldots,y_n)\mapsto (y_1,\ldots,y_n)$.
\ed

\subsection{Theta functions}

Let us set
$$\theta=\sum_{i\in\Z} z_i.$$
Note that $\theta$ is invariant with respect to the $\Z$-action.
Using the relation \eqref{X-z-rel} and the equivalent relation $z_{i+1}=Y_{i+1}z_i$ we can rewrite this series as
$$\theta=z_0\cdot \bigl(1+[X_0+X_0X_{-1}+X_0X_{-1}X_{-2}+\ldots] + [Y_1+Y_1Y_2+Y_1Y_2Y_3+\ldots]\bigr). 
$$
Recall that $z_0$ gives a trivialization of $L|_{U_0}$ and the functions $(X_i)_{i\le 0}$ and $(Y_i)_{i\ge 1}$
are regular on $U_0$. Furthermore, the relations \eqref{XYt-rel} show that
the infinite sum defining $\th$ becomes finite on any finite order thickening of the central fiber in $U_0$.
By $\Z$-invariance this is true over all $\TT_{\infty,n}$, and so $\th$
defines a global section of $L$ over the formal neighborhood $\hat{\TT}=\hat{\TT}_{\infty,n}$ of the central fiber $\TT_0$.

\begin{lem} The sections $\si_i$ of the central fiber (see \eqref{sigma-i-eq}) extend
to some sections $\si^\th_i:\Spf(\Z[[t]]_n)\to \hat{\TT}$ given by 
$$-X_i=1+(\tau^*)^{-i}s(t)$$ 
on $\hat{\TT}\cap V_i$ for some formal series $s(t)\in\Z[[t]]_n$ with no
constant term, such that $\th$ defines an isomorphism 
\begin{equation}\label{L-si-th-eq}
\OO_{\hat{T}}(\sum_{i\in \Z}\si^\th_i)\simeq L|_{\hat{T}}.
\end{equation}
\end{lem}

\Pf . Restricting to the open subset $V_0$ we can use the trivialization of $L$ given by $z_0$ and
view $\th$ as a function on the formal neighborhood of the central fiber:
$$\th=\sum_{i\in\Z} z^{w_i}=\sum_{i\in\Z} z^{e_i(0)}X_0^{-i}=
1-u-t_0u^{-1}+t_{-1}u^2+t_0^2t_1u^{-2}-t_{-2}t_{-1}^2u^3+\ldots,$$
where $u=-X_0$.
There is a unique formal series $s(t_1,\ldots,t_n)$ with no constant term
such that $u=1+s$ is a zero of the above function. Namely, if we write $s=s_1+s_2+\ldots$,
where $s_d$ is homogeneous of degree $d$ then the equation
$$\sum_{i\in\Z} (-1)^i z^{e_{-i}(0)}(1+s)^i=-s-t_0(1+s)^{-1})+t_{-1}(1+s)^2+\ldots=0$$
will give the recursive formulas for $s_d$. E.g., we get $s_1=t_{-1}-t_0$, $s_2=(t_0+2t_{-1})(t_{-1}-t_0)$, etc.

The fact that $\si^\th_i$ are the only simple zeros of $\th$ follows easily by computing the restriction to the central
fiber:
$$\th|_{\TT_0\cap U_i}=1+X_i+Y_{i+1}.$$
This shows that $\th$ is invertible away from the loci $X_i=-1$ on $\TT_0\cap V_i$,
and that $X_i=-1$ are simple zeros of $\th|_{\TT_0\cap V_i}$.
\ed


We are interested in the formal curve
$$\hat{\TT}_n:=\hat{\TT}/n\Z$$
over $\Z[[t]]_n$. Recall that we have a lifting of the $\Z$-action to $L$, so $L$ descends to a line bundle on
$\hat{\TT}_n$.

For each $i\in\Z/n$ let us set
$$\th_i=\sum_{j\in\Z} z_{i+nj}.$$
These are sections of $L$ over $\hat{\TT}$, invariant with respect to the action of $n\Z\sub\Z$,
so they descend to sections of $L$ on $\hat{\TT}_n$. Similarly, the sections $\si^\th_i$ are $n\Z$-equivariant,
so we can view them as $\Z[[t]]_n$-points of $\hat{\TT}_n$ (with $i\in \Z/n\Z$), and
the isomorphism \eqref{L-si-th-eq} descends to
\begin{equation}\label{L-si-th-eq-bis}
\OO_{\hat{\TT}_n}(\si^\th_0+\ldots+\si^\th_{n-1})\simeq L|_{\hat{\TT}_n}.
\end{equation}

We also have a natural generator $\om$ of the relative dualizing sheaf for the family $\TT_{\infty,n}\to \Spec(\Z[t]_n)$,
such that $\om|_{V_i}=dX_i/X_i$, which induces a generator of the relative dualizing sheaf of
$\hat{\TT}_n$ over $\Z[[t]]_n$. 

\begin{lem}\label{alg-form-lem} The formal scheme $\hat{\TT}_n$ is obtained from the usual 
scheme $T_n$ over $\Z[[t]]_n$, and the line bundle $L$, its sections $\th_i$, the sections $\si^\th_i$, and the element
$\om$, all come from the corresponding data over $T_n$.
\end{lem}

\Pf . The isomorphism \eqref{L-si-th-eq-bis} shows that the
restriction of $L$ to the central fiber, which is the standard $n$-gon, has degree one on every component.
Hence, this restriction is ample. Thus, the assertion follows from Grothendieck's existence theorem
(see \cite[5.1.4, 5.4.5]{EGA3}).
\ed

\begin{defi}\label{n-Tate-def} 
The {\it $n$-Tate curve} is 
the data $(T_n,\si^\th_0,\ldots,\si^\th_{n-1},\om)$ over $\Spec(\Z[[t]]_n)$ defined in Lemma \ref{alg-form-lem}.
\end{defi}

Using the fact that the specialization $t_i=0$ gives 
the $n$-gon curve $G_n$, which a family in $\wt{\UU}_{1,n}^{sns}(\Z)$, 
one can easily deduce that the $n$-Tate curve $(T_n,\si^\th_0,\ldots,\si^\th_{n-1},\om)$ is a family in
$\wt{\UU}_{1,n}^{sns}(\Z[[t]]_n)$. 

Note that by \cite[Prop.\ 1.1.5]{LPol}, for $n\ge 5$ the ring of functions on the moduli space $\wt{\UU}_{1,n}^{sns}$
(which is an affine scheme over $\Z$) is generated by the functions $c_{ij}$ (see \eqref{c-ij-eq}).
Thus, the $n$-Tate curve is determined by the corresponding formal power series $c_{ij}\in \Z[[t]]_n$.
Let us show how to express them in terms of theta functions (this result is not used anywhere else in the paper).

First, we claim that rational functions $\th_i/\th$ have poles only along $\si^\th_i$ and $\si^\th_{i+1}$. Indeed, this is checked
easily by computing the restrictions of these functions to the central fiber: we have
$$\frac{\th_i}{\th}|_{\TT_0\cap V_i}=\frac{1}{X_0+1},$$
$$\frac{\th_i}{\th}|_{\TT_0\cap V_{i+1}}=\frac{X_1}{X_1+1},$$
and the restrictions of $\th_i/\th$ to $\TT_0\cap V_j$ for $j\neq i,i+1$, are zero.

Now we can use our rational functions $\th_i/\th$ to compute the coordinates
$c_{ij}$ from Section \ref{results-sec} for the $n$-Tate curve $T_n$ equipped with the marked points
$p_i=\si^\th_i$ (where $i\in\Z/n$).

Let us denote by $\pa$ the global relative vector field on $\TT_{\infty,n}$ such that
$\pa|_{U_i}=X_i\pa_{X_i}$. Set
$$R_0:=\Res_{p_0}(\frac{\th_0}{\th}\om)=\frac{\th_0 (p_0)}{\pa \th (p_0)}
=\frac{\sum (-1)^{ni} z^{e_{-ni}(0)}(1+s)^{ni}}
{\sum_i (-1)^i i z^{e_{-i}(0)}(1+s)^i}.$$
Note that this is an invertible element of $\Z[[t]]_n$ (equal to $1$ modulo the maximal ideal).
Then the rational function
$$h_{01}:=\frac{\th_0}{R_0 \th}$$
satisfies the conditions of Section \ref{results-sec}: it belongs to $\OO(p_0+p_1)$ and
$h_{01}\om$ has residue $1$ at $p_0$.
Similarly, the rational functions
$$h_{i,i+1}:=\frac{\th_i}{R_i\th},$$
where 
$$R_i:=\frac{\th_i (p_i)}{\pa \th(p_i)}=-\frac{\th_i (p_{i+1})}{\pa \th (p_{i+1})},$$
satisfy similar properties with respect to $p_i$ and $p_{i+1}$.

Now the coordinates $c_{ij}$ on the moduli space are determined by 
$$b_{ij}:=h_{i,i+1}(p_j), \text{  where }\  j\neq i,i+1, \text{ and}$$
$$b_i:=(h_{i-1,i}+h_{i,i+1})|_{p_i}.$$
Namely, since $h_{1i}=h_{12}+h_{23}+\ldots+h_{i-1,i}$, we have
$$c_{ij}=\begin{cases} \sum_{r=1}^{i-1} b_{r,j}, & i<j\\
b_j+\sum_{1\le r<i, r\neq j-1,j} b_{r,j}, & 1<j<i.\end{cases}$$
It remains to express $b_{ij}$ and $b_i$ in terms of the theta functions.

\begin{lem} We have
$$b_{ij}=h_{i,i+1}(p_j)=\frac{\th_i(p_j) \pa\th(p_i)}{\th_i(p_i)\th(p_j)},$$
$$b_{i+1}=(h_{i,i+1}+h_{i+1,i+2})|_{p_{i+1}}=\pa \log \frac{\th_{i+1}}{\th_i} (p_{i+1}),$$
where $j\neq i,i+1$.
\end{lem}

\Pf . The first formula is straightforward. For the second we use l'Hopital's rule:
$$(h_{i,i+1}+h_{i+1,i+2})|_{p_{i+1}}=\frac{R_{i+1}\th_i+R_i\th_{i+1}}{R_iR_{i+1}\th}|_{p_{i+1}}=
\frac{R_{i+1}\pa\th_i+R_i\pa\th_{i+1}}{R_iR_{i+1}\pa\th}|_{p_{i+1}},$$
and then use the definition of $R_i$ and $R_{i+1}$ to rewrite this as
$$-\frac{\pa \th_i(p_{i+1})}{\th_i(p_{i+1})}+\frac{\pa\th_{i+1}(p_{i+1})}{\th_{i+1}(p_{i+1})}.$$
\ed

\subsection{Multiplication of theta functions}\label{mult-theta-sec}

First, let us determine a basis of global sections of $L^m$ for $m\ge 1$ on $\hat{T}$.

For each $p\in\Q$ let us set
$$w_p=-pf+\sum_{i=0}^{n-1}n\phi(\frac{i-p}{n})e_i.$$

\begin{lem}\label{basis-Lm-sections-lem} 
For $m\ge 1$ the sections $(z^{mw_p})_{p\in \frac{1}{m}\Z}$ form a $\Z[t]_n$-basis 
(resp., $\Z[[t]]_n$-basis) of
$H^0(\TT_{\infty,n},L^m)$ (resp., of $H^0(\hat{\TT},L^m)$).
\end{lem}

\Pf . The $\Z$-basis of $H^0(\TT_{\infty,n},L^m)$ is formed by $z^{mu}$ with $u\in\De\cap \frac{1}{m}M$.
Now Lemma \ref{phi-lem}(ii) implies that
$$\De\cap \frac{1}{m}M=\sqcup_{p\in \frac{1}{m}\Z} (w_p+\Z_{\ge 0}e_0+\ldots+\Z_{\ge
0}e_{n-1}),$$
and the assertion follows.
\ed

Note that the $\Z$-action on $L^m$ corresponds to the affine automorphism 
$v\mapsto m\tau(v/m)$ of $m\De\cap M$. This automorphism sends $mw_p$ to $mw_{p-1}$.

Now for $m\ge 1$, $p\in \frac{1}{m}\Z$, let us set
\begin{equation}
\th_{m,p}=\sum_{i\in \Z} z^{mw_{-p+in}}=\sum_{i\in \Z} z^{mn[(\frac{p}{n}+i)f+\phi(\frac{p}{n}+i)e_0+\phi(\frac{p+1}{n}+i)e_1+\ldots+
\phi(\frac{p+n-1}{n}+i)e_{n-1}]}.
\end{equation}
Then $\th_{m,p+n}=\th_{m,p}$ and each $\th_{m,p}$ is invariant with respect to the action of $n\Z\sub \Z$. Furthermore, as an 
easy consequence of Lemma \ref{basis-Lm-sections-lem} we get that
$(\th_{m,p})$, for $p\in \frac{1}{m}\Z/n\Z$, is a $\Z[[t]]_n$-basis of the space of $n\Z$-invariant
sections $H^0(\hat{\TT},L^m)^{n\Z}$. By \cite[5.1.4]{EGA3}, the latter space is identified with $H^0(T_n,L^m)$.

\begin{prop} \label{Bsidecomp} For $m_1,m_2\ge 1$, $p_1\in \frac{1}{m_1}\Z$, $p_2\in \frac{1}{m_2}\Z$, one has 
$$\th_{m_1,p_1}\th_{m_2,p_2}=\sum_{k\in\Z}\th_{m_1+m_2,E(p_1,p_2+kn)}\cdot
\prod_{j=0}^{n-1} t_j^{n\la(\frac{p_1+j}{n},\frac{p_2+j}{n}+k)},$$
where
$$E(a,b)=\frac{m_1a+m_2b}{m_1+m_2},$$
$$\la(a,b)=m_1\phi(a)+m_2\phi(b)-(m_1+m_2)\phi(E(a,b)).$$
\end{prop}

\Pf . The calculation is very similar to the one in \cite[Sec.\ 8.4.2]{Gross}.
\ed

\section{Relative Fukaya categories of genus 1 curves and homological mirror symmetry}
\label{relFuk}

We recall the definition of \emph{the relative Fukaya category} $\mathcal{F}(M,D)$ in
the case where $M=\mathbb{T}$ is a symplectic $2$-torus and $D$ is the divisor consisting of $n$ points on $\T$. The relative Fukaya category of a pair $(M,D)$
was introduced in \cite{SeidelICM} and was further studied in
\cite{Seidelquartic}. We follow  closely the exposition provided in \cite{LP}
where the case $M=\T$ is a symplectic $2$-torus and $D = \{z_1 \}$ is a single
point was discussed in detail. As the construction given there applies to the
mildly generalised situation where $D$ is the union of $n$ marked points, we will not give full details here.

Let $\T$ be a closed, orientable surface of genus 1; $\omega$ a symplectic form on $\T$. Let $z_1,
\ldots, z_n$ be $n$ marked points on $\T$, and $\T_0 = \T \backslash \{ z_1, \ldots, z_n \}$ be the
$n$-punctured torus. We shall fix a primitive $\theta$ for $\omega|_{\T_0}$ and give $\T_0$ a Liouville structure. We
will also fix an unoriented (real) line-field $l$ on $\T$. Such line fields form a torsor
for $C^{\infty}(\T, \mathbb{R}P^1)$ and the connected components can be
identified with $H^1(\T; \mathbb{Z})$. 

One concrete way to fix these data is as follows (cf.\ \cite{HKK}). Let us consider $\T = \mathbb{C}
/ ( \mathbb{Z} \oplus i \mathbb{Z}) $ as
a Riemann surface and $D$ as a divisor on $\T$. Now, consider a holomorphic one-form $\alpha \in
H^0(\T, \Omega_{\T}^{1,0})$. The square $\Omega = \alpha \otimes \alpha \in H^0(\T,
(\Omega_{\T}^{1,0})^{\otimes  2})$ determines a
 non-vanishing quadratic form, which gives a flat Riemannian metric $|\Omega|$
 on $\T$ and a horizontal foliation of tangent vectors $v$ with $\Omega(v,v)
 >0$. The flat Riemannian metric determines an area form $\omega$ and the horizontal
 foliation determines a grading structure on $\T$, i.e., a section of the
 projectivized tangent bundle of $\T$, which we view as an unoriented line field $l   \subset
 T(\T)$. (Note that if one is only interested in working with $\T_0$, one could start with a
 holomorphic one-form $\alpha \in H^0(\T, \Omega^{1,0}_{\T} (D))$ giving rise to a more general grading
 structure on the Fukaya category $\mathcal{F}(\T_0)$, which does not extend to
 $\mathcal{F}(\T,D)$.)
 
 In addition, to be able to work with exact Lagrangians, we also need to fix a primitive
 $\theta$ of  $\omega|_{\T_0}$. This amounts to giving a (real) vector field $Z$ on $\T_0$ which is {\it Liouville}, i.e.,
 satisfies $\mathcal{L}_Z \omega = d(\iota_Z \omega)= \omega$. We can choose this vector field $Z$
 conveniently so as to make our favourite objects given in Figure \ref{figure3} exact, see Prop.
 \ref{exactlags} and Prop.  \ref{primitive} below.
 Recall that a closed curve $L$ is an exact Lagrangian if and only if $\int_L \iota_Z \omega =0$.

Starting from this data, one constructs the relative Fukaya category ${\mathcal
F}(\T,D)$, an $A_\infty$-category linear over $\Z[[t_1,t_2,\ldots,t_n ]]$,
well-defined up to quasi-equivalence (in particular, independent of the primitive $\theta$ and the
line field $l$). 

Recall that the objects of ${ \mathcal F} (\T, D)$ are compact exact Lagrangian
submanifolds $ L \subset \T_0$ which are equipped with an orientation, a spin
structure and a grading (a grading is a homotopy from $l|_{L}$ to $TL$ in $T(\T_0)|_{L}$). Since
$\dim_\R \T_0 = 2$, an oriented Lagrangian submanifold is just
an oriented simple closed curve on $\T_0$. It is well known that an oriented simple closed curve which is not homotopic into a neighborhood of
a puncture (in particular, not null-homotopic), is smoothly isotopic to a 
unique oriented exact Lagrangian up to Hamiltonian isotopy. Furthermore, since we have required
that the grading structure on $\T_0$ is restricted from a grading structure on $\T$, for an
oriented exact Lagrangian in $\T_0$ to have a grading, it is necessary and sufficient that the
underlying curve $L \subset \T_0$ is \emph{non-separating} (see \cite{seidelgraded} for gradings; in
particular the proof of Prop. 2.12 in \cite{seidelgraded} is relevant here). Note also that an oriented exact
Lagrangian can be equipped with either the trivial or the non-trivial spin structure. We often refer to
an object of $\mathcal{F}(\T,D)$ by specifying an exact
Lagrangian $L \subset \T_0$, but suppressing the choice of an orientation, a spin structure and a grading.

The morphism space $hom(L, L')$ is the free $\Z[[t_1,t_2,\ldots, t_n]]$-module
on the intersections $L \cap L'$.  Given a sequence of exact Lagrangians
$L_0,\ldots, L_k$ in $\T_0$, one constructs the $A_\infty$-structure maps \[
\m_k: hom(L_{k-1},L_k) \otimes \ldots \otimes hom(L_0, L_1) \to
hom(L_0,L_k)[2-k] \] defined by counts of solutions to a family of
inhomogeneous Cauchy-Riemann equations $u: S \to (T, L_1 \cup \ldots \cup L_k)$
on the closed unit disk $S$ with $(k+1)$ boundary punctures, which are weighted
by \[ \epsilon(u) t_1^{u \cdot z_1} t_2^{u \cdot z_2} \ldots t_n^{u \cdot z_n}
\] where $u \cdot z_i$ denotes the intersection number of $u$ with $z_i$ and
$\epsilon(u)$ is a sign (see \cite[Sec. 7]{SeidelGenus2} for a formula).

Note that objects of ${ \mathcal F} (\T,D)$ are defined as submanifolds of $\T_0
= \T \backslash D$, however the $A_\infty$-structure is defined by counting maps $u$ that
intersect the divisor $D$. One can define an $A_\infty$-category linear over
$\Z$ where one requires $u \cdot z_i =0 $ for all $i=1,\ldots n$. We
call this \emph{the exact Fukaya category} of the $n$-punctured torus and denote it by
${ \mathcal F}(\T_0)$. Note that by definition
$$\FF(\T_0)=\FF(\T_0,D)\ot_{\Z[[t_1,\ldots,t_n]]} \Z,$$
where we use the homomorphism $\Z[[t_1,\ldots,t_n]]\to\Z$, sending all $t_i$ to zero.
Thus, one should think of ${ \mathcal F}(\T,D)$ as a deformation of ${\mathcal F}(\T_0)$.

We will also write $D^\pi{\mathcal F}(\T,D)$ for the split-closed
triangulated closure of ${\mathcal F}(\T,D)$, which is called the derived Fukaya
category of the pair $(\T,D)$. An explicit model for this triangulated category
is provided by the split-closure of twisted complexes (see \cite{SeidelBook}).

The exact Fukaya category $\mathcal{F}(\T_0)$ has only compact Lagrangians as objects. In fact, there
is an enlargement of this that allows non-compact objects called the \emph{wrapped Fukaya category}
and denoted by $\mathcal{W}(\T_0)$. It is a $\Z$-linear $A_\infty$-category containing $\mathcal{F}(\T_0)$
as a full $A_\infty$-subcategory. The objects of $\mathcal{W}(\T_0)$ are properly embedded
eventually conical exact Lagrangian submanifolds of $\T_0$ equipped with orientations, spin structures and gradings.
The morphism spaces are $hom_{{\WW}(\T_0)}(L,L')$ are cochain complexes $CW^*(L,L')$ computing the
wrapped Floer cohomology (see \cite{AAEKO} for a working definition in this dimension or \cite{abouzgen} for a general definition).  

Finally, we note that a symplectomorphism $\phi : \T_0 \to \T_0$ is exact if $[\phi^* \theta -
\theta] = 0 \in H^1(\T_0)$. Exact symplectomorphisms, equipped with a grading structure, act on $\mathcal{F}(\T_0)$ and
$\mathcal{W}(\T_0)$ as they send exact Lagrangians to exact Lagrangian, and symplectomorphisms act
on Fukaya categories. If $\phi : \T \to \T$ is a symplectomorphism which fixes $D$ pointwise and such that
$\phi|_{\T_0} : \T_0 \to \T_0$ is an exact symplectomorphism, then we also get an action on the
relative Fukaya category $\mathcal{F}(\T,D)$. We also note that exact
symplectomorphisms of $\T_0$ up to exact isotopy form a group, and since
$\dim_{\R}(\T_0)=2$, the natural map to the mapping class group of $\T_0$ yields an
isomorphism by an application of Moser's theorem (cf. \cite{seidelmcg}). Therefore, to compare exact
symplectomorphisms, we may apply techniques from mapping class groups (cf. \cite{primer}).

An example of an exact symplectomorphism equipped with a grading is the (right-handed) 
Dehn twist $\tau_K$ around a (spherical) object $K$ of $\FF(\T_0)$.
We will use the corresponding Dehn twist exact triangle for $L\in\FF(\T,D)$
(see \cite{seidelLES}):
\begin{equation} \label{DTT-eq}
\xymatrix{
    \mathit{HF}^*(K,L) \otimes K \ar[r]^-{\mathrm{ev}} & L \ar[d] \\
                    & \tau_{K}(L). \ar[ul]^{[1]}
    } 
\end{equation}
Note that since $\tau_K$ is an exact symplectomorphism, the Lagrangian $\tau_{K}(L)$ is
exact. Furthermore, the orientation, spin structure and the grading structures on $K$ and $L$,
induces the same structures on $\tau_{K}(L)$ in a canonical fashion for the exact triangle \eqref{DTT-eq} to hold. It is worth
highlighting that if $K$ and $L$ are equippied with non-trivial spin structures, then $\tau_{K}(L)$
should be equipped with the non-trivial spin structure. (See \cite{LPshort} for an explicit verification
of this exact triangle in the case of once-punctured torus.)

\subsection{Generating objects for the relative Fukaya category}\label{generators-Aside-sec}

\begin{figure}[htb!]
\centering
\begin{tikzpicture} [scale=1.1]
      \tikzset{->-/.style={decoration={ markings,
        mark=at position #1 with {\arrow[scale=2,>=stealth]{>}}},postaction={decorate}}}
         \draw (0,0) -- (5,0);
         \draw[red, ->-=.5]  (0,0) -- (6,0);
         \draw [red, ->-=.5] (0,6) -- (0,0);
         \draw [red, ->-=.5] (1,6) -- (1,0);
         \draw [red, ->-=.5] (4,6) -- (4,0);
         \draw [red, ->-=.5] (5,6) -- (5,0);
         
         \draw (0,6) -- (6,6);
         \draw (6,0) -- (6,6);
       
         \node at (0,5) {$\star$};
          \node at (1,5) {$\star$};
          \node at (4,5) {$\star$};
           \node at (5,5) {$\star$};
           \node at (0.5,0) {$\star$};

\node at (-0.3,2)   {\footnotesize $L_n$};
         \node at (0.65,2)   {\footnotesize $L_{n-1}$};
         \node at (3.7,2)   {\footnotesize $L_{2}$};
         \node at (4.7,2)   {\footnotesize $L_1$};
         \node at (3,-0.3)   {\footnotesize $L_0$};

\draw[red, thick, fill=red] (1.7,2) circle(.02);
\draw[red, thick, fill=red] (2.4,2) circle(.02);
\draw[red, thick, fill=red] (3.1,2) circle(.02);

\draw[black, thick, fill=black] (2.2,0.5) circle(.01);
\draw[black, thick, fill=black] (2.6,0.5) circle(.01);
\draw[black, thick, fill=black] (3.0,0.5) circle(.01);

\draw[thick, fill=black] (0.5,0.5) circle(.03);
\draw[thick, fill=black] (1.5,0.5) circle(.03);
\draw[thick, fill=black] (3.5,0.5) circle(.03);
\draw[thick, fill=black] (4.5,0.5) circle(.03);
\draw[thick, fill=black] (5.5,0.5) circle(.03);

\node at (0.5,0.7) {\footnotesize $z_n$};
\node at (1.5,0.7) {\footnotesize $z_{n-1}$};
\node at (3.5,0.7) {\footnotesize $z_{3}$};
\node at (4.5,0.7) {\footnotesize $z_2$};
\node at (5.5,0.7) {\footnotesize $z_1$};

\end{tikzpicture}
\caption{Symplectic torus with $n+1$ oriented Lagrangians, $n$ vertical, $1$ horizontal}
\label{figure2}
\end{figure}
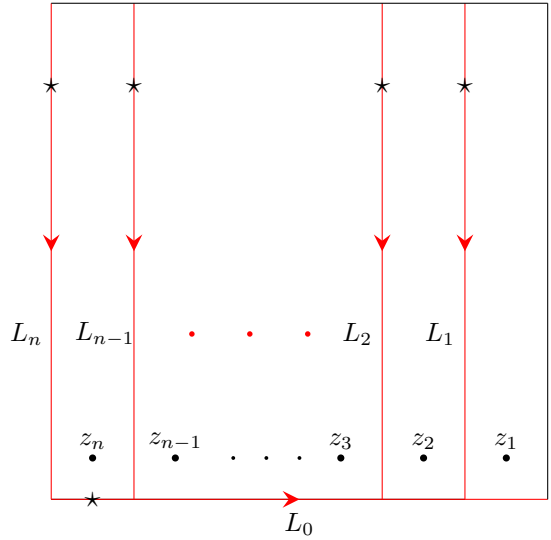

    We consider $(n+1)$ objects 
    \[ L_0, L_1, \ldots L_n    \in
    { \mathcal F}(\T,D)\] with the underlying oriented
    Lagrangians  given by exact representatives of simple closed curves
    that are depicted in Figure \ref{figure2}. In Cartesian coordinates $(x,y) \in \R^2/\Z^2$, we
    have \[ L_0 = \{ (x,0) : x \in \R/\Z \}, \ \  L_i = \{ (n-i)/n , y) : y \in \R/\Z \} \text{
    for $i =1,\ldots n$}. \]
We equip these with non-trivial spin structures
    which we keep track of with a marked point $\star \in L_i$ signifying the non-trivial double
    cover of $L_i$. These are needed in the calculation of signs $\epsilon(u)$ for polygons $u$ that
    contribute to the $A_\infty$-structure maps. Namely, if the intersection points at the corners
    of $u$ have \emph{even} Floer cohomology indices then $\epsilon(u) = (-1)^s$ where $s$ is the
    number of stars on the boundary. In our explicit computations, we will be in this situation. In
    general, $\epsilon(u)$ depends on orientations, spin structures and the indices of the corners.  For a complete
    description of gradings and of signs $\epsilon(u)$, see \cite[Sec. 7]{SeidelGenus2}, and also
    \cite{LP}.

    As depicted in Figure \ref{figure2}, we specify the divisor $D = \{z_1,\ldots,z_n \}$ by
    letting
    \[ \{z_1,z_2,\ldots, z_n \} := \{ (n-i)/n + \epsilon, \epsilon) : i=1,\ldots, n \} \] 
    in $\mathbb{R}^2/\mathbb{Z}^2$, for sufficiently small $\epsilon>0$. Ultimately, the exact locations of the points $z_i$ do not matter as long as there is one point on
each connected component of $\T \setminus \{ L_0,L_1,\ldots, L_n \}$. The specific choice that we
have (for sufficiently small $\epsilon$) is so that the formulae that we will get in explicit
calculations would match the identity of Prop. \ref{Bsidecomp}.

\begin{lem} \label{relgenerate} The objects $L_0, L_1, \ldots L_n$ split-generate the split-closed triangulated
    category $D^\pi{\mathcal F}(\T,D)$.
\end{lem}
\Pf.  It is well known that Dehn twists around the curves $L_0,\ldots,L_n$ generate the pure mapping
class group of $(\mathbb{T},D)$ (\cite[Sec. 4.4.4]{primer}). Furthermore, the
pure mapping
class group acts transitively on the set of oriented non-separating simple closed curves, as
follows from the classification of surfaces (\cite[Sec. 1.3.1]{primer}). Hence, using the exact triangles
associated with Dehn twists (see \eqref{DTT-eq}), we
deduce that any object $L$ of $\mathcal{F}(\T,D)$, where the underlying spin structure is
non-trivial, is generated by the collection $L_0,L_1,\ldots, L_n$. Now, given an arbitrary object
$L$ with a trivial spin structure, we claim that $L \oplus L[2]$ is generated by $L_0,L_1,\ldots,
L_n$. To see this, let $L'$ be the object with the same underlying oriented curve and same grading structure as $L$, but with
the non-trivial spin-structure. Let $L^{''}$ be an object whose underlying curve intersect $L'$
at a unique point, and also equipped with a non-trivial spin structure. Now $(L',L'')$ is an
$(A_2)$-configuration in $\T_0$. Hence, a neighborhood of
$L' \cup L^{''}$ is a torus $P$ with one boundary component $\partial P$, embedded in $\T_0$. Now, looking at the mapping class
group of $P$, we have the relation \[ (\tau_{L'} \tau_{L^{''}})^6 \simeq \tau_{\partial P}. \]
In particular, we can observe that $(\tau_{L'} \tau_{L^{''}})^6$ sends the curve $L$ back to itself.
However, as proven in Lemma 5.9 of \cite{seidelgraded}, this automorphism acts non-trivially on the
grading. Namely, we have 
\[ (\tau_{L'} \tau_{L^{''}})^6 (L) \simeq L[2]. \]
By the above argument, we know that $L'$ and $L^{''}$ are generated by the collection
$L_0,L_1,\ldots, L_n$. Hence, we can combine the exact triangles of Dehn twists to obtain an exact triangle
between $L$, $(\tau_{L'} \tau_{L^{''}})^6(L) \simeq L[2]$ and a complex built out of $L'$ and
$L^{''}$. Now, for grading reasons, the map between $L$ and $L[2]$ has to vanish from which we conclude that $L \oplus L[2]$ is
generated by $L_0, L_1,\ldots,L_n$, hence $L$ is split-generated by $L_0, L_1 \ldots,L_n$.
\ed

In view of this generation result, one studies the derived Fukaya category of $(\T,  D)$ via the
$A_\infty$-algebra over $\Z[[t_1,t_2,\ldots t_n]]$,
\[ {\mathscr A} = \bigoplus_{i,j=0}^{n} hom_{{\mathcal F}(\T,D)}(L_i,L_j). \]
We also consider the exact Fukaya category ${\mathcal F}(\T_0)$ and correspondingly, we have the $A_\infty$-algebra 
over $\Z$,
\begin{equation}\label{A0-algebra-eq} 
{\mathscr A}_0 = \bigoplus_{i,j=0}^{n} hom_{{\mathcal F}(\T_0)} (L_i,L_j). 
\end{equation}

Note that the proof in Lemma \ref{relgenerate} also gives that the collection 
$(L_0,\ldots,L_n)$ split-generates $\mathcal{F}(\T_0)$. 

\begin{figure}[htb!]
\centering
\begin{tikzpicture} [scale=1.1]
      \tikzset{->-/.style={decoration={ markings,
        mark=at position #1 with {\arrow[scale=2,>=stealth]{>}}},postaction={decorate}}}
         \draw (0,0) -- (5,0);
         \draw[red, ->-=.5]  (0,0) -- (6,0);
         \draw [red, ->-=.5] (0,6) -- (0,0);
         \draw [red, ->-=.5] (1,6) -- (1,0);
         \draw [red, ->-=.5] (4,6) -- (4,0);
         \draw [red, ->-=.5] (5,6) -- (5,0);

         \draw [blue, ->-=.5] (1.5,0.5) -- (0.5,0.5);
         \draw [blue, ->-=.5] (4.5,0.5) -- (3.5,0.5);
         \draw [blue, ->-=.5] (5.5,0.5) -- (4.5,0.5);
         \draw [blue, ->-=1] (0.5,0.5) -- (0,0.5);
         \draw [blue ] (6,0.5) -- (5.5,0.5);
         \draw [blue, ->-=.5] (5.5,6) -- (5.5,0);

         \draw (0,6) -- (6,6);
         \draw (6,0) -- (6,6);
        
         \node at (-0.3,2)   {\footnotesize $L_n$};
         \node at (0.65,2)   {\footnotesize $L_{n-1}$};
         \node at (3.7,2)   {\footnotesize $L_{2}$};
         \node at (4.7,2)   {\footnotesize $L_1$};
         \node at (3,-0.3)   {\footnotesize $L_0$};

         \node at (-0.2,0.4)   {\footnotesize $\hat{L}_n$};
         \node at (1,0.25)   {\footnotesize $\hat{L}_{n-1}$};
         \node at (3.8,0.25)   {\footnotesize $\hat{L}_{2}$};
         \node at (4.8,0.25)   {\footnotesize $\hat{L}_{1}$};
         \node at (5.7,2)   {\footnotesize $\hat{L}_0$};

         \node at (0,5) {$\star$};
          \node at (1,5) {$\star$};
          \node at (4,5) {$\star$};
           \node at (5,5) {$\star$};
           \node at (0.5,0) {$\star$};

\draw[red, thick, fill=red] (1.7,2) circle(.02);
\draw[red, thick, fill=red] (2.4,2) circle(.02);
\draw[red, thick, fill=red] (3.1,2) circle(.02);

\draw[black, thick, fill=black] (2.2,0.5) circle(.01);
\draw[black, thick, fill=black] (2.6,0.5) circle(.01);
\draw[black, thick, fill=black] (3.0,0.5) circle(.01);

\draw[thick, fill=black] (0.5,0.5) circle(.03);
\draw[thick, fill=black] (1.5,0.5) circle(.03);
\draw[thick, fill=black] (3.5,0.5) circle(.03);
\draw[thick, fill=black] (4.5,0.5) circle(.03);
\draw[thick, fill=black] (5.5,0.5) circle(.03);

\node at (0.5,0.7) {\footnotesize $z_n$};
\node at (1.5,0.7) {\footnotesize $z_{n-1}$};
\node at (3.5,0.7) {\footnotesize $z_{3}$};
\node at (4.5,0.7) {\footnotesize $z_2$};
\node at (5.5,0.7) {\footnotesize $z_1$};
\end{tikzpicture}
\caption{$n+1$ non-compact Lagrangians (blue) dual to $n+1$ compact Lagrangians (red).}
\label{figure3}
\end{figure}
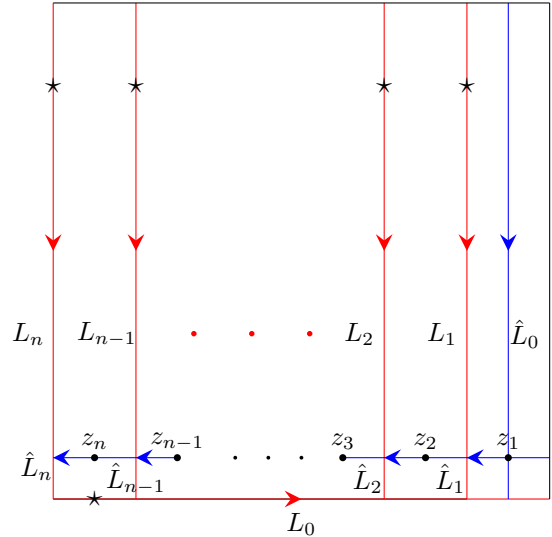

As for the wrapped Fukaya category, we will consider dual objects. Namely, we consider $(n+1)$ objects 
\[ \hat{L}_0, \hat{L}_1, \ldots \hat{L}_n     \in{
\mathcal W}(\T_0)\] with
    the underlying oriented non-compact arcs given by exact
    representatives of free homotopy classes that are depicted
    in Figure \ref{figure3}. In Cartesian coordinates $(x,y) \in \mathbb{R}^2/ \Z^2$, we have that 
    \[ \hat{L}_0 = \{ ((n-1)/n + \epsilon, \epsilon + t) : t \in (0,1) + \Z \}, \] \[ \hat{L}_i =
    \{ ( \epsilon + (n-i-t)/n, \epsilon) : t \in (0,1) + \Z \} \text{ for } i=1,\ldots, n. \] 
Note that there exists a unique isomorphism class of a spin structure on these arcs once they are oriented.

\begin{lem} \label{ncgener} The objects $\hat{L}_0, \hat{L}_1, \ldots \hat{L}_n$
    generate the triangulated category $D^b {\mathcal W}(\T_0)$.
\end{lem}
\Pf. We shall apply the argument of Theorem A.1 of \cite{AAEKO} (which in turn is based on
Prop.\ 18.17 of \cite{SeidelBook}). Namely, the once-punctured torus has a Lefschetz fibration over
$\C$ given by a double covering branched over 3 points. There is an unbranched covering of the
once-punctured torus by the $n$-punctured torus. This, in turn, gives a Lefschetz fibration on the $n$-punctured torus. Now Theorem A.1 of \cite{AAEKO} gives that a basis of
thimbles for this Lefschetz fibration generates the derived category $D^b \mathcal{W}(C)$. 

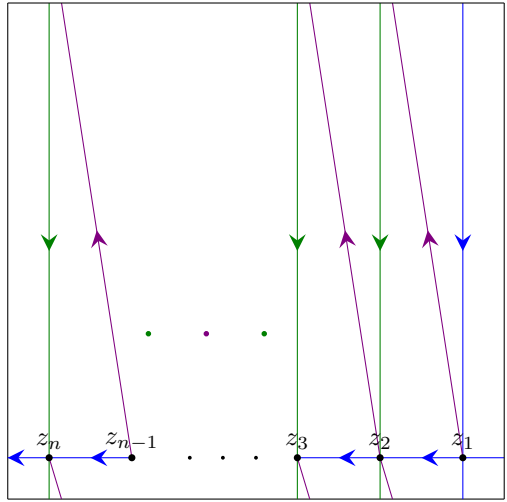
\begin{figure}[htb!]
\centering
\begin{tikzpicture} [scale=1.1]
      \tikzset{->-/.style={decoration={ markings,
        mark=at position #1 with {\arrow[scale=2,>=stealth]{>}}},postaction={decorate}}}
         \draw (0,0) -- (6,0);
         \draw (0,0) -- (0,6);

         \draw [blue, ->-=.5] (1.5,0.5) -- (0.5,0.5);
         \draw [blue, ->-=.5] (4.5,0.5) -- (3.5,0.5);
         \draw [blue, ->-=.5] (5.5,0.5) -- (4.5,0.5);
         \draw [blue, ->-=1] (0.5,0.5) -- (0,0.5);
         \draw [blue ] (6,0.5) -- (5.5,0.5);
         \draw [blue, ->-=.5] (5.5,6) -- (5.5,0);

         \draw [green!50!black, ->-=.5] (0.5,6) -- (0.5,0);
         \draw [green!50!black, ->-=.5] (3.5,6) -- (3.5,0);
         \draw [green!50!black, ->-=.5] (4.5,6) -- (4.5,0);

         \draw [violet, ->-=.5] (5.5,0.5) -- (4.65,6);
         \draw [violet, =.5] (4.65,0) -- (4.5,0.5);
         \draw [violet, ->-=.5] (4.5,0.5) -- (3.65,6);
         \draw [violet, =.5] (3.65,0) -- (3.5,0.5);
         \draw [violet, ->-=.5] (1.5,0.5) -- (0.65,6);
         \draw [violet, =.5] (0.65,0) -- (0.5,0.5);

         \draw (0,6) -- (6,6);
         \draw (6,0) -- (6,6);

\draw[green!50!black, thick, fill=green!50!black] (1.7,2) circle(.02);
\draw[violet, thick, fill=violet] (2.4,2) circle(.02);
\draw[green!50!black, thick, fill=green!50!black] (3.1,2) circle(.02);

\draw[black, thick, fill=black] (2.2,0.5) circle(.01);
\draw[black, thick, fill=black] (2.6,0.5) circle(.01);
\draw[black, thick, fill=black] (3.0,0.5) circle(.01);

\draw[thick, fill=black] (0.5,0.5) circle(.03);
\draw[thick, fill=black] (1.5,0.5) circle(.03);
\draw[thick, fill=black] (3.5,0.5) circle(.03);
\draw[thick, fill=black] (4.5,0.5) circle(.03);
\draw[thick, fill=black] (5.5,0.5) circle(.03);

\node at (0.5,0.7) {\footnotesize $z_n$};
\node at (1.5,0.7) {\footnotesize $z_{n-1}$};
\node at (3.5,0.7) {\footnotesize $z_{3}$};
\node at (4.5,0.7) {\footnotesize $z_2$};
\node at (5.5,0.7) {\footnotesize $z_1$};
\end{tikzpicture}
\caption{Generators given by a basis of thimbles} 
\label{figure3bis}
\end{figure}

In our case, a choice of such a basis of thimbles is depicted in Figure \ref{figure3bis}. This
construction gives $3n$ objects which is more than what one needs. Indeed, we will use explicit exact sequences to show that the
(blue) curves that correspond to $\hat{L}_i$, for $i=0,\ldots,n$, generate the other ones. For this purpose, let us
label the vertical (green for $i>1$) curves as $G_i$, for $i=1,\ldots,n$, where $G_i$ has ends approaching to
the puncture $z_i$. Note that $G_1 = \hat{L}_0$. Let us also label the negatively sloped (violet)
curves, that connect $z_i$ to $z_{i+1}$, as $P_i$, for $i=1,\ldots, n$. We claim that we have the
following exact triangles for $i=1,\ldots,n$ (where $G_{n+1}=G_1$): 
\begin{equation} \label{wexact}
\xymatrix{
    \hat{L}_i \ar[r]  &  P_i     \ar[d] \\
                    & G_i
 \ar[ul]^{[1]}
    } 
\ \ \ , \ \ \ \xymatrix{
    \hat{L}_i \ar[r]  &  P_i     \ar[d] \\
                      & G_{i+1} 
 \ar[ul]^{[1]}
    } 
\end{equation}
(Note that the objects on the top of the two exact triangles are the same, however, the
degree 0 morphisms between them are different. Compare these to the exact sequences \eqref{O(-1)-O-qi-seq} and 
\eqref{O(-1)-O-qi-bis-seq}
appearing later in the proof of Prop. \ref{Db-Coh-gen-prop}). 

Since $G_1= \hat{L}_0$ by definition, the above exact sequences suffice to show that $\{ \hat{L}_i
\}$, for $i=0,\ldots, n$, generate all the objects $P_i$ and $G_i$, and hence, the wrapped category
$D^b \mathcal{W}(\T_0)$. 

It remains to establish the existence of the exact triangles \eqref{wexact}. We will give the proof for the one on the left.
The proof for the other one is similar. Let $\Sigma_i$ be a Liouville subdomain in $\mathbb{T}_0$ that is obtained by removing a
collar neighborhood of the puncture at $z_i$. In other words, $\mathbb{T}_0$ can be symplectically
identified as the completion $\widehat{\Sigma}_i$ of $\Sigma_i$ at its compact boundary (by putting back the collar neighborhood). Our strategy will be to use the restriction functors on wrapped Fukaya categories constructed by
Abouzaid and Seidel in \cite{AbouzSeidel}, which in our case gives an exact functor
\[ \mathcal{W}(\widehat{\Sigma}_i) \to \mathcal{W}(\T_0) \]
such that on objects we just intersect the underlying Lagrangian in $\widehat{\Sigma}_i$ with $\Sigma_i$ and then extend it to $\T_0$ in the obvious conical way. We will establish an exact triangle in $\widehat{\Sigma}_i$ associated to a (negative) Dehn twist
around a spherical object and the desired exact triangle in $\mathcal{W}(\T_0)$ will be obtained as
the image under this restriction functor.

\begin{figure}[htb!]
\centering
\begin{tikzpicture} [scale=1.1]
      \tikzset{->-/.style={decoration={ markings,
        mark=at position #1 with {\arrow[scale=2,>=stealth]{>}}},postaction={decorate}}}

\draw[thick] (0,0) circle(1);
\draw [green!50!black, ->-=.7] (0,1) to[in=160,out=200] (0,-1);
\draw [blue,  ->-=.8] (0,0) to (-1,0);
\draw [violet,  ->-=.4] (0,0) to (-0.7,0.7);
\draw[black, thick, fill=black] (0,0) circle(.03);

\draw[thick] (3.5,0) circle(1);
\draw [green!50!black, ->-=.7] (3.5,1) to (3.5,-1);
\draw [blue,  ->-=.8] (3.5,0) to (3.5-1,0);
\draw [violet,  ->-=.4] (3.5,0) to (3.5-0.7,0.7);
\draw[black, thick, fill=black] (3.5,0) circle(.03);
\draw[thick] (0,0) circle(1);

\end{tikzpicture}
\caption{Neighborhood of the puncture $z_i$ and the extensions of $\hat{L}_i$, $P_i$ and $G_i$
in $\widehat{\Sigma}_i$ (left) and $\mathbb{T}_0$ (right). } 
\label{figure3bisbis}
\end{figure}
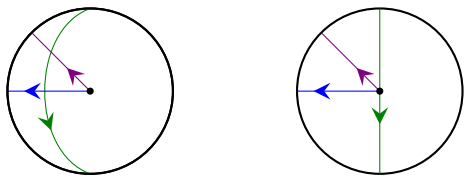

Now, we have the restriction of the objects $\hat{L}_i$, $P_i$ and $G_i$ to $\Sigma_i$. We
shall extend them to $\widehat{\Sigma}_i$ so that $\hat{L}_i$ and $P_i$ are extended as in
$\mathbb{T}_0$ but we modify $G_i$ so that it is extended to a compact curve $\overline{G}_i
\subset \widehat{\Sigma}_i$, see Figure \ref{figure3bisbis} where the neighborhood of the puncture
$z_i$ is drawn. In order to ensure that $\overline{G}_i$ is actually an exact Lagrangian in
$\widehat{\Sigma}_i$, one may need to isotope $G_i$ in $\Sigma_i$ but this is unproblematic.
(Note that in Figure \ref{figure3bisbis}, we completed $G_i$ to a compact curve that goes around
the puncture in a counter-clockwise manner, the other choice is used in proving the 
other exact triangle in (\ref{wexact})). 

On $\widehat{\Sigma}_i$ the Lagrangian $\overline{G}_i$ is a spherical object, so we can apply
the Dehn twist exact triangle associated to the (negative) Dehn twist around $\overline{G}_i$.
This gives an exact triangle
\begin{equation}
\xymatrix{
    \hat{L}_i \ar[r]  &  P_i     \ar[d] \\
                      & \overline{G}_i
 \ar[ul]^{[1]}
    } 
\end{equation}
since it is easy to see that $\tau_{\overline{G}_i}^{-1}(P_i)$ is isotopic to $\hat{L}_i$ in $\widehat{\Sigma}_i$. 

Finally, note that we have arranged it so that after restricting $\hat{L}_i, P_i$ and
$\overline{G}_i$ to $\Sigma_i$ and extending them conically to $\mathbb{T}_0$, we get back the
objects $\hat{L}_i, P_i$ and $G_i$ in $\mathcal{W}(\T_0)$ (this is clear from Figure \ref{figure3bisbis}). 
Hence, the image of the above exact
triangle in $\mathcal{W}(\T_0)$ under the restriction functor is 
\begin{equation}\label{LPG-exact-triangle-eq}
\xymatrix{
    \hat{L}_i \ar[r]  &  P_i     \ar[d] \\
                      & G_i
 \ar[ul]^{[1]}
    } 
\end{equation}
which is precisely the left triangle in \eqref{wexact}.
\ed

\begin{rem} \label{pitfall} An alternative approach would be to try to check the split-generation result due to  Abouzaid
    \cite[Thm. 1.1]{abouzgen}. In order to apply this result, one needs to show that the open-closed map:\[ \mathcal{OC}:    HH_{*-1}( \langle \hat{L}_i
\rangle ) \to SH^*(\T_0) \] hits the unit in $SH^0(\T_0)$ (we shifted the grading so that
$\mathcal{OC}$ is degree preserving). Now, the unit in $SH^0(\T_0)$ is represented by the
fundamental class. Hence, it is enough to find holomorphic polygons with boundary on the $\{
\hat{L}_i \}$ whose total image covers the whole manifold $\T_0$ with multiplicity 1. At first sight, this seems easy to arrange
since the complement $\T_0 \setminus \{ \hat{L}_0, \hat{L}_1,\ldots, \hat{L}_n \}$ consists of a
disjoint union of (open) polygons. However, even though the Hochschild chain that appears at the
corners of these polygons is sent by the open-closed map to the unit in $SH^0(\T_0)$, there is
no a priori guarantee that this chain is actually a cocycle and so it may not give an element of
$HH_*(\langle \hat{L}_i \rangle)$. We thank Hansol Hong for warning us about this dangerous
pitfall. 
\end{rem} 

In view of Lemma \ref{ncgener} one can study the derived Fukaya category of $\WW(\T_0)$ via the
$A_\infty$-algebra over $\Z$,
\begin{equation}\label{B-algebra-eq}
 {\mathscr B} = \bigoplus_{i,j=0}^{n} hom_{{\mathcal W}(\T_0)}(\hat{L}_i,\hat{L}_j). 
 \end{equation}

We note that the intersection pattern of the compact Lagrangians $L_i$ with the non-compact
Lagrangians $\hat{L}_i$ immediately gives the following duality property.

\begin{prop} \label{exactlags} One has 
    $$\Hom_{\WW(\T_0)}(L_i,\hat{L}_j)=\Hom_{\WW(\T_0)}(\hat{L}_j,L_i)=0 \text{ for } i\neq j,$$
    $$\Hom_{\WW(\T_0)}(L_i,\hat{L}_i)\simeq \Hom_{\WW(\T_0)}(\hat{L}_i,L_i)[1]\simeq \Z$$ 
for $i=0,\ldots,n$. Furthermore, the composition map
$$\Hom_{\WW(\T_0)}(\hat{L}_i,L_i)\ot \Hom_{\WW(\T_0)}(L_i, \hat{L}_i)[1]\to \Hom_{\FF(\T_0)}(L_i,L_i)[1]$$
is an isomorphism.
\end{prop} 
\Pf .  We only need to check that we can actually find a primitive $\theta$ of $\omega|_{\T_0}$
that makes the Lagrangians $L_i$ and $\hat{L}_i$, as depicted in Figure \ref{figure3},
exact. This is equivalent to exhibiting a Liouville vector field $Z$ on $\T_0$ such that $\int_{L_i} \iota_Z \omega = \int_{\hat{L}_i} \iota_Z \omega =0$ for all $i=0,\ldots, n$. By observing that $\T_0$ retracts
onto a neighborhood of $\bigcup_i L_i$ which can be locally identified with a neighborhood of
a plumbing of $T^*L_i$'s, it is easy to see that there is a vector
field $Z$ which vanishes along $L_i$ and is tangent to $\hat{L}_i$. Hence, the Lagrangians
$L_i$ and $\hat{L}_i$ can be made exact as drawn in Figure \ref{figure3}. 
\ed

This duality property will play a key role in proving that $\mathscr{A}_0$ and $\mathscr{B}$, equipped with certain
augmentations, are Koszul dual $A_\infty$-algebras (see Sec.\ \ref{koszul}). 

Let $T : \T \to \T$ be the symplectomorphism given by the translation 
\begin{equation}\label{T-translation-eq} 
T(x,y) = (x+1/n, y) 
\end{equation}

Then, $T$ preserves $L_0$, sends $L_i$ to $L_{i+1}$ for $i=1,\ldots (n-1)$, $L_n$ to $L_1$, and
preserves $D = \{ z_1,\ldots,z_n\}$. Since the homology classes of $L_i$ for $i=0,\ldots,n$ give a basis of $H_1(\T_0)$, 
we derive that for a primitive $\theta$ of $\omega|_{\T_0}$ for which $L_i$ are exact
Lagrangians, the induced symplectomorphism $T : \T_0 \to \T_0$ is exact, i.e., $[T^*\theta - \theta] = 0 \in
H^1(\T_0)$. In
particular, $T$ acts on the wrapped Fukaya category $\WW(\T_0)$.  

\begin{figure}[htb!]
\centering
\begin{tikzpicture} [scale=1.1]
      \tikzset{->-/.style={decoration={ markings,
        mark=at position #1 with {\arrow[scale=2,>=stealth]{>}}},postaction={decorate}}}
        \draw [red, ->-=.5] (0,0) -- (6,0);
         \draw (0,6) -- (0,0);
         \draw [red, ->-=.5] (0.5,6) -- (0.5,0);
         \draw [red, ->-=.5] (3.5,6) -- (3.5,0);
         \draw [red, ->-=.5] (4.5,6) -- (4.5,0);
         \draw [red, ->-=.5] (5.5,6) -- (5.5,0);

         \draw (0,6) -- (6,6);
         \draw (6,0) -- (6,6);
        
         \node at (1.1,2)   {\footnotesize $T^{n-1} \hat{L}_0$};
         \node at (4.85,2)   {\footnotesize $T \hat{L}_0$};
         \node at (3.95,2)   {\footnotesize $T^2 \hat{L}_0$};

         \node at (3,-0.3)   {\footnotesize $L_0$};

         \node at (5.7,2)   {\footnotesize $\hat{L}_0$};

           \node at (0.3,0) {$\star$};

\draw[red, thick, fill=red] (1.8,2) circle(.02);
\draw[red, thick, fill=red] (2.5,2) circle(.02);
\draw[red, thick, fill=red] (3.2,2) circle(.02);

\draw[black, thick, fill=black] (1.2,0.5) circle(.01);
\draw[black, thick, fill=black] (2.1,0.5) circle(.01);
\draw[black, thick, fill=black] (3.0,0.5) circle(.01);

\draw[thick, fill=black] (0.5,0.5) circle(.03);
\draw[thick, fill=black] (3.5,0.5) circle(.03);
\draw[thick, fill=black] (4.5,0.5) circle(.03);
\draw[thick, fill=black] (5.5,0.5) circle(.03);

\node at (0.5,0.7) {\footnotesize $z_n$};
\node at (3.5,0.7) {\footnotesize $z_{3}$};
\node at (4.5,0.7) {\footnotesize $z_2$};
\node at (5.5,0.7) {\footnotesize $z_1$};

\end{tikzpicture}
\caption{Another set of generators for $\WW(\T_0)$} \label{figure4}
\end{figure}
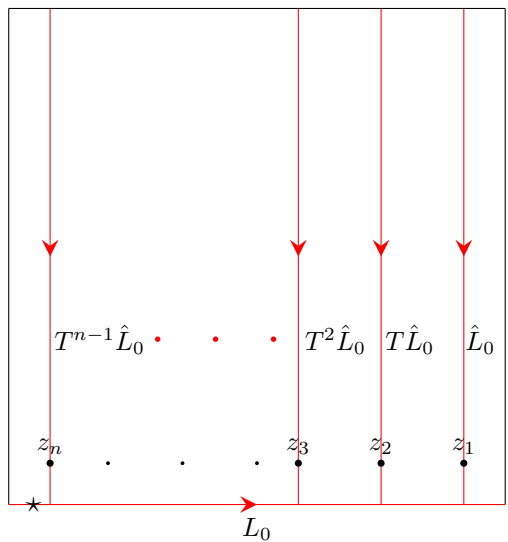

Another set of generators for the wrapped category is given in the following lemma (see Figure \ref{figure4}).  

\begin{lem}\label{another-generators-W-lem} $D^\pi(\WW(\T_0))$ is split-generated by $L_0$ and the objects $T^i\hat{L}_0$, $i=0,\ldots,n-1$.
\end{lem}

\Pf. The proof is similar to the proof of Lemma \ref{ncgener}. In the notation of the proof of
Lemma \ref{ncgener}, we have $T^i \hat{L}_0 = G_{i+1}$ for $i=0,1\ldots, n-1$. 

Note that for each $i$ we have that $G_i$ intersects $L_0$ at a unique point. As in the proof
of Lemma \ref{ncgener}, we can consider compact curves $\overline{G}_i \subset
\widehat{\Sigma}_i$. Let us do this for all $i=1,\ldots, n$, and write
$\widehat{\Sigma}$ for the completion obtained in this manner for all $i$ simultaneously.
Next, consider the Dehn twist exact triangle corresponding to
the composition $\tau = \tau_{\overline{G}_1} \circ \ldots \circ
\tau_{\overline{G}_n}$ around the disjoint curves $\overline{G}_i$ for $i=1, \ldots n$. This gives the following exact triangle in $\mathcal{W}(\widehat{\Sigma})$:
\begin{equation}
\xymatrix{
    \bigoplus_{i=1}^n \overline{G}_i \ar[r]  &  L_0[1] \ar[d] \\
                                              & \bigoplus_{i=1}^n \overline{P}_i \
 \ar[ul]^{[1]}
    } 
\end{equation}
where $\overline{P}_i$ are compact curves which restrict to the non-compact curves $P_i[1]$ in $\Sigma$
(note that $P_i[1]$ has the same underlying non-compact Lagrangian as $P_i$ but equipped with the opposite
orientation) .
Hence, as in Lemma \ref{ncgener}, we can use the exact restriction functor
$\mathcal{W}(\widehat{\Sigma}) \to \mathcal{W}(\mathcal{\T}_0)$ to obtain an exact triangle (which is an analog
of the exact sequence \eqref{surgery} appearing later on the B-side):
\begin{equation}
\xymatrix{
    L_0 \ar[r]  &  \bigoplus_{i=1}^n P_i \ar[d] \\
                & \bigoplus_{i=1}^n G_i \
 \ar[ul]^{[1]}
    } 
\end{equation}
(Note that one could alternatively use the Lagrangian surgery exact triangle due to
Fukaya-Oh-Ohta-Ono \cite{FOOO} to obtain this result more directly.)  

This proves that $L_0$ together with $G_i$, for $i=1,\ldots n$, split-generate $P_i$. Now the exact triangles
\eqref{LPG-exact-triangle-eq} show that $G_i$ and $P_i$ generate $\hat{L}_i$. Since by Lemma \ref{ncgener}, $\hat{L}_i$,
for $i=0,\ldots, n$, generate $\mathcal{W}(\T_0)$, the result follows.
\ed

Additively, wrapped Floer cohomology of $T^i\hat{L}_0$ is easy to compute, as it is easy to see
that the Floer differential vanishes.

\begin{lem}\label{sympl-dim-Hom-lem} For all $i=0,\ldots, n-1$, the non-compact Lagrangians
    $T^i\hat{L}_0$ can be equipped with a grading structure such that for any commutative ring $R$
    one has an isomorphism of $R$-modules 
    \[ \Hom^*_{\WW(\T_0)\ot R}(T^i\hat{L}_0,T^i\hat{L}_0) \cong R \langle u_i, v_i \rangle / (u_i^2, v_i^2) \]
    where $\deg(u_i)= \deg(v_i)= 1$. 
\end{lem}
\Pf. Recall that the wrapped Floer complex $CW^*(T^i \hat{L}_0, T^i \hat{L}_0 )$ is generated by time 1
chords of the flow generated by a Hamiltonian $H : \T_0 \to \R$ which is quadratic at infinity (in
other words, $H(r,\theta) = r^2$ at the cylindrical ends; this implies that near the ends the flow is given by the
clockwise rotation). For simplicity, we require that $H|_{T^i\hat{L}_0}$ is a Morse function with a unique minimum. The minimum gives a generator of
$\Hom^0 (T^i\hat{L}_0, T^i\hat{L}_0)$. Let us label the shortest (non-constant) time 1 chords
$u_i$ and $v_i$ (corresponding to left and right semicircles at the cylindrical end oriented clockwise). Since we required the line field $l$ on $\T_0$ to extend to $\T$,
the rotation number of a small simple closed loop around $i^{th}$ puncture with respect to the
induced trivialization is 1. It follows that the Maslov indices of $u_i$ and $v_i$ have to obey
the equality
\[ \deg(u_i) + \deg(v_i) = 2, \]
as composing $u_i$ and $v_i$, we get a loop that goes once around the puncture. In addition, because of the rigid polygonal region with the boundary $T^i \hat{L}_0$ and $T^{i+1} \hat{L}_0$ and $L_0$, we have the constraint
\[ \deg(u_i) + \deg(v_{i+1}) =2. \] 
A priori, these are the only restrictions that we have on the degrees. To pin down the degrees exactly, we
need to choose a grading structure on $T^i\hat{L}_0$, and we do this so that $\deg(u_i)=1$,
which then forces that $\deg(v_i)=1$. Any other chord is obtained by composing $u_i$ and $v_i$ in
alternating fashion, hence it follows that there is a graded isomorphism of vector spaces such that
$CW^*(T^i\hat{L}_0, T^i\hat{L}_0) \cong R \langle u_i,v_i \rangle / (u_i^2, v_i^2)$. 

It remains to see that there is no differential. This could either be deduced directly, or we
could argue as follows. We could alternatively choose a grading structure on $T^i\hat{L}_0$ such
that $\deg(u_i)=2$ and $\deg(v_i)=0$. Then the complex $CW^*(T^i\hat{L}_0, T^i\hat{L}_0)$ would be concentrated
in even degree, hence there could not be any rigid holomorphic curve contributing to the
differential. This geometric fact remains the same for any choice of grading structures,
hence it follows that the Floer differential vanishes irrespective of how we grade the Lagrangians
$T^i\hat{L}_0$. 
\ed

We will always grade our Lagrangian $T^i\hat{L}_0$ so that the graded isomorphism given in Lemma
\ref{sympl-dim-Hom-lem} holds. Later, we will see that in fact the isomorphism given in Lemma \ref{sympl-dim-Hom-lem} is an isomorphism of graded rings.

\subsection{Calculations in the relative Fukaya category}

Following the strategy in \cite{LP}, we will not compute the $A_\infty$-algebra $\mathscr{A}$
directly. The presence of higher products (cf. \cite{LPshort}) makes such a direct computation
difficult. On the other hand, the cohomology algebra $A= H^*\mathscr{A}$ is easy to compute.

\begin{prop} 
    There exists a choice of gradings on the Lagrangians $(L_i)$, such that we have an isomorphism of graded algebras
    $$H^*\mathscr{A} \cong E_{1,n}\ot \Z[[t_1,\ldots,t_n]],$$
    where the algebra $E_{1,n}$ is given by \eqref{E1n-eq}.
 
\end{prop}
\Pf. The additive isomorphism is a direct consequence of the intersection pattern of the curves
$L_0,L_1,\ldots, L_n$ and the choice of orientations and gradings. Namely, the given
orientations imply that the rank 1 free module $HF^*(L_0,L_i)$ is supported in even degree for all
$i=1,\ldots, n$. We can pick the grading structure on $L_0$ arbitrarily, and then choose
the grading structure on $L_i$ so that $HF^*(L_0,L_i) \cong HF^0(L_0,L_i)$. This implies that
the rank 1 module $HF^*(L_i,L_0)$ is supported in degree 1, i.e., $HF^*(L_i,L_0) \cong
HF^1(L_i,L_0)$. Note that since $L_i$ are exact, we have that $HF^*(L_i,L_i) \cong H^*(L_i) \cong
H^*(S^1)$ as a ring. 
The rest of the algebra structure is determined by the following pairings, for $i=1,\ldots,n$:
\begin{align} \label{PD} HF^1(L_i, L_0) \otimes HF^0(L_0,L_i) &\to HF^1(L_0,L_0), \\
    \label{PD2}  HF^0(L_0,L_i) \otimes HF^1(L_i,L_0) &\to HF^1(L_i,L_i). 
\end{align}
Note that, as explained in \cite[Sec.7]{SeidelGenus2}, the only contributions to these products
are given by constant holomorphic triangles (and their moduli space is regular). We can either compute their contribution explicitly or observe that as
a consequence, we have that
\[ H^*\mathscr{A} \cong H^*\mathscr{A}_0 \otimes  \Z[[t_1,\ldots,t_n]]. \]
Now, all of the maps (\ref{PD}) and (\ref{PD2}) have to be non-degenerate pairings by the general Poincar\'e duality property of Floer
cohomology (see \cite[Sec. 12e]{SeidelBook}). In $\mathscr{A}_0$ each of these maps is of the
form: \[ \Z \otimes \Z \to \Z. \]
This has to be non-degenerate after reducing mod $p$ for all $p$, hence we see that the map has
to be $(1,1) \to \pm 1$ and the sign does not matter up to isomorphism, so we deduce that
$H^* \mathscr{A}_0$ is isomorphic to $E_{1,n}$. 
\ed

Thus, by Theorem \ref{M1n-ainf-thm} (resp., by Proposition \ref{n=2-ainf-surj-prop}, if $n=2$), 
there exists a family of curves $(C_{mirror},p_1,\ldots,p_n,\om)$ in $\wt{\UU}_{1,n}^{sns}$ over $\Z[[t_1,\ldots,t_n]]$
such that the $A_\infty$-algebra $\mathscr{A}$ is gauge-equivalent to the $A_\infty$-structure on $E_{1,n}\ot \Z[[t_1,\ldots,t_n]]$,
associated with this family via \eqref{curve-to-ainf-map}.

Under this correspondence, we can identify $L_0 \leftrightarrow \mathcal{O}_C$
and $L_i \leftrightarrow \mathcal{O}_{p_i}$. Therefore, to identify the $A_\infty$-algebra $\mathscr{A}$
up to an $A_\infty$-equivalence, it suffices to determine the curve $(C_{mirror},p_1,\ldots, p_n)$. We will do this, by
computing its homogeneous coordinate ring which can be done at the level of the cohomological
category $H^*(\mathcal{F}(\T,D))$:
\[ R_{C_{mirror}} = \bigoplus_{N \geq 0} H^0({\OO}_C (N (p_1+\ldots+p_n))) \cong \bigoplus_{N
    \geq 0} HF^0(L_0, (\tau_{L_1} \circ \tau_{L_2} \circ \ldots \tau_{L_n})^N (L_0) ),\]
where $\tau_{L_1}, \ldots, \tau_{L_n}$ are Dehn twists around the curves $L_1,
\ldots, L_n$, respectively. To make sure that these Dehn twists act on the category
$\mathcal{F}(\T,D)$, we need to give models for them that are \emph{exact}
symplectomorphisms. Furthermore, we prefer to have linear models so that we can explicitly compute
$R_{C_{mirror}}$. Rather than exhibiting explicit models for each $\tau_{L_i}$, we will exhibit
a model for $\tau_{L_1} \circ \tau_{L_2} \circ \ldots \circ \tau_{L_n}$. 

To this end, consider the symplectomorphism $\rho : \T \to \T$ given by
\[ \rho(x,y) = (x, y-nx). \]

This is almost a model for $\tau$, however, it does not preserve the divisor
$D$, so we do not get a symplectomorphism of $\T_0$. To fix this, we observe that the divisor $D
= \{z_1,\ldots,z_n \}$ on $\T$ is sent by $\rho$ to $\{w_1,\ldots, w_n \}$ where \[ w_i = \rho(z_i)=
(\frac{n-i}{n}+\epsilon, \epsilon - (n-i)- n\epsilon) = (\frac{n-i}{n}+\epsilon, \epsilon -
n\epsilon) \in \R^2/\Z^2. \] 

Note that for $\epsilon$ sufficiently small, for all $i$, $z_i$ and $w_i$ are very close to each
other as points in $\T$.
Let $U_i$ be an $\epsilon^2$-neighborhood of the segment $[z_i, w_i]$. We can find a compactly supported symplectomorphism $\delta_i$ of $U_i$ sending $w_i$
to $z_i$ for all $i$. We then set 
\begin{equation}\label{tau-sympl-eq}
\tau = (\delta_1 \circ \ldots \circ \delta_n) \circ \rho.
\end{equation}
Note that $\tau$ defines a symplectomorphism of $\T_0$, such that
    \[ \tau(x,y)= (x, y-nx) \]
outside of an arbitrarily small (depending on $\epsilon$) neighbourhood of $D$.

The following proposition ensures that for a suitable choice of the primitive $\theta$, $\tau$ will
be an exact symplectomorphism of $\T_0$, that is a model for the composition of the Dehn twists around $L_i$, 
which is linear outside of a small neighborhood of $D$. 

\begin{prop} \label{primitive} There exists a primitive $\theta = \iota_Z \omega$ of
    $\omega|_{\T_0}$ such that the
    Lagrangians $L_0$ and $\tau(L_0)$ are exact. Furthermore, for the same choice of $\theta$, the symplectomorphism
    $\tau: \T_0 \to \T_0$ given by \eqref{tau-sympl-eq} is exact and is symplectically isotopic to
    the composition of Dehn twists $\tau_{L_1} \circ \tau_{L_2} \circ \ldots \circ \tau_{L_n}$.
 \end{prop}

\begin{figure}[htb!]
\centering
\begin{tikzpicture} [scale=1.1]
      \tikzset{->-/.style={decoration={ markings,
        mark=at position #1 with {\arrow[scale=2,>=stealth]{>}}},postaction={decorate}}}
         \draw (0,0) -- (6,0);
         \draw (6,6) -- (0,6);
         \draw (6,0) -- (6,6);
         \draw (0,0) -- (0,6);

         \draw[red, ->-=.5]  (0,0) -- (6,0);
         \draw [red, ->-=.5] (0,6) -- (1,0);
         \draw [red, ->-=.5] (1,6) -- (2,0);
         \draw [red, ->-=.5] (2,6) -- (3,0);
         \draw [red, ->-=.5] (4,6) -- (5,0);
         \draw [red, ->-=.5] (5,6) -- (6,0);

         \node at (3,-0.3)   {\footnotesize $L_0$};

\draw[red, thick, fill=red] (3.0,2) circle(.02);
\draw[red, thick, fill=red] (3.7,2) circle(.02);
\draw[red, thick, fill=red] (4.4,2) circle(.02);

\draw[black, thick, fill=black] (3.2,0.5) circle(.01);
\draw[black, thick, fill=black] (3.6,0.5) circle(.01);
\draw[black, thick, fill=black] (4,0.5) circle(.01);

\draw[thick, fill=black] (0.5,0.4) circle(.03);
\draw[thick, fill=black] (1.5,0.4) circle(.03);
\draw[thick, fill=black] (2.5,0.4) circle(.03);
\draw[thick, fill=black] (4.5,0.4) circle(.03);
\draw[thick, fill=black] (5.5,0.4) circle(.03);

\node at (0.6,0.6) {\footnotesize $z_n$};
\node at (1.6,0.6) {\footnotesize $z_{n-1}$};
\node at (2.6,0.6) {\footnotesize $z_{n-2}$};
\node at (4.6,0.6) {\footnotesize $z_2$};
\node at (5.6,0.6) {\footnotesize $z_1$};

\end{tikzpicture}
\caption{Effect of exact symplectomorphisms $\tau$ on $L_0$ } 
\label{alexander}
\end{figure}
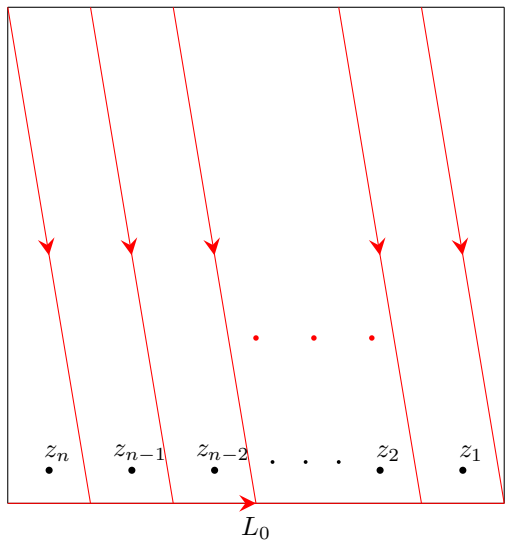

\Pf . To make $L_0$ and $\tau({L}_0)$ exact Lagrangians, we need to exhibit a
Liouville vector field $Z$ on $\T_0$  such that $\int_{L_0} \iota_Z \omega =
\int_{\tau ({L}_0)} \iota_Z \omega =0$. By observing that $\T_0$ retracts
onto a neighborhood of $L_0 \cup \tau (L_0)$ which can be locally identified with a neighborhood of
a plumbing of $T^*L_0$ and $T^*(\tau(L_0))$'s, it is easy to see that there is a Liouville vector
field $Z$ which vanishes along $L_0$ and $\tau(L_0)$ (see Figure
\ref{alexander}). 

Next, recall that a symplectomorphism $\phi: \T_0 \to \T_0$ is exact  if $[\phi^* \theta
- \theta ] = 0 \in H^1(\T_0)$. On the other hand, the compact Lagrangians $L_i$ form a basis of
$H_1(\T_0)$. Hence, to show that $\tau$ is exact, it suffices to check that
$\int_{\tau(L_i)} \iota_Z \omega = \int_{L_i} \iota_Z \omega$ for all $i=0,\ldots, n$. We have chosen the
primitive $\theta = \iota_Z \omega$
above so that both integrals are zero for $i=0$. On the other hand, for $i=1,\ldots, n$, we note
that $\tau(L_i) = L_i$, since $L_i = \{ ((n-i)/n, y) : y \in \R/\Z \}$ for $i=1,\ldots, n$
and $\tau(L_i) = \{ ((n-i)/n , y- (n-i)) : y \in \R/ \Z \}$ are the same curves.

Finally, we will apply Alexander's method (see \cite{primer}) using $L_0, L_1, \ldots, L_n$ as test
curves to show that $\tau$ is isotopic to the composition of Dehn twists around $L_i$. 
Since $\tau({L}_i) = L_i$ for $i=1,\ldots n$, it suffices to check that $\tau({L}_0)$ is isotopic to
the curve that is obtained by Dehn twisting $L_0$ around $L_1,\ldots, L_n$. 
But this is clear from the depiction of $\tau(L_0)$ in Figure \ref{alexander}. \ed

For a fixed $m_0$, we can choose $\epsilon$ small enough so that for all $1 \leq m \leq m_0$, $\tau^m : \T_0 \to \T_0$ is given by 
\[ \tau^m(x,y) = (x, y-mnx) \]
outside an arbitrary small neighborhood of $D$ and is a model for the $m^{th}$ power of the composition of the Dehn
twists around $L_1,\ldots, L_n$. Therefore, we can use these linear models to compute triangle products.

As explained in \cite{LP}, the ring
structure on $\bigoplus_{N \geq 0 } HF^* (L_0, \tau^{N} L_0)$ is determined by the Floer triangle products:
\[ \m_2: HF^*(\tau^{m_1}L_0, \tau^{m_1+m_2}L_0) \otimes HF^*(L_0, \tau^{m_1}L_0) \to HF^*(L_0,
\tau^{m_1+m_2}(L_0)).  \]

This is what we compute next. First, let us observe that since the intersection points of $L_0$
and $\tau^{m}L_0$ all live in degree 0 by our grading choices, the Floer cohomology
$HF^*(L_0, \tau^{m}L_0)$ is freely generated by $L_0 \cap \tau^{m}L_0$. 
Note that $L_0$ intersects $\tau^{m}L_0$ at the points $x + \mathbb{Z}$, where
\[ -nmx \in \mathbb Z. \]
So, $x \in \{ 0,  \frac{1}{nm}, \frac{2}{nm}, \ldots, \frac{nm-1}{nm} \} + \mathbb{Z}$. Thus, we get a one-to-one correspondence between the $x$-coordinates of the intersection
points 
$L_0 \cap \tau^{m}L_0$ and elements of the set $\frac{1}{m}\mathbb{Z}/n\mathbb{Z}$, 
so we write
\[ HF^*(L_0, \tau^{m}L_0) = \bigoplus_{p \in \frac{1}{m}\mathbb{Z}/n\mathbb{Z}}
\mathbb{Z}[[t_1,\ldots,t_n]] x_{m,p}. \]

Note also that the Dehn twist $\tau^{m_1}$ gives an identification of $HF(L_0,
\tau^{m_2} L_0)$ and $HF(\tau^{m_1} L_0, \tau^{m_1+m_2} L_0)$ by mapping the intersection
points $L_0 \cap \tau^{m_2} L_0$ bijectively onto $\tau^{m_1}L_0 \cap
\tau^{m_1+m_2}L_0$. Thus, we have
\[ HF^*(\tau^{m_1}L_0, \tau^{m_1+m_2}L_0) \cong \bigoplus_{p \in
    \frac{1}{m_2}\mathbb{Z}/n\mathbb{Z}} \mathbb{Z}[[t_1,\ldots,t_n]] \tau^{m_1}x_{m_2,p}.
\]

The ring structure on $R_{C_{mirror}} \cong \bigoplus_{N \geq 0}
HF^0(L_0, \tau^{N}L_0)$ is given by the triangle counts
\[ x_{m_2, p_2} \cdot x_{m_1, p_1} = \m_2 (\tau^{m_1} x_{m_2,p_2}, x_{m_1,p_1}). \]

In the following proposition we compute this explicitly.

\begin{prop} Let $x_{m_i,p_i} \in HF^*(L_0, \tau^{m_i} L_0)$, where $p_i \in \frac{1}{m_i}
    \mathbb{Z}/ n\mathbb{Z}$. For $\epsilon \ll 1/ n(m_1+m_2)$, we have:
\[ x_{m_2,p_2} \cdot x_{m_1,p_1} = \m_2 ( \tau^{m_1}(x_{m_2,p_2}), x_{m_1,p_1}) = \sum_{k \in \Z }
x_{m_1+m_2,E(p_1, p_2+ kn)} \Pi_{j=1}^{n} t_j^{n\lambda(\frac{p_1+j}{n}, \frac{p_2+j}{n}+k)} \]

where $$E(a,b)=\frac{m_1a+m_2b}{m_1+m_2},$$
$$\la(a,b)=m_1\phi(a)+m_2\phi(b)-(m_1+m_2)\phi(E(a,b))$$ and $\phi$ is the piecewise linear function given by (\ref{phi-def}). 

\end{prop}
\Pf. This is analogous to the computation given in \cite{LP} with the exception that there are
$n$ different marked points and we need to keep track of intersection numbers of triangles with
respect to these. On the other hand, one can compute these intersection numbers one point at a time,
hence Brion's formula \cite{brion} (see also \cite{Barvinok}) applied in \cite{LP} can be used here as well.

For $k \in \Z$, let us set $p_3 = E(p_1, p_2+kn)$ and $m_3 = m_1 + m_2$. For $p_i \in \frac{1}{m_i}
\Z/ n\Z$ and fixed $j \in \{1,\ldots, n\}$, let us set 
\begin{align*} 
    \frac{p_1+j}{n} &= q_1 + r_1/m_1 n, \text{ where } q_1, r_1 \in \Z \text{ and } 0 \leq r_1 < m_1 n, \\
    \frac{p_2+kn+j}{n} &= q_2 + r_2/m_2 n, \text{ where } q_2, r_2 \in \Z \text{ and } 0 \leq r_2 < m_2 n, \\
    \frac{p_3+j}{n} &= q_3 + r_3/m_3 n, \text{ where } q_3, r_3 \in \Z \text{ and } 0 \leq r_3 < m_3 n.
\end{align*} 
    Then, one can compute 
\begin{align*}
&n\lambda(\frac{p_1+j}{n}, \frac{p_2+j}{n}+k) = \\
& n \left(m_1 \frac{q_1(q_1-1)}{2} + m_2
\frac{q_2(q_2-1)}{2} - m_3\frac{q_3(q_3-1)}{2} \right) + r_1 q_1 + r_2 q_2 -r_3 q_3.
\end{align*}

The count of triangles contributing to $\m_2 (\tau^{m_1}x_{m_2,p_2},x_{m_1,p_1})$ can be enumerated
as embedded triangles, $T(p_1,p_2+kn)$ for $k\in \mathbb{Z}$, in the universal cover $\R^2$. 
The first vertex of $T(p_1,p_2+kn)$ is a lift of $x_{m_1,p_1}$ which we can fix to be the point $(\frac{p_1}{n}, 0)$.
The second vertex is a lift of $\tau^{m_1}(x_{m_2,p_2})$, which lies on the line
of slope $-m_1n$ that passes through $(\frac{p_1}{n}, 0) $, so it has to be of the form 
\[ \left( \frac{p_2}{n}+k , -m_1(p_2+nk-p_1) \right) , \ \ \ k \in \Z. \]
Finally, the third vertex is a lift of $x_{m_1+m_2, E(p_1,p_2+kn)}$ and has the coordinates
\[ (\frac{ m_1 p_1 + m_2 (p_2+kn)}{(m_1+m_2)n}, 0) = (\frac{E(p_1,p_2 +kn)}{n} , 0) , \ \ \ k \in \Z. \]

\begin{figure}[htb!]
\centering
\begin{tikzpicture} [scale=1.1]
      \tikzset{->-/.style={decoration={ markings,
        mark=at position #1 with {\arrow[scale=2,>=stealth]{>}}},postaction={decorate}}}
        
        \node at (-2,6.3) {\footnotesize $L_0$};
        \node at (2,4) {\footnotesize $\tau^{m_1+m_2} L_0$};
        \node at (-2,4) {\footnotesize $\tau^{m_1} L_0$};

        \draw (-4,6) -- (0,6);
        \node at (-5,6) {\footnotesize $(\frac{p_1}{n}, 0)$};
        \draw (0,6) -- (2,2);
         \node at (3,1.6) {\footnotesize 
         $\left( \frac{p_2}{n}+k , -m_1(p_2+nk-p_1) \right)$ };     
        \draw (-4,6) -- (2,2);
        \node at (2,6) {\footnotesize $(\frac{E(p_1,p_2+kn)}{n}, 0)$};

\end{tikzpicture}
\caption{Contributions to the product $x_{m_2,p_2} \cdot x_{m_1,p_1}$} 
\label{figure5}
\end{figure}
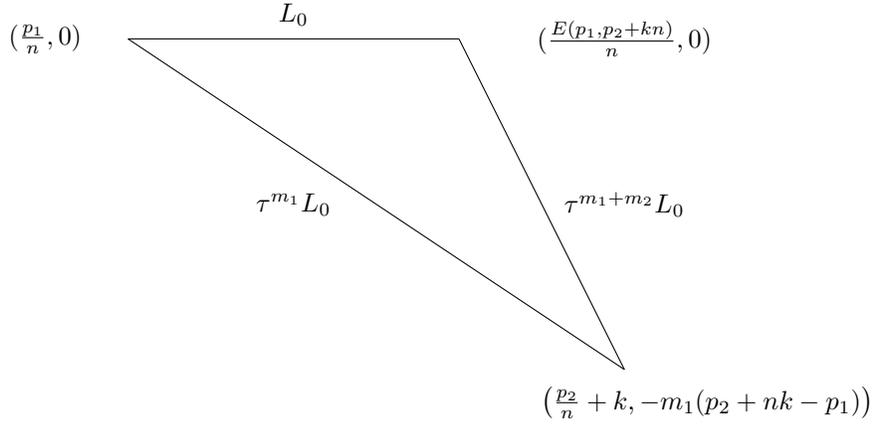

Figure \ref{figure5} shows the type of triangles that appear in the computation.

For a convex region $C \subset \mathbb{R}^2$ and $j=1,\ldots,n$ and $\epsilon>0$, let us consider
the 2-variable Laurent series recording the count of points:
\[ F_{C,j} (x,y) = \sum_{ \{(a,b) \in \mathbb{Z}^2 : (a- \frac{j}{n} + \epsilon, b +\epsilon) \in
C\} }  x^a y^b. \]

Note that $F_{T(p_1,p_2+kn),j}(1,1)$ records the intersection number of the triangle $T(p_1,p_2+kn)$
with the base-point at $z_j = ((n-j)/n + \epsilon, \epsilon)$. We want to show that
$$F_{T(p_1,p_2+kn),j}(1,1) = n\lambda(\frac{p_1+j}{n}, \frac{p_2+j}{n}+k).$$

It is slightly more convenient to work with perturbed lattice points $(\Z+\epsilon)^2$, hence we
rewrite this as:  \[ F_{C,j} (x,y) = \sum_{ \{(a,b) \in \mathbb{Z}^2 : (a + \epsilon, b +\epsilon) \in C(j)\} }  x^a y^b \]
where $C(j)= C + {(\frac{j}{n}, 0)}$ is a horizontal translate of $C$.

To apply Brion's formula, we consider the convex cones $C_1, C_2 , C_3$ with
sides parallel to the sides of the triangle $T(p_1,p_2+kn)+ (\frac{j}{n},0)$ and
the tip points at the vertices $(\frac{p_1+j}{n}, 0)$, $\left(
    \frac{p_2+j}{n}+k ,
-m_1(p_2+nk-p_1) \right)$ and $(\frac{p_3+j}{n}, 0)$. Each of $F_{i,j} := F_{C_i,j}$
is a rational function of $x$ and $y$ and Brion's formula gives:
\[ F_{T(p_1,p_2+kn),j} = F_{1,j} + F_{2,j} + F_{3,j}. \]

To compute $F_{T(p_1,p_2+kn),j}(1,1)$ it suffices to specialize to $y=1$, so let us set $G_{i,j} =
{F_{i,j}}|_{{y=1}}$. We can compute the generating function $G_{i,j}$ by first counting the points in a primitive parallelogram $P_i$ of the cone and then tiling the cone. For example, $P_1$
is the parallelogram with the corners
\[ (\frac{p_1+j}{n},0), (\frac{p_1+j}{n}+1, 0),
    (\frac{p_1+j}{n}+1, -nm_1),
(\frac{p_1+j}{n}+2 , -nm_1). \]

The resulting count can be expressed as follows:
\[ G_{1,j} = g_{1,j}(x) \frac{1}{(1-x)^2} \ \ \ , \ \ \ G_{2,j} = g_{2,j}(x)
    \frac{x^2}{(1-x)^2}
\ \ \ , \ \ \ G_{3,j} = g_{3,j}(x) \frac{-x}{(1-x)^2}, \]
            where $g_{i,j}$ are counts of points in the primitive parallelograms. One can compute these explicitly:
            \[ g_{1,j}(x) =  (nm_1-r_1) x^{q_1+1} + r_1 x^{q_1+2}, \]
            \[ g_{2,j}(x) =  (nm_2-r_2) x^{q_2-1} + r_2 x^{q_2}, \]
            \[ g_{3,j}(x) = (nm_3 - r_3) x^{q_3} + r_3 x^{q_3+1}.  \]
            Let $h(x) = g_{1,j}(x)+ x^2 g_{2,j}(x) - x g_{3,j}(x)$. Then, we have that $h(1)=0$
            and $h'(1)=0$ since $r_1 + r_2 - r_3 = n (m_3 q_3 - m_2 q_2 - m_1 q_1)$. Hence,
            we can compute:
\[ F_{T(p_1,p_2+kn),k}(1,1) = \frac{1}{2} h''(1) = n \lambda( \frac{p_1+j}{n},
\frac{p_2+j}{n}+k). \] \ed

Comparing the above computation of $R_{C_{mirror}}$ with the formulas for multiplication of theta functions
(see Proposition \ref{Bsidecomp}) we derive the following key result.

\begin{cor}\label{curve-isom-cor}
There exists a $\Z[[t_1,\ldots,t_n]]$-linear isomorphism of the homogeneous ring $R_{C_{mirror}}$ with the ring
$\bigoplus_{N\ge 0}H^0(T_n,L^{\ot N})$, where $T_n$ is the $n$-Tate curve with its natural polarization $L$. 
In particular, we get an isomorphism of curves $C_{mirror}\simeq T_n$ over $\Z[[t_1,\ldots,t_n]]$.
\end{cor}

\subsection{Generating objects on the $B$ side}\label{gen-B-sec}


Let $G_n=\cup_{i=1}^n C_i$ be the standard $n$-gon curve over $\Z$.
We assume that the components $C_i\simeq \P^1$ are glued so that the point $0\in C_i$ is identified with the point $\infty\in C_{i+1}$. For each $i$ let $p_i\sub C_i$ be the $\Z$-point with the coordinate $1$, 
and let $\pi_i:\P^1\to G_n$ denote the natural map with the image $C_i$. Note that in the case $n=1$ this becomes
the normalization map. Let us denote by $q_i$ the node at the intersection of $C_i$ and $C_{i+1}$ (so $q_i$ is a closed
subscheme, isomorphic to $\Spec(\Z)$).
We identify indices with elements of $\Z/n$, so $q_0=q_n$.

For a commutative ring $R$ we set $G_{n,R}=G_n\times \Spec(R)$. We still denote by $p_i$ and $q_i$ the 
$R$-points of $G_{n,R}$ obtained from the similar $\Z$-points of $G_n$.

\begin{lem}\label{Perf-C-gen-lem}
Let $(\CC,p_1,\ldots,p_n)$ be a flat proper family of pointed curves (where the marked points are smooth and distinct)
over $\Spec(R)$, where $R$ is a Noetherian ring. Assume that $\OO(p_1+\ldots+p_n)$ is ample on every fiber. Then 
the perfect derived category $\Perf(\CC)$ is split-generated by the objects $(\OO_\CC, \OO_{p_1},\ldots,\OO_{p_n})$. 
\end{lem}

\Pf . Using the twist functors with respect to $\OO_{p_i}$ we obtain all the line bundles $L^m$, where
$L=\OO_C(p_1+\ldots+p_n)$ and $m\in\Z$. But $L$ is ample, and it is well known that all powers of an ample
line bundle generate the perfect derived category (see e.g., \cite[Thm.\ 4]{O-generators}).
\ed

\begin{cor}\label{Perf-C-gen-cor}
For every Noetherian ring $R$,
the category $\Perf(G_{n,R})$ is split-generated by the objects 
\begin{equation}\label{F-collection-eq}
F_0:=\OO_{G_{n,R}} \ \text{ and } \ F_i:=\OO_{p_i}, \ i=1,\ldots,n.
\end{equation}
\end{cor}

Next, we consider generators of the full derived category $D^b(\Coh G_{n,R})$.
Let us define a collection of objects $(\hat{F}_0,\ldots,\hat{F}_n)$ in $D^b(\Coh G_{n,R})$ as follows:
\begin{equation}\label{Fhat-collection-eq}
    \hat{F}_i:=(\pi_{i})_* \OO(-1)[1], \ i=1,\ldots,n, \ \ \hat{F}_0:=\OO_{q_0}.
\end{equation}

\begin{prop}\label{Db-Coh-gen-prop} Assume $R$ is a regular ring. Then
    the category $D^b(\Coh G_{n,R})$ is split-generated by the objects $(\hat{F}_0,\ldots,\hat{F}_n)$.
Furthermore, $D^b(\Coh G_{n,R})$ is generated (as a triangulated
category) by the sheaves $P\ot_R\hat{F}_i$, $i=0,\ldots,n$, where $P$ ranges over finitely generated projective $R$-modules.
\end{prop}

\Pf . Below $P$ denotes a finitely generated projective $R$-module.
Set $C=G_{n,R}$, and
let $\CC\sub D^b(\Coh C)$ (resp., $\CC'\sub D^b(\Coh C)$) 
be the thick subcategory split-generated by $\hat{F}_i$, $i=0,\ldots,n$
(resp., the triangulated subcategory generated by $(P\ot_R\hat{F_i})$). Note that $\CC'\sub \CC$.

\noindent
{\bf Step 1}. One has $P\ot_R (\pi_{i})_* \OO\in \CC'$, $P\ot_R \OO_{q_i}\in \CC'$ for $i=0,\ldots,n-1$. 
This follows easily by considering the exact sequences 
\begin{equation}\label{O(-1)-O-qi-seq}
    0\to (\pi_{i})_* \OO(-1)\to (\pi_{i})_* \OO\to \OO_{q_i}\to 0,
\end{equation}
\begin{equation}\label{O(-1)-O-qi-bis-seq}
    0\to (\pi_{i+1})_* \OO(-1)\to (\pi_{i+1})_* \OO\to \OO_{q_i}\to 0,
\end{equation}
and tensoring them with $P$.
Indeed, for $i=0$ these sequences give $P\ot_R(\pi_{0})_* \OO\in \CC'$, $P\ot_R(\pi_{1})_* \OO\in \CC'$.
Next, the first sequence for $i=1$ gives $P\ot_R\OO_{q_1}\in \CC'$, and then we continue using the induction on $i$. 

\noindent
{\bf Step 2}. One has $\Perf(C)\sub \CC$.
Indeed, the exact sequence
\begin{equation}\label{surgery} 0\to \OO_C\to \bigoplus_{i=1}^n (\pi_{i})_* \OO \to
\bigoplus_{i=1}^n\OO_{q_i}\to 0, \end{equation}
together with Step 1, shows that $\OO_C\in\CC$. Also, the exact sequences
$$
0\to (\pi_{i})_* \OO(-1)\to (\pi_{i})_* \OO\to \OO_{p_i}\to 0
$$
show that $\OO_{p_i}\in \CC$. But $\Perf(C)$ is split-generated by $\OO_C$ and $\OO_{p_i}$ by 
Corollary \ref{Perf-C-gen-cor}.

\noindent
{\bf Step 3}. The category $D^b(\Coh C)$ is split-generated by $\Perf(C)$ and by all the structure sheaves of nodes 
$\OO_{q_i}$, $i=1,\ldots,n$, hence $\CC=D^b(\Coh C)$. Indeed, this is proved similarly to \cite[Prop.\ 2.7]{O-compl}.
Namely, let $\TT\sub D^b(\Coh C)$ be the thick subcategory, split generated by $\Perf(C)$ and $(\OO_{q_i})$,
and let $j:U\hra C$ be the complement to $Z=\cup_{i=1}^n q_i$.
Note that any coherent sheaf supported on $Z$ is obtained as an iterated extension of 
coherent sheaves of the form $q_{i*}F$ for some $F$ on $\Spec(R)$. Since $R$ is regular, each $q_{i*}F$ belongs to
the thick subcategory generated by $\OO_{q_i}$. Let $D^b_Z(\Coh C)\sub D^b(\Coh C)$ be the subcategory of
complexes with cohomology supported on $Z$. Then we obtain the inclusion 
$$D^b_Z(\Coh C)\sub \TT.$$
On the other hand, it is well known that the quotient $D^b(\Coh C)/D^b_Z(\Coh C)$ is naturally equivalent to
$D^b(\Coh U)$. Hence, the projection $D^b(\Coh C)\to D^b(\Coh C)/\TT$ factors as a composition
$$D^b(\Coh C)\rTo{j^*} D^b(\Coh U)\to D^b(\Coh C)/\TT.$$
Now given a coherent sheaf $F$ on $C$ we can find an exact sequence of the form
$$0\to G\to P_N\to P_{N-1}\to \ldots P_1\to P_0\to F\to 0$$
where $P_i$ are vector bundles on $C$, and $N>d$, where $d$ is the global dimension of $R$. 
Then the induced morphism $\a:F\to G[N+1]$ has a cone in $\Perf(C)$, so it becomes an isomorphism
in $D^b(\Coh C)/\TT$. On the other hand, $j^*\a=0$ since $\Coh(U)$ has homological dimension $d+1$. 
Therefore, $\a$ becomes zero in $D^b(\Coh C)/\TT$, which gives that $F\in\TT$, as required.

\noindent
{\bf Step 4}. Now we can prove that $\CC'=D^b(\Coh C)$, i.e., $D^b(\Coh C)$
is generated by objects of the form $P\ot_R \hat{F}_i$.
By the previous steps we know that $\CC'$ is dense in $D^b(\Coh C)$. Hence, by Thomason's theorem
(\cite{thomason}), it is enough to see that
the classes of our objects generate the Grothendieck group of $D^b(\Coh C)$. Let us consider
the localization sequence
$$\ldots\to K_0(D^b_Z(\Coh C))\to K_0(D^b(\Coh C))\to K_0(D^b(\Coh U))\to 0.$$
Since $R$ is regular, the objects of the form $P\ot \OO_{q_i}$ generate $D^b_Z(\Coh C)$.
Thus, by Step 1, we have $D^b_Z(\Coh C)\sub\CC'$.
It remains to check that the images of the classes of $(P\ot\hat{F}_i)_{i=1,\ldots,n}$ generate
$K_0(D^b(\Coh U))$. But $U$ is the disjoint union of $n$ copies of $\Spec(R[t,t^{-1}])$,
so we have an identification
$$K_0(D^b(\Coh U))\simeq \bigoplus_{i=1}^n K_0(R[t,t^{-1}]).$$
It remains to observe that the map $[P]\mapsto [P\ot\hat{F}_i|_U]$ 
corresponds to the standard map
\begin{equation}\label{Rt-isom}
K_0(R)\to K_0(R[t,t^{-1}])
\end{equation}
followed by the inclusion as the $i^{th}$ component of the above direct sum. But \eqref{Rt-isom} is known to be
an isomorphism, since $R$ is regular. 
\ed

\begin{cor}\label{Db-Coh-gen-cor} 
If $R$ is a regular ring with $\hat{K}_0(R)=0$ (here $\hat{K}_0\sub K_0$ is the kernel of the rank homomorphism)
then the category $D^b(\Coh G_{n,R})$ is generated by the objects $(\hat{F}_i)_{i=0,\ldots,n}$ as
a triangulated category. In particular, this is true if $R=\Z$ or $R$ is a field.
\end{cor}

\Pf . Indeed, in this case every finitely generated projective module over $R$ is stably free, so it has a $2$-term resolution
by free modules.
\ed

Our next result is that the collections $(F_0,\ldots,F_n)$ and $(\hat{F}_0,\ldots,\hat{F}_n)$
(see \eqref{F-collection-eq} and \eqref{Fhat-collection-eq}) are dual.

\begin{prop}\label{dual-generators-prop} 
Let $R$ be a Noetherian ring, and consider the objects $(F_i)$ and $(\hat{F}_i)$ on $G_{n,R}$.
One has 
$$R\Hom(F_i,\hat{F}_j)=R\Hom(\hat{F}_j,F_i)=0 \text{ for } i\neq j,$$
$$R\Hom(F_i,\hat{F}_i)\simeq R\Hom(\hat{F}_i,F_i)[1]\simeq R$$ 
for $i=0,\ldots,n$.
\end{prop}

\Pf . Recall that since $F_i$ are perfect objects, by Serre duality, 
$$R\Hom(F_i,\hat{F}_j)\simeq R\Hom(R\Hom(\hat{F}_j,F_i)[1],R).$$
Thus, the assertion follows from 
$$R\Hom(\OO_C,(\pi_{i})_* \OO(-1))=R\Ga(\P^1,\OO(-1))=R\Hom(\OO_{p_i},\OO_{q_1})=0,$$
$$R\Hom((\pi_{i})_* \OO(-1),\OO_{p_i})=0 \text{ for } i\neq j,$$
$$R\Hom(\OO_C,\OO_{q_1})\simeq R\Hom((\pi_{i})_* \OO(-1),\OO_{p_i})\simeq R.$$
\ed

\begin{lem}\label{Ext-nodes-lem} 
For $C=G_{n,R}$ and any node $q_i\in C$ 
one has an isomorphism of graded algebras 
$$\Ext^*_C(\OO_{q_i},\OO_{q_i})\simeq R\langle u,v\rangle/(u^2,v^2),$$
where $\deg(u)=\deg(v)=1$.
\end{lem}

\Pf . First, we claim that this algebra does not change if we replace $C$ by the affine curve $\Spec(A)$, where $A=R[x,y]/(xy)$.
Indeed, note that for a locally free sheaf $V$ one has $\Ext^{>0}(V,\OO_{q_i})=0$.
Hence, we can use any locally free resolution $V_\bullet$ of $\OO_{q_i}$ to compute $\Ext^*_C(\OO_{q_i},\OO_{q_i})$. 
Since the similar assertion holds for the $\Ext^*$ computed over the affine open $\Spec(A)$, our claim follows.
Thus, we have to understand the algebra $\Ext^*_A(R,R)$. But the algebra $A$ is Koszul since it is given by monomial
relations $xy=yx=0$, so the $\Ext$-algebra is just given by the quadratic dual.
\ed

\subsection{Proofs of Theorems A and B(i)} $\phantom{x}$

\noindent
{\bf Proof of Theorem A:} We follow the same strategy as in the proof of \cite[Thm.\ A]{LP}.
Recall that by Theorem \ref{M1n-ainf-thm} (resp., by Proposition \ref{n=2-ainf-surj-prop}, if $n=2$), 
there exists a family of curves $(C_{mirror},p_1,\ldots,p_n,\om)$ in $\wt{\UU}_{1,n}^{sns}$ over $\Z[[t_1,\ldots,t_n]]$, such that
the corresponding $A_\infty$-structure on $E_{1,n}\ot \Z[[t_1,\ldots,t_n]]$ is gauge equivalent to the
$A_\infty$-endomorphisms algebra $\AA$ of the object $L_0\oplus\ldots\oplus L_n$ of $\FF(\T,D)$.

We have already seen that the objects $(L_0,\ldots,L_n)$ split-generate the derived category of $\FF(\T,D)$,
while the objects $(F_0,\ldots,F_n)$ split-generate $\Perf(C_{mirror})$ (see Lemmas \ref{relgenerate} and \ref{Perf-C-gen-lem}).
Hence, the equivalence of the corresponding $A_\infty$-algebras implies the equivalence between the
$A_\infty$-categories $D^\pi\FF(\T,D)$ and $\Perf(C_{mirror})$.

Finally, Corollary \ref{curve-isom-cor} implies that $C_{mirror}$ is isomorphic to the $n$-Tate curve (as a family of curves
over $\Z[[t_1,\ldots,t_n]]$). 
\ed 

As a corollary (from the above proof) we obtain a proof of Theorem B(i) stated in the introduction:

\begin{cor}\label{Fuk-Perf-cor} For any commutative Noetherian ring $R$ there is an $R$-linear equivalence of 
the split-closed derived category of $\FF(\T_0)\ot R$ with $\Perf(G_{n,R})$, sending the objects
$(L_0,L_1,\ldots,L_n)$ to $(\OO, \OO_{p_1},\ldots,\OO_{p_n})$.
\end{cor}

\Pf . We use the fact that the construction of the $A_\infty$-structure on $E_{1,n}\ot A$ associated with a family of curves
in $\wt{\UU}^{sns}_{1,n}$ over $\Spec(A)$ is compatible with the base changes $A\to A'$.
The $n$-Tate curve is such a family over $\Z[[t_1,\ldots,t_n]]$. Applying the base change 
$\Z[[t_1,\ldots,t_n]]\to\Z\to R$ gives $G_{n,R}$. Since the relevant objects still split-generate $D^\pi(\FF(\T_0)\ot R)$
and $\Perf(G_{n,R})$, respectively, this implies the result.
\ed

{\bf Another proof of Theorem B(i):}  We can give another proof of Corollary
\ref{Fuk-Perf-cor}, which still uses the connection to the moduli space $\wt{\UU}_{1,n}^{ns}$ but
replaces the computations leading to Corollary \ref{curve-isom-cor} with Theorem \ref{B-wheel-char-thm}. 

Assume first that $R=k$ is a field.
Then we can check that conditions (i) and (ii)' of Theorem \ref{B-wheel-char-thm} hold for
$\mathscr{A}_0\ot k$. Indeed, condition (i) follows from \cite[Thm.\ 8]{LPshort}
when $\text{char}(k) \neq 2$ or $3$ and from the main result of \cite{LP} in arbitrary
characteristic. As for condition (ii)', for each $i=1,\ldots, n$, let $N_i$ be a Weinstein
neighborhood of the union of $L_0$ and $L_i$. Notice that $N_i$ is a Liouville subdomain in $\T_0$
whose completion $\hat{N}_i$ is a once-punctured torus. Now, the work \cite{AbouzSeidel} of Abouzaid and Seidel
gives for each $i$ a restriction functor
\[ \mathcal{W}(\T_0) \to \mathcal{W}(\hat{N}_i). \] 
At the level of objects, this functor intersects a Lagrangian in $\T_0$ with the neighborhood $N_i$ and extends it to the completion in an obvious way. In the case at hand, these restriction
functors are easy to understand at the level of objects. In particular, compact Lagrangians contained in
$N_i$ go to their obvious representatives in $\hat{N}_i$. By abuse of notation, we do not distinguish these notationally. 

By pre-composing with the full and faithful inclusion of $\mathcal{F}(\T_0) \to
\mathcal{W}(\T_0)$, we obtain functors
\[ r_i: \mathcal{F}(\T_0) \to \mathcal{W}(\hat{N}_i). \]
Now, it follows easily from Lemma A.2 in \cite{AAEKO} that for $i\neq j$, $r_i(L_j)$ is a non-compact Lagrangian given by a cotangent fibre to $L_0$
inside $\hat{N}_i$. So, $r_i$ sends $\langle L_0, L_i \rangle$ into
$\mathcal{F}(\hat{N}_i)$, while it sends $L_j$, for $i\neq j$, to an object of
$\mathcal{W}(\hat{N}_i)$ which has infinite-dimensional endomorphisms. Hence, the functors $r_i$, for $i=1,\ldots n$, distinguish the
subcategories $\langle L_0 , L_i \rangle$ of $D^\pi\mathcal{F}(\T_0)$, split-generated by $L_0$ and $L_i$. This completes the verification of condition (ii)' of Theorem
\ref{B-wheel-char-thm}, hence we deduce that $\AA_0\ot k$ is equivalent to the $A_\infty$-algebra
associated with the standard $n$-gon over $k$. 

Next, let us consider the case $R=\Z$. As in the proof of Theorem A, we see that $\AA_0$
is equivalent to the $A_\infty$-algebra coming from a family of curves
$(C,p_1,\ldots,p_n,\om)$ in $\wt{\UU}_{1,n}^{sns}$ over $\Z$. Furthermore, the above
argument shows that the base change of $(C,p_1,\ldots,p_n)$ with respect to any homomorphism
$\Z\to k$, where $k$ is a field, gives the standard $n$-gon over $k$ (with one marked point on each component).
Hence, the family $(C,p_1,\ldots,p_n)$ corresponds to a morphism from $\Spec(\Z)$ to an affine neighborhood
of the standard $n$-gon (with marked points) in the moduli space of $n$-pointed stable curves without automorphisms.
We have another such morphism, which corresponds to the standard $n$-gon over $\Z$.
Since these morphisms agree on the generic point of $\Spec(\Z)$, they are in fact the same, so
$C$ is the standard $n$-gon over $\Z$. 
Thus, the $A_\infty$-algebra $\AA_0$ is equivalent to the one coming from $G_{n,\Z}$.
As in Corollary \ref{Fuk-Perf-cor}, this implies the equivalence of $\AA_0\ot R$ with the $A_\infty$-structure
associated with $G_{n,R}$, and hence the equivalence of $D^\pi(\FF(\T_0)\ot R)$ with $\Perf(G_{n,R})$.
\ed

\subsection{Characterization of the generators on the B side}


\begin{lem}\label{DbCoh-lem} 
Let $R$ be a Noetherian ring, and let
$X$ be a projective scheme over $\Spec(R)$. Assume that $\FF\in D(\Qcoh X)$ satisfies the
following property: $R\Hom(P,F)$ is in $D^b(R)$ for every perfect complex $P$ on $X$. Then $\FF\in D^b(\Coh X)$.
\end{lem}

\Pf . Let $i:X\to \P^n_R$ be a closed embedding. It is enough to prove that $i_*\FF$ is in $D^b(\Coh \P^N_R)$.
Note that by our assumption for every $m\in\Z$ we have
$$R\Hom(\OO(m),i_*F)\simeq R\Hom(i^*\OO(m),\FF)\in \Perf(R).$$
Now let us consider the following sequence of mutations of $\FF_0=i_*\FF$:
$$\FF_1\to R\Hom(\OO,\FF_0)\ot_R \OO\to \FF_0\to \ldots,$$
$$\FF_2\to R\Hom(\OO(-1),\FF_1)\ot_R \OO(-1)\to \FF_1\to \ldots,$$
$$\ldots$$
$$\FF_{n+1}\to R\Hom(\OO(-n),\FF_n)\ot_R \OO\to \FF_n\to \ldots$$
Note that all the middle terms are in $D^b(\Coh \P^N_R)$ (this can be seen by induction).
Then we have $\Hom(\OO(m),\FF_{n+1})=0$ for $m=0,-1,\ldots,-n$.
Since the category $D(\Qcoh\P^n_R)$ is generated by $\OO,\OO(-1),\ldots,\OO(-n)$ (see
\cite[Thm.\ 2.1.2, 3.1.1]{B-VdB}), it follows that $\FF_{n+1}=0$. Now the above triangles show that $i_*\FF\in D^b(\Coh \P^n_R)$.
\ed

\begin{lem}\label{char-nodes-lem} 
Let $R$ be a Noetherian ring, and let $C=G_{n,R}$.
Let $F\in D(\Qcoh C)$ be such that $\Ext^*(\OO_{p_i},F)=0$ for all $i=1,\ldots,n$,
and $\Ext^*(\OO_C,F)\simeq R$, concentrated in degree $0$.
Then there exists an $R$-point $p:\Spec(R)\to C$ such that $F\simeq\OO_p$. 
Furthermore, assume that $R=k$ is a field.
Then $\Ext^*(F,F)$
is infinite-dimensional (equivalently, $\Ext^1(F,F)$ is $2$-dimensional) if and only if $p$ is one of the nodes.
\end{lem}

\Pf . First, since $\OO_{p_i}$ and $\OO_C$ generate $\Perf(C)$ (see Corollary \ref{Perf-C-gen-cor}),
we see that our assumptions imply that $R\Hom(P,F)$ is in $\Perf(R)$ for every $P\in\Perf(C)$.
Hence, using Lemma \ref{DbCoh-lem} we derive that $F$ is in $D^b(\Coh C)$. 

Next, let us consider the case when $R$ is a field.
Then using the spectral sequences
of the form
$$\bigoplus_k \Ext^q_C(\und{H}^{k-p}(P),\und{H}^k(F))\implies \Ext^{p+q}(P,F)$$
for $P\in\Perf(C)$, as in \cite[Lem.\ 8.9]{LP}\footnote{The relevant spectral 
sequences converge due to the fact that $R\und{\Hom}(P,F)\in D^b(\Coh C)$.} 
we see that that $F$ is supported on 
the open affine subset $C\setminus \{p_1,\ldots,p_n\}$, and then deduce that it is of the form $\OO_p$.

The last assertion follows from the fact that over a field, $\Ext^*(\OO_q,\OO_q)$ is infinite-dimensional 
(resp., $\Ext^1(\OO_q,\OO_q)$ is $2$-dimensional) for each node
$q$ (see Lemma \ref{Ext-nodes-lem}), 
whereas $\Ext^*(\OO_p,\OO_p)$ is finite-dimensional (resp., $\Ext^1(\OO_p,\OO_p)$ is $1$-dimensional)
for each smooth point $p$.

Now let us consider the case of general $R$. Since $C$ is flat over $R$, for every $P\in \Perf(C)$
the formation of $R\Hom(P,F)$ is compatible with the base change. Thus, we deduce that for every homomorphism
$R\to k$, with $k$ a field, the object $F\otimes^{\dL}_R k$ is a structure sheaf of a point on $C\times_{\Spec(R)}\Spec(k)$.
This implies that $F$ is in fact the push-forward of a line bundle on $\Spec(R)$ with respect to some
section $p:\Spec(R)\to C$. Finally, the condition $H^0(C,F)\simeq R$ implies that this line bundle is trivial.
\ed

\begin{lem}\label{char-dual-to-point-lem} Let $C=G_{n,k}$, where $k$ is a field.

\noindent
(i) Assume $n>1$.
Let $F\in D(\Qcoh C)$ be such that $\Ext^*(\OO,F)=0$, $\Ext^*(\OO_{p_i},F)=0$ for $i=2,\ldots,n$,
and $\Ext^*(\OO_{p_1},F)$ is one-dimensional, concentrated in degree $1$. Then $F\simeq
    \hat{F}_1=(\pi_{1})_* \OO(-1)$.

\noindent
(ii) Now let $n=1$. Then any $F\in D(\Qcoh C)$ such that $\Ext^*(\OO,F)=0$ and $\Ext^*(\OO_{p_1},F)$ is one-dimensional,
sitting in degree $1$, is either a nontrivial line bundle of degree $0$ on $C$ or isomorphic to $\pi_*\OO(-1)$, where
$\pi:\P^1\to C$ is the normalization map.
\end{lem}

\Pf . First, as in Lemma \ref{char-nodes-lem}, we see that such $F$ is in $D^b(\Coh C)$. Next, we observe that
for $P=\OO$ or $P=\OO_{p_i}$ the spectral sequence
$$E_2^{rs}=\Ext^r(P,\und{H}^s(F)) \implies \Ext^{r+s}(P,F)$$
degenerates at $E_2$, since $\Ext^r(P,\und{H}^s(F))=0$ for $r\neq 0,1$.
Thus, each cohomology sheaf $G=\und{H}^s(F)$ still satisfies $H^*(C,G)=0$ and $\Ext^*(\OO_{p_i},G)=0$ for
$i=2,\ldots,n$. Note that the condition $H^0(C,G)=0$ implies that $G$ is torsion free. Hence, from the vanishing
of $\Ext^*(\OO_{p_i},G)$ for $i=2,\ldots,n$ we deduce that $G$ is supported on $C_1$.

Assume first that $n>1$.
Then we claim that $G$ is necessarily of the form $V\ot (\pi_{1})_* \OO(-1)$ for some vector space $V$, 
which will imply our statement due
to the condition $\dim \Ext^*(\OO_{p_1},G)=1$. Indeed, any torsion free coherent sheaf
supported on $C_1$ is necessarily of the form $G=(\pi_{1})_* \wt{G}$, as follows from the classification of torsion
free modules over $k[[x,y]]/(xy)$ (see \cite[Sec.\ 7]{BD}, \cite{Bass}). Now we have $H^*(C,G)=H^*(\P^1,\wt{G})=0$, hence, 
$\wt{G}\simeq V\ot \OO(-1)$.

Now assume $n=1$. Then using the condition $\dim \Ext^*(\OO_{p_1},G)=1$ we see that $G$ is a rank $1$ torsion free
sheaf on $C$. It is well known that such $G$ is either a line bundle or has form $\pi_*L$, where
$L$ is a line bundle on $\P^1$ (see e.g., \cite[Sec.\ 1.2]{Jarvis}).
Now the condition $H^*(C,G)=0$ implies that
in the former case $G$ has to be a nontrivial line bundle of degree $0$, while in the latter case $G\simeq \pi_*\OO(-1)$.
\ed

\subsection{Fully faithfulness of the Yoneda functor $\WW(\T_0)\to D(\mod \text{-}\FF(\T_0))$}

In this section $R$ denotes a commutative Noetherian ring.

If $\CC$ is an enhanced triangulated category, $\CC_0\sub \CC$ a full subcategory, then
by {\it Yoneda functor} for the pair $(\CC,\CC_0)$ we mean the functor
$$\CC\to D(\mod \text{-}\CC_0)$$
sending $X$ to the module $X_0\mapsto \Hom_{dg}(X_0,X)$.

The following result is well known (see \cite[Thm.\ 8.9]{Toen}).

\begin{lem}\label{B-Qcoh-lem} 
    The Yoneda functor $D(\Qcoh G_{n,R})\to D(\mod \text{-}\Perf(G_{n,R}))$ is an equivalence.
\end{lem}

This implies that the Yoneda functor $D^b(\Coh G_{n,R})\to D(\mod \text{-} \Perf(G_{n,R}))$ is fully faithful.
We want to prove the following analog of this result on the symplectic side.

\begin{thm}\label{A-fully-faithful-thm} Assume that $R$ is a regular ring.
Then the Yoneda functor 
$$Y:\WW(\T_0)\ot R\to D(\mod \text{-}\FF(\T_0)\ot R)$$ is fully faithful.
\end{thm}

The proof will substantially use the equivalence
$$\FF(\T_0)\ot R\simeq \Perf(G_{n,R}),$$
sending $(L_0,L_1,\ldots,L_n)$ to $(\OO, \OO_{p_1},\ldots,\OO_{p_n})$
(see Corollary \ref{Fuk-Perf-cor}).
Namely, using this equivalence and Lemma \ref{B-Qcoh-lem}
we can identify $D(\mod \text{-}\FF(\T_0))$ with $D(\Qcoh G_{n,R})$, and

view the Yoneda functor for the pair $(\WW(\T_0), \FF(\T_0))$
as a functor
\begin{equation}\label{Yoneda-W-Qcoh-eq}
Y:\WW(\T_0)\otimes R\to D(\Qcoh G_{n,R}),
\end{equation}
whose restriction to $\FF(\T_0)$ is fully faithful and sends $(L_0,L_1,\ldots,L_n)$ to $(\OO, \OO_{p_1},\ldots,\OO_{p_n})$.

We will use the following notation below: for a coherent sheaf $F$ on $G_{n,R}$, for which we have already constructed
an object $\Lambda_F$ in $\WW(\T_0)\otimes R$ 
such that $Y(\Lambda_F)\simeq F$, we will write $[F]=\Lambda_F$. 

Note that the functor $Y$ extends naturally to
$D^b (\WW(\T_0)\ot R)$ and that for $F$ in $\Perf(G_{n,R})$
we already have an object $[F]$ in $D^\pi(\FF(\T_0)\otimes R)$ such that $Y([F])=F$, due to Corollary \ref{Fuk-Perf-cor}. 
For example, $[\OO]=L_0$, $[\OO_{p_i}]=L_i$.

\begin{lem}\label{Y-hat-L0-lem} 
Assume that $\Spec(R)$ is connected. Then there exists $i_0$ such that $Y(\hat{L}_0)\simeq \OO_{q_{i_0}}$.
\end{lem}

\Pf . By Lemma \ref{char-nodes-lem}, there exists an $R$-point $p:\Spec(R)\to G_{n,R}$ such that
$Y(\hat{L}_0)\simeq \OO_p$. Furthermore, for every homomorphism $R\to k$, where $k$ is a field,
by the same Lemma, the $k$-point associated with $p$ is a node of $G_{n,k}$. This implies that 
the image of $p$ is contained in $\sqcup_{i=1}^n q_i(\Spec(R))$. Since $\Spec(R)$ is connected,
the assertion follows.
\ed

\begin{rem} Later we will prove that $Y(\hat{L}_0)\simeq\OO_{q_n}$.
\end{rem}

Recall that we have an exact symplectomorphism $T : \T_0 \to \T_0$ given by translation (see \eqref{T-translation-eq}). 
This is an exact symplectomorphism of $\T_0$ which is of
contact type at infinity. Such an exact symplectomorphism acts on all the data used to construct
$\mathcal{F}(\T_0)$ and $\mathcal{W}(\T_0)$, so we get the following result.

\begin{lem} \label{cyclicsym} The Yoneda functor $Y : \WW(\T_0)\ot R \to D(\mod\text{-}\FF(\T_0)\ot R)$ is $T$-equivariant.
\ed
\end{lem}

\noindent
{\it Proof of Theorem \ref{A-fully-faithful-thm}}.
It is enough to consider the case when $\Spec(R)$ is connected. 
In this case by Lemma \ref{Y-hat-L0-lem}, we have $[\OO_{q_{i_0}}]=\hat{L}_0$. Since $Y$ is compatible with
the cyclic symmetry (see Lemma \ref{cyclicsym}), 
the objects $[\OO_{q_1}],\ldots,[\OO_{q_n}]$ are obtained from $\hat{L}_0$ by the action of $T$
(see \eqref{T-translation-eq}).
Hence, by Lemma \ref{another-generators-W-lem}, $\WW(\T_0)\ot R$ is split-generated by 
$\FF(\T_0)\ot R$ and by $[\OO_{q_1}],\ldots,[\OO_{q_n}]$.
Since we already know that $Y$ is fully faithful when restricted to $\FF(\T_0)\ot R$,
it is enough to prove that for any $a\in\FF(\T_0)$ and $i,j\in [1,n]$,
the following natural maps are isomorphisms:
\begin{equation}\label{Hom-a-q-eq}
\Hom^*(a,[\OO_{q_i}])\to \Hom^*(Y(a),\OO_{q_i}),
\end{equation}
\begin{equation}\label{Hom-q-a-eq}
\Hom^*([\OO_{q_i}],a)\to \Hom^*(\OO_{q_i},Y(a)),
\end{equation}
\begin{equation}\label{Hom-q-q-eq}
\Hom^*([\OO_{q_i}],[\OO_{q_j}])\to \Hom^*(\OO_{q_i},\OO_{q_j}).
\end{equation}
Note that \eqref{Hom-a-q-eq} is an isomorphism by the definition of $Y$.
To see that \eqref{Hom-q-a-eq} is an isomorphism, it is enough to consider the case when
$a$ is one of the generators $L_i$ of $\FF(\T_0)$. Since in this case $Y(a)$ is a spherical object on $G_{n,R}$, by Serre duality, we
have a perfect pairing
$$\Hom^i(\OO_{q_i}, Y(a))\ot \Hom^{1-i}(Y(a), \OO_{q_i})\to \Hom^1(Y(a),Y(a))\simeq R.$$
We have a similar perfect pairing on the symplectic side (see Proposition \ref{exactlags}). 
Thus, the assertion follows from the fact that
\eqref{Hom-a-q-eq} is an isomorphism. 

Note that $\Hom^*([\OO_{q_i}],[\OO_{q_j}])=0$ for $i\neq j$.
Thus, it remains to prove that \eqref{Hom-q-q-eq} is an isomorphism for $i=j$. 
By Lemmas \ref{sympl-dim-Hom-lem} and \ref{Ext-nodes-lem}, 
both sides in \eqref{Hom-q-q-eq} for $i=j$ are free $R$-modules of the same rank. 
Hence, it is enough to prove the same assertion for the morphism \eqref{Hom-q-q-eq} reduced modulo any
maximal ideal in $R$.
Since the formation of the maps \eqref{Hom-q-q-eq} is compatible with the change of scalars
$R\to R'$, we can assume for the rest of the proof that $R=k$ is a field. We set $C=G_{n,k}$.

By cyclic symmetry, it is enough to consider the case of $i=i_0$.
Let us write $q=q_{i_0}$ for brevity.
We know by Lemma \ref{Ext-nodes-lem} that the algebra $\Hom^*(\OO_q,\OO_q)$ is generated in degree $1$.
Thus, it suffices to check that the map
$$\Hom^1([\OO_q],[\OO_q])\to \Hom^1(\OO_q,\OO_q)$$
is an isomorphism. 

Let $J_q\in D^b(\Coh C)$ denote the ideal sheaf of $q$, and let us define $[J_q]$ from the exact triangle 
\[ [J_q]\to [\OO_C]\xrightarrow{\alpha} [\OO_q]\xrightarrow{\delta} [J_q][1] \] 
in $\WW(\T_0)$, where $\a$ corresponds to $1\in H^0(C,\OO_q)$ under the isomorphism 
$\Hom^0([\OO_C],[\OO_q])\simeq \Hom^0(\OO_C,\OO_q)$.
Note that since $Y$ is an exact functor, the image of the above exact triangle under $Y$ is the standard triangle
\[ J_q\to \OO_C\to \OO_q\rightarrow J_q[1]. \] 
In particular, $Y([J_q])\simeq J_q$. 

To describe the object $[J_q]$ concretely as a Lagrangian, 
we note that the Dehn twist triangle \eqref{DTT-eq} with
            $K = L_0 = [\OO_C]$ and $L = \hat{L}_0 =[\OO_q]$ gives an isomorphism 
            \[ [J_q] \simeq \tau_{L_0}(\hat{L}_0)[-1]. \]

Note that we have a commutative square
\begin{diagram}
\Hom^0([J_q],[\OO_q]) &\rTo{\circ\de}& \Hom^1([\OO_q],[\OO_q])\\
\dTo{}&&\dTo{}\\
\Hom^0(J_q,\OO_q) &\rTo{}& \Hom^1(\OO_q,\OO_q)\\
\end{diagram}
We claim that the horizontal arrows in this square are  
isomorphisms. Indeed, for the bottom arrow this follows from the long exact sequence
$$0\to \Hom^0(\OO_q,\OO_q)\xrightarrow{\sim} \Hom^0(\OO_C,\OO_q)\to \Hom^0(J_q,\OO_q)\to \Hom^1(\OO_q,\OO_q)\to \Hom^1(\OO_C,\OO_q)=0.
$$
For the top arrow this follows similarly from the exact triangle defining $[J_q]$.
Therefore, it is enough to check that the Yoneda functor induces an isomorphism
$$\Hom^0([J_q],[\OO_q])\to \Hom^0(J_q,\OO_q).$$

By Lemma \ref{J-sympl-lem} below, there exists a line bundle $\mathcal{L}$ on $C$ (defined over $k$)
such that the pairings 
\begin{equation}\label{line-bundle-pairing-A}
    \Hom^0([J_q],[\OO_q])\otimes_k \Hom^0([\mathcal{L}],[J_q])\to
    \Hom^0([\mathcal{L}],[\OO_q])\simeq k,
\end{equation}
\begin{equation}\label{line-bundle-pairing-B}
    \Hom^0(J_q,\OO_q)\otimes_k \Hom^0(\mathcal{L},J_q)\to \Hom^0(\mathcal{L},\OO_q)\simeq k,
\end{equation}
given by the composition, are perfect pairings between $2$-dimensional vector spaces (note that $[\LL]$ is in 
$D^\pi(\FF(\T_0)\ot k)$).
Thus, we have a commutative square
\begin{diagram}
    \Hom^0([J_q],[\OO_q])&\rTo{}&  \Hom_k\bigl(\Hom^0([\mathcal{L}],[J_q]),\Hom^0([\mathcal{L}],[\OO_q])\bigr)\\
\dTo{Y}&&\dTo{Y}\\
    \Hom^0(J_q,\OO_q)&\rTo{}&  \Hom_k\bigl(\Hom^0(\mathcal{L},J_q),\Hom^0(\mathcal{L},\OO_q)\bigr)
\end{diagram}
in which the right vertical arrow is an isomorphism by the definition of $Y$ (since $[\mathcal{L}]$ is in $D^\pi(\FF(\T_0)\otimes k)$)
and the horizontal arrows are isomorphisms.
Hence, we deduce that the left vertical arrow is also an isomorphism, which is what we needed.
\ed

\begin{lem}\label{J-sympl-lem} Let $C=G_{n,k}$, where $k$ is a field.
For $n>1$ 
    let $\mathcal{L}$ be a line bundle that has degrees $-1$ and $-2$ when restricted to the components
of $C$, passing through $q$, and is trivial on all the other components. In the case $n=1$ let
$\mathcal{L}$ be a line bundle of
degree $-3$ on $C$. Then the pairing 
\eqref{line-bundle-pairing-B} is perfect.
Furthermore, for $\mathcal{L}=\OO_C(-2p_{i_0}-p_{i_0+1})$ (resp., $\mathcal{L}=\OO_C(-3p_1)$ if $n=1$) there exists an 
object $[\mathcal{L}]\in\FF(\T_0)$ such that $Y([\mathcal{L}])\simeq \mathcal{L}$ and the pairing
\eqref{line-bundle-pairing-A} is perfect.
\end{lem}

\Pf . First, let us check the assertion on the B side in the case $n>1$.
Let $C^+$ and $C^-$ be the components of $C$, passing through $q$, and let $q^+\in C^+$ and $q^-\in C^-$
be the points that glue into $q$. Let also $p^+\in C^+$ and $p^-\in C^-$ be the points corresponding to the other nodes
on $C^+$ and $C^-$ (we can assume that under an 
isomorphism $C^\pm\simeq \P^1$ the point $q^\pm$ corresponds to $\infty$,
while the point $p^\pm$ corresponds to $0$). Since $\mathcal{L}$ is trivial on all the other components, we can identify 
$\Hom^0(\mathcal{L},J_q)\simeq H^0(C,\mathcal{L}^{-1}\ot J_q)$ with the vector space
\begin{align*}
V=\{ & (s^+,s^-)\in H^0(C^+, \mathcal{L}^{-1}|_{C^+})\oplus H^0(C^-,\mathcal{L}^{-1}|_{C^-}) \ |\ s^+(q+)=0, \\
&s^-(q^-)=0, s^+(p^+)=s^-(p^-)\}.
\end{align*}
Without loss of generality we can assume that $\mathcal{L}|_{C^-}\simeq\OO(-1)$.
Then the projection to $s^+$ induces an isomorphism of $V$ with $H^0(C^+,\mathcal{L}^{-1}(-q^+))$.
Now the space $\Hom(J_q,\OO_q)$ has a basis $(e^+,e^-)$ corresponding to the decomposition of the completion
of $J_q$ at $q$ into $xk[[x]] \oplus yk[[y]]$, so that the pairing of $V$ with $e^+$ (resp., $e^-$) is given by
$s^+\mod \mg_{q^+}^2$ (resp.,  $s^-\mod \mg_{q^-}^2$).  
Note that since $s^-(q^-)=0$, the functional $s^-\mod\mg_{q^-}^2$ on $V$ is proportional (with nonzero scalar)
to $s^-(p^-)=s^+(p^+)$. Thus, the two linear maps on $V$ given by the pairing with 
$e^+$ and $e^-$ can be identified with the maps 
$$s^+\mapsto s^+\mod \mg_{q^+}^2, \ \ \ s^+\mapsto s^+(p^+).$$
Since these two maps form a basis of the dual of $V$, the assertion follows.

In the case $n=1$, the assertion similarly reduces to the fact that 
$$s\mapsto s\mod \mg_{\infty}^2, \ \ \ s\mapsto s\mod \mg_0^2$$
form a basis of the dual of $H^0(\P^1, \OO(3)(-0-\infty))$.

Now let us prove the analogous statement on the A side. 
Recall that each $1$-spherical object $S$ gives rise to the corresponding dual twist functor $T'_S$, such that
there is an exact triangle
$$T'_S(F)\to F\to \Hom(F,S)^\vee\otimes S\to\ldots$$
In the case when $S=\OO_{p_i}$ we have $T'_S(F)\simeq F(-p_i)$. 
Hence, $\mathcal{L}=\OO_C(-2p_{i_0}-p_{i_0+1})$ can be obtained from $\OO_C$ by applying three reflection functors
of this type. 

On the Fukaya category side, we have the (negative) Dehn twist exact triangle:
\begin{equation} \label{eq:exactness2}
\xymatrix{
    \tau_{K}^{-1}(L)\ar[r]  &  L     \ar[d]^-{\mathrm{ev}^{\vee}}  \\
                    & \mathit{HF}^* (L,K)^{\vee} \otimes K 
 \ar[ul]^{[1]}
    } 
\end{equation}
where $\tau_{K}^{-1}(L)$ is the exact Lagrangian in $\T_0$ which is
            obtained by a left-handed Dehn twist of $L$ around $K$ (equipped with the induced
            orientation, spin structure and grading).

This implies that the corresponding object $[\LL]$ of the Fukaya category is
given by $\tau_{L_{i_0}}^{-2} \tau_{L_{i_0 +1}}^{-1}(L_0)$ for $n\geq 2$ (resp., by $\tau_{L_1}^{-3}(L_0)$ for $n=1$).

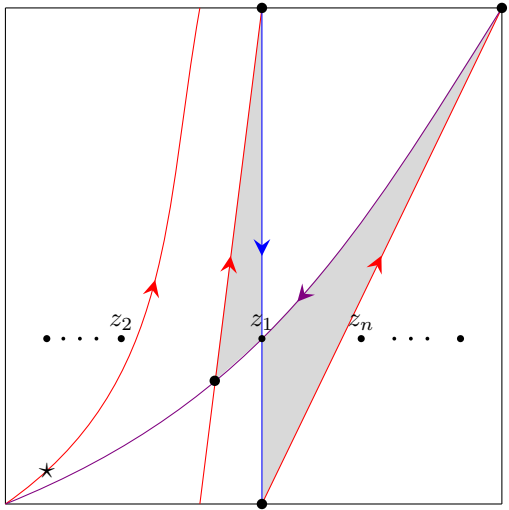
\begin{figure}[htb!]
\centering
\begin{tikzpicture} [scale=1.1]
      \tikzset{->-/.style={decoration={ markings,
        mark=at position #1 with {\arrow[scale=2,>=stealth]{>}}},postaction={decorate}}}
         \draw (0,0) -- (6,0);
         \draw (6,6) -- (0,6);
         \draw (6,0) -- (6,6);
         \draw (0,0) -- (0,6);
      
      \begin{scope} 
          \clip (2.53,1.49) -- (3.1,6) -- (3.1,2) ;
           \clip[preaction={draw,fill=gray!30}] (0,0) rectangle (8,8);
   \end{scope}
  \begin{scope}        
      \clip (3.1,2) to[in=-123,out=46] (6,6) -- (3.1,0) ;
           \clip[preaction={draw,fill=gray!30}] (0,0) rectangle (8,8);
 \end{scope}

         \draw [red, ->-=.5] (0,0) to[in=-100,out=35] (2.35,6);

         \draw[red, ->-=.5]  (2.35,0) -- (3.1,6);
         \draw[red, ->-=.5]  (3.1,0) -- (6,6);
         \draw[blue, ->-=.5] (3.1,6) -- (3.1,0);
         \draw[violet, ->-=.5] (6,6) to[in=22,out=-122] (0,0); 

           \node at (0.5,0.4) {$\star$};

\draw[black, thick, fill=black] (0.7,2) circle(.01);
\draw[black, thick, fill=black] (0.9,2) circle(.01);
\draw[black, thick, fill=black] (1.1,2) circle(.01);
\draw[black, thick, fill=black] (4.7,2) circle(.01);
\draw[black, thick, fill=black] (4.9,2) circle(.01);
\draw[black, thick, fill=black] (5.1,2) circle(.01);

\draw[thick, fill=black] (0.5,2) circle(.03);
\draw[thick, fill=black] (1.4,2) circle(.03);
\draw[thick, fill=black] (3.1,2) circle(.03);
\draw[thick, fill=black] (4.3,2) circle(.03);

\draw[thick, fill=black] (5.5,2) circle(.03);

\node at (1.4,2.2) {\footnotesize $z_2$};
\node at (3.1,2.2) {\footnotesize $z_1$};
\node at (4.3,2.2) {\footnotesize $z_n$};

\draw[thick, fill=black] (6,6) circle(.05);
\draw[thick, fill=black] (3.1,6) circle(.05);
\draw[thick, fill=black] (3.1,0) circle(.05);
\draw[thick, fill=black] (2.53,1.49) circle(.05);

\end{tikzpicture}
\caption{The curves corresponding to $[\mathcal{L}]$ (red), $[J_q[1]]$ (violet) and $\hat{L}_0$
(blue)} 
\label{fig7}
\end{figure}

In Figure \ref{fig7}, assuming that $n\geq 2$,
we drew the isotopy classes of the curves corresponding to $[\mathcal{L}]$ (in red),
$\hat{L}_0$ (in blue) and $[J_q[1]]$ (in violet). For
$n=1$ the same picture can be used except that only the marked point labelled $z_1$ should be left.

Now, the non-degeneracy of the pairing (\ref{line-bundle-pairing-A}) is equivalent to the non-degeneracy of the pairing
\begin{equation}
    \Hom^{1}([J_q][1],[\OO_q])\otimes \Hom^{-1}([\mathcal{L}],[J_q][1])\to \Hom^0([\mathcal{L}],[\OO_q])\simeq k.
\end{equation}
This can be computed via the Floer $\mathfrak{m}_2$-products given by triangles with boundary on
$([\mathcal{L}], [J_q][1], [\OO_q])$ as in Figure \ref{fig7}. We see that there are precisely two
triangles (shaded in Figure \ref{fig7}) contributing to this product. From this we conclude that
the pairing (\ref{line-bundle-pairing-A}) is perfect, as required.  \ed

\subsection{Equivalence of the wrapped Fukaya category with $D^b(\Coh G_{n,R})$}

Recall that the objects $\hat{F}_0,\ldots,\hat{F}_n$ in $D^b(\Coh G_{n,R})$ were defined by \eqref{Fhat-collection-eq}.

\begin{lem}\label{Y-gen-lem} 
One has $Y(\hat{L}_i)\simeq \hat{F}_i$, $i=0,\ldots,n$.
\end{lem}

\Pf . For $i=1,\ldots,n$ this follows from Lemma \ref{char-dual-to-point-lem}. Furthermore,
Lemma \ref{char-nodes-lem} implies that $Y(\hat{L}_0)$ is isomorphic to $\OO_{q_i}$ for one of the nodes $q_i$.
By cyclic symmetry, it is enough to distinguish $\OO_{q_n}$ from the other nodes.
Now we use the property that $\Hom^1(\hat{F}_i,\OO_{q_n})=0$ for $i\neq 1,n$, while for $i=1,n$, this space is nonzero.
\ed

\begin{thm}\label{wrapped-equivalence-thm} 
For a regular ring $R$
there exists an exact equivalence $D^\pi(\WW(\T_0)\ot R)\to D^b(\Coh G_{n,R})$, extending the equivalence of
Corollary \ref{Fuk-Perf-cor}, and sending $\hat{L}_i$ to $\hat{F}_i$ for $i=0,\ldots,n$.
If in addition $\hat{K}_0(R)=0$ (say, $R=\Z$ or $R$ is a field), then $D^\pi(\WW(\T_0)\ot R)=D^b(\WW(\T_0)\ot R)$.
\end{thm}

\Pf . Set $C=G_{n,R}$.
Theorem \ref{A-fully-faithful-thm} together with Lemma \ref{another-generators-W-lem} imply that
the Yoneda functor induces an equivalence of $D^\pi(\WW(\T_0)\ot R)$ with the subcategory of $D(\Qcoh C)$, split-generated by
$Y(\hat{L}_i)$, $i=0,\ldots,n$. But by Lemma \ref{Y-gen-lem} together with Proposition \ref{Db-Coh-gen-prop},
this subcategory is equivalent to $D^b(\Coh C)$. The last assertion follows from Lemma \ref{ncgener} and Corollary \ref{Db-Coh-gen-cor}.
\ed


\section{Koszul duality}
\label{koszul} 

\subsection{Koszul duality for $A_\infty$-algebras}

In this section we work over a field $k$.

We consider the Koszul duality picture for $A_\infty$-algebras, similar to the one for the dg-categories, considered
by Keller in \cite[Sec.\ 10]{Keller-derived-dg}.
It involves a slightly more general notion of augmented $A_\infty$-algebras than the one considered in 
\cite[Sec.\ 3.5]{Keller-ainf} and \cite{LPWZ}.

\begin{defi} (i) Let $K=\bigoplus_{i=0}^n k$, and let $A$ be an $A_\infty$-algebra over $K$. 
Let $S_i$, $i=0,\ldots,n$, be the simple $K$-modules (so that $S_i$ corresponds to the $i$th summand in $K$).
A {\it left augmentation} on $A$ is a collection $(\wt{S}_i)$ of left $A_\infty$-modules over $A$ such
that $H^*\wt{S}_i\simeq S_i$ as $K$-modules. Thus, we get an $A_\infty$-module
$\wt{K}:=\bigoplus_{i=0}^n \wt{S}_i$
over $A$. Similarly, we define a {\it right augmentation} using right $A_\infty$-modules.

\noindent
(ii) Let $A$ be a right augmented $A_\infty$-algebra over $K$. The {\it Koszul dual of} $A$ is the left augmented dg-algebra
$$E(A):=R\Hom_A(\wt{K},\wt{K}),$$
where the augmentation is given by the natural $E(A)$-module structure on $\wt{K}$. Similarly,
the Koszul dual of a left augmented dg-algebra $A$ is the right augmented dg-algebra $E(A)=R\Hom_A(\wt{K},\wt{K})^{op}$.
\end{defi}

Note that if $\wt{K}$ is a left augmentation of $A$ then $R\Hom_K(\wt{K},K)$ is a right augmentation of $A$ and vice versa.

We view augmented $A_\infty$-algebras up to a natural equivalence,
extending the $A_\infty$-equivalence of $A_\infty$-algebras over $K$. 
The operation of passing to the Koszul dual is well-defined on equivalence classes.

\begin{rems} 
1. For any (left or right) augmented $A_\infty$-algebra $A$ there is an $A_\infty$-morphism $A\to E(E(A))$. However,
in general it is not a quasi-isomorphism, and the double Koszul dual $E(E(A))$ should be viewed as some kind 
of completion of $A$, see \cite{DG}, \cite{Efimov}.

\noindent
2. In \cite[Sec.\ 3.5]{Keller-ainf} and \cite{LPWZ} the authors 
work with a stronger notion of augmentation: they require the existence of
a surjection $A\to K$ in the category of $A_\infty$-algebras. In the main example of interest for us, considered below, 
such a surjection may not exist.
\end{rems}


Recall that we have the $A_\infty$-algebra $\AA_0$ associated with the generators $(L_i)$ of the exact Fukaya category
$\FF(\T_0)$ and the $A_\infty$-algebra $\BB$ associated with the generators $(\hat{L}_i)$ of the wrapped Fukaya category
$\WW(\T_0)$ (see \eqref{A0-algebra-eq}, \eqref{B-algebra-eq}). 
Since we now work over a field $k$, we consider the $A_\infty$-algebras
$$\AA_{0,k}:=\AA_0\ot k, \ \ \BB_k:=\BB\ot k.$$

Recall that the cyclic group $\Z/n\Z$ acts on $\FF(\T_0)$ (resp., $\WW(\T_0)$) by means of the transformation $T$
(see \eqref{T-translation-eq}). This action preserves $L_0$ and cyclically permutes the remaining
generators $L_1,\ldots,L_n$, hence, we get
an induced action on the $A_\infty$-algebra $\AA_{0,k}$. Note however, that the $\Z/n\Z$-action on $\WW(\T_0)$ does not
preserve the generator $\bigoplus_{i=0}^n \hat{L}_i$ (the object $\hat{L}_0$ is mapped to similar objects associated with other
punctures).

Using Lemmas \ref{char-nodes-lem} and \ref{char-dual-to-point-lem}, as well as the equivalence of Corollary \ref{Fuk-Perf-cor},
we derive the following fact about $A_\infty$-modules over $\AA_{0,k}$.

\begin{prop}\label{aug-prop} 
There is a unique left (resp., right) augmentation on $\AA_{0,k}$, up to the twist
by the $\Z/n\Z$-action on $\AA_{0,k}$, such that $\dim\Ext^1_{\AA_{0,k}}(\wt{S}_i,\wt{S}_i)>1$ for every $i$
(equivalently, $\Ext^*_{\AA_{0,k}}(\wt{S}_i,\wt{S}_i)$ is infinite-dimensional). \ed
\end{prop}

We are going to check that $\BB_k$ and $\AA_{0,k}$ are Koszul dual of each other, with respect to some augmentations,
such that the left augmentation on $\AA_{0,k}$ is the one given in Proposition \ref{aug-prop}.
Note that $n$ choices of such an augmentation correspond to $n$ different choices of defining the object
$\hat{L}_0$ in our generating set for $\WW(\T_0)$.

\subsection{The proof of Koszul duality on the $B$-side}

Let us set $C=G_{n,k}$, and consider the objects
$$F:=\bigoplus_{i=0}^n F_i, \ \ \hat{F}:=\bigoplus_{i=0}^n \hat{F}_i$$
in $D^b(\Coh C)$ (see Sec.\ \ref{gen-B-sec}).
Due to equivalences of Corollary \ref{Fuk-Perf-cor} and Theorem \ref{wrapped-equivalence-thm},  we have
$$\AA_{0,k}:=R\Hom_{\Perf(C)}(F,F), \ \ \BB_k=R\Hom_{D^b(\Coh C)}(\hat{F},\hat{F}).$$
Here we apply the homological perturbation to get the corresponding $A_\infty$-algebras over $K$.

Let us define the left augmentation of $\AA_{0,k}$ and the right augmentation of $\BB_k$ using 
$$\wt{K}:=R\Hom_{D^b(\Coh C)}(\hat{F},F),$$
which has a natural $\AA_{0,k}-\BB_k$-bimodule structure.
More precisely, we use the corresponding
$A_\infty$-bimodule obtained by the homological perturbation.

\begin{prop}\label{Koszul-B-side-prop} 
One has $A_\infty$-equivalences
$$\AA_{0,k}\simeq R\Hom_{\BB_k}(\wt{K},\wt{K}), \ \ \BB_k\simeq R\Hom_{\AA_{0,k}}(\wt{K},\wt{K})^{op}.$$
Thus, $\AA_{0,k}$ and $\BB_k$ are the Koszul duals of each other.
\end{prop}

\Pf .
Since $\hat{F}$ generates $D^b(\Coh C)$, the functor 
$$R\Hom_{D^b(\Coh C)}(\hat{F},?): D^b(\Coh C)\to D(\mod \text{-}\BB_k)$$
is fully faithful, so we get an equivalence 
$$\AA_{0,k}\simeq R\Hom_{D^b(\Coh C)}(F,F)\simeq R\Hom_{\BB_k}(\wt{K},\wt{K}).$$
On the other hand, we claim that the functor
$$R\Hom_{D^b(\Coh C)})(?,F): D^b(\Coh C)^{op}\to D(\AA_{0,k}\text{-}\mod),$$
sending $\hat{F}$ to $\wt{K}$, is also fully faithful. Indeed, this follows easily from Serre duality and the fact that
$R\Hom_{D^b(\Coh C)})(F,?)$ is fully faithful (with the essential
image consisting of right $\AA_{0,k}$-modules with finite-dimensional cohomology). Thus, we get
$$\BB_k\simeq R\Hom_{D^b(\Coh C)}(\hat{F},\hat{F})\simeq R\Hom_{\AA_{0,k}}(\wt{K},\wt{K})^{op}.$$
\ed

\begin{rem}\label{Koszul-rem} 
A Koszul duality result between $\mathcal{F}(M)$ and $\mathcal{W}(M)$ was proven in
    \cite{EL} in the case when $M$ is a plumbing of the cotangent bundles $T^*Q$'s according to a
    plumbing tree, where $Q=S^2$. The Koszul duality result above can be seen as an analogue of this for $Q=S^1$.
    It is interesting to note that the result in \cite{EL} was inspired by the
    classical Koszul duality result between the dg-algebras $C^*(Q)$ and  $C_{-*}(\Omega Q)$ for $Q$ a \emph{simply-connected} manifold. The relevance of these two dg-algebras comes
    from the fact that they are quasi-isomorphic to the endomorphism algebras of generators of
    $\mathcal{F}(T^*Q)$ and $\mathcal{W}(T^*Q)$ respectively (see \cite{EL} for details). In the $2$-dimensional case 
    we see that even though the
    corresponding Koszul duality result fails for $T^*S^1$ as $S^1$ is not simply-connected, Koszul duality between
    $\mathcal{F}(M)$ and $\mathcal{W}(M)$ holds for $M=\T_0$, a 
    plumbing of $T^*S^1$'s according to the star-shaped tree as in Figure \ref{fig0}. 
\end{rem}

\begin{rem} In Lemma \ref{ncgener} we prove that the objects
$\hat{L}_0, \hat{L}_1, \ldots, \hat{L}_n$ generate the wrapped Fukaya category 
$\mathcal{W}(\T_0)$. It is likely that this can also be checked via the generation criterion from
\cite{abouzgen} (see Remark \ref{pitfall}). This is the statement that the natural
map
\[ HH_{*-1}(\mathcal{W}(\T_0)) \to SH^*(\T_0) \]
relating the Hochschild homology of the wrapped Fukaya category to the symplectic cohomology of
$\T_0$ hits the identity element. When this generation criterion is
satisfied, then it was proven by Ganatra \cite{ganatra} that the natural maps \[
HH_{*-1}(\mathcal{W}(\T_0)) \to SH^*(\T_0) \to HH^*(\mathcal{W}(\T_0)) \] are all isomorphisms. On
the other hand, in the case at hand, it is easy to compute $SH^*(\T_0)$ explicitly via a Morse-Bott
type spectral sequence (see Ex. 3.3 of \cite{seidelbiased}). In fact, we have
\[ SH^0(\T_0) = k\ , \ SH^1(\T_0) = k^{n+1}, \ SH^d(\T_0) =k^n (d \geq 2). \] 
In particular, we deduce that $\text{dim} HH^2(\mathcal{W}(\T_0)) = n$. 

Now, we have seen above that there is a Koszul duality
between the exact Fukaya category $\mathcal{F}(\T_0)$ and the wrapped Fukaya category
$\mathcal{W}(\T_0)$. By the result of Keller \cite{keller} this implies that we
have an isomorphism between $HH^*(\mathcal{F}(\T_0))$ and
$HH^*(\mathcal{W}(\T_0))$. This gives that $\text{dim}HH^2 (\mathcal{F}(\T_0)) = \text{dim}
HH^2(\mathscr{A}_0, \mathscr{A}_0)= n$, so the above calculation matches our calculation of
$HH^2(G_{n,k})$ in Lemma \ref{HH-wheel-lem}. 
\end{rem}

\end{document}